\documentclass{article}
\usepackage{amsthm,amsmath,amssymb,mathrsfs,afterpage,url,bm,graphicx,color}
\RequirePackage[colorlinks,citecolor=blue,urlcolor=blue]{hyperref}
\def\P{{\rm P}} 
\def\E{{\rm E}} 
\def\std{{\rm std}}
\def\hit{{\rm hit}}
\def\dbl{{\rm dbl}}
\def\spl{{\rm spl}}
\allowdisplaybreaks

\title{\textbf{Snackjack: A toy model of blackjack}}
\author{Stewart N. Ethier\footnote{Department of Mathematics, University of Utah, 155 S. 1400 E., Salt Lake City, UT 84112, USA. Email: \url{ethier@math.utah.edu}.} \ and Jiyeon Lee\footnote{Department of Statistics, Yeungnam University, 280 Daehak-Ro, Gyeongsan, Gyeongbuk 38541, South Korea. Email: \url{leejy@yu.ac.kr}.}}
\date{}
\begin{document}
\maketitle

\begin{abstract}
Snackjack is a highly simplified version of blackjack that was proposed by Ethier (2010) and given its name by Epstein (2013).  The eight-card deck comprises two aces, two deuces, and four treys, with aces having value either 1 or 4, and deuces and treys having values 2 and 3, respectively.  The target total is 7 (vs.\ 21 in blackjack), and ace-trey is a natural.  The dealer stands on 6 and 7, including soft totals, and otherwise hits. The player can stand, hit, double, or split, but split pairs receive only one card per paircard (like split aces in blackjack), and there is no insurance.

We analyze the game, both single and multiple deck, deriving basic strategy and one-parameter card-counting systems.  Unlike in blackjack, these derivations can be done by hand, though it may nevertheless be easier and more reliable to use a computer.  More importantly, the simplicity of snackjack allows us to do computations that would be prohibitively time-consuming at blackjack.  We can thereby enhance our understanding of blackjack by thoroughly exploring snackjack.\medskip

\end{abstract}

\section{Introduction}\label{intro}

According to Marzuoli~\cite{M09}, 

\begin{quote}
Toy models in theoretical physics are invented to make simpler the modelling of complex physical systems while preserving at least a few key features of the originals. Sometimes toy models get a life of their own and have the chance of emerging as paradigms.
\end{quote}
For example, the simple coin-tossing games of Parrondo form a toy model of the rather complex flashing Brownian ratchet in statistical physics (see, e.g.,~\cite{EL18}).  Our aim here is to explore a toy model of the game of blackjack, primarily as a way of gaining insight.

One possible toy model of blackjack is the contrived game of \textit{red-and-black} in which one can bet, at even money, that the next card dealt will be red.  This game has been studied by Thorp and Walden~\cite{TW73}, Griffin~\cite{G76}, Ethier and Levin~\cite{EL05}, and others.  Its simplicity allows the card counter to play perfectly, and analysis is straightforward.  However, because the game is vastly simpler than blackjack and rather unlike blackjack, the insights it offers are limited.

Epstein, in the first edition of \textit{The Theory of Gambling and Statistical Logic}~\cite[p.~269]{Ep67}, proposed the game of \textit{grayjack}, a simplified version of blackjack, ``offering an insight into the structure of the conventional game.''  Grayjack uses a 13-card deck comprising one ace, two twos, two threes, two fours, two fives, and four sixes, with aces having value 1 or 7, and the other cards having their nominal values.  The target total is 13, and ace-six is a natural.  The dealer stands on 11, 12, and 13, including soft totals, and otherwise hits.  The player can stand, hit, double, or split, just as in blackjack, but there is no resplitting.  The problem with grayjack, as a toy model of blackjack, is that its analysis is only marginally simpler than that of blackjack itself.  Both require elaborate computer programs.  Fifty years after grayjack was introduced, its basic strategy was still unpublished and perhaps even unknown (but see Appendix~A).

Ethier, in \textit{The Doctrine of Chances}~\cite[Problem 21.19]{Et10}, proposed an even simpler toy model of blackjack, which was renamed \textit{snackjack} by Epstein in the latest edition of \textit{TGSL}~\cite[p.\ 291]{Ep13}.  The eight-card deck comprises two aces, two deuces, and four treys, with aces having value either 1 or 4, and deuces and treys having values 2 and 3, respectively.  The target total is 7, and ace-trey is a natural.  The dealer stands on 6 and 7, including soft totals, and otherwise hits. The player can stand, hit, double, or split, but split pairs receive only one card per paircard (like split aces in blackjack), and there is no  insurance.  Unlike blackjack or grayjack, snackjack can be analyzed by hand, though it may nevertheless be easier and more reliable to use a computer.

\begin{table}[htb]
\caption{\label{3games}Salient features of blackjack, grayjack, and snackjack.}
\catcode`@=\active \def@{\hphantom{0}}
\catcode`#=\active \def#{\hphantom{$\,^1$}}
\tabcolsep=.2cm
\begin{center}
\begin{tabular}{cccc}
\hline
\noalign{\smallskip}
        & blackjack & grayjack & snackjack \\
\noalign{\smallskip}
\hline
\noalign{\smallskip}
deck size        & 52                     & 13                & 8                 \\                
\noalign{\medskip}
deck             & four aces;             & one ace;          & two aces;         \\
composition      & four each of 2--9;     & two each of 2--5; & @@two deuces;@@   \\
                 & 16 tens                & four sixes           & four treys         \\
\noalign{\medskip}
ace value        & 1 or 11                & 1 or 7            & 1 or 4            \\
\noalign{\medskip}
target total     & 21                     & 13                & 7                 \\
\noalign{\medskip}
natural          & ace-ten                & ace-six           & ace-trey          \\
\noalign{\medskip}
dealer stands    & 17--21                 & 11--13            & 6--7              \\
                 & @(incl.~soft)@         & @(incl.~soft)@    & @(incl.~soft)@    \\
\noalign{\smallskip}
\hline
\end{tabular}
\end{center}
\end{table}

Our aim here is to thoroughly analyze the game of snackjack, both single and multiple deck.  We do not propose snackjack as a new casino game; rather, we believe that it offers insight into the more complex game of blackjack.  Specifically, the simplicity of snackjack allows us to do computations that would be prohibitively time-consuming at blackjack.  In addition, snackjack has pedagogical value:  The very complex theory of blackjack becomes a bit more transparent when viewed through the lens of this simple toy model.

In Section~\ref{rules} we give a self-contained description of the rules of snackjack.  Section~\ref{comparisons} tries to justify our claim that snackjack, but not blackjack or grayjack, can be analyzed by hand.  Section~\ref{BS methodology} describes our methodology for deriving basic strategy, which is consistent with what was used for single-deck blackjack in~\cite[Section 21.2]{Et10}.  It can also be adapted to grayjack.  In Section~\ref{BS results} we apply this methodology to snackjack, both single and multiple deck.  Section~\ref{FTCC} explores some of the consequences for snackjack of the fundamental theorem of card counting~\cite{EL05,TW73}.  In Section~\ref{bet variation} we investigate card counting at snackjack and its application to bet variation.  In all discussions of card counting, we emphasize the 39-deck (312-card) game.  In Section~\ref{strategy variation} we explore the potential for gain by varying basic strategy.  Finally, Section~\ref{conclusions} summarizes what snackjack tells us about blackjack.

\section{Detailed rules of snackjack}\label{rules}

Snackjack is played with a single eight-card deck comprising two aces, two deuces, and four treys, or with multiple such decks mixed together.  Aces have value either 1 or 4, and deuces and treys have values 2 and 3, respectively.  Suits do not play a role.  A hand comprising two or more cards has value equal to the total of the values of the cards.  The total is called \textit{soft} if the hand contains an ace valued as 4, otherwise it is called \textit{hard}. 

Each player competes against the dealer.  (In the single-deck game there can be only one player.)  After making a bet, each player receives two cards, usually face down (but it does not actually matter), and the dealer receives two cards, one face down (the \textit{downcard} or \textit{hole card}) and one face up (the \textit{upcard}).  If the player has a two-card total of 7 (a \textit{natural}) and the dealer does not, the player wins and is paid 3 to 2.  If the dealer has a natural and the player does not, the player loses his bet.  If both player and dealer have naturals, a push is declared.  If the dealer's upcard is an ace or a trey, he checks his downcard to determine whether he has a natural before proceeding.  There is no insurance bet in snackjack.

If the dealer and at least one player fail to have naturals, play continues.  Starting with the player to the dealer's left and moving clockwise, each player completes his hand as follows.   He must decide whether to \textit{stand} (take no more cards) or to \textit{hit} (take an additional card).  If he chooses the latter and his new total does not exceed 7, he must make the same decision again and continue to do so until he either stands or \textit{busts} (his total exceeds 7).  If the player busts, his bet is lost, even if the dealer subsequently busts.  The player has one or two other options after seeing his first two cards.  He may \textit{double down}, that is, double his bet and take one, and only one, additional card.  If he has two cards of the same value, he may \textit{split} his pair, that is, make an additional bet equal to his initial one and play two hands, with each of his first two cards being the initial card for one of the two hands and each of his two bets applying to one of the two hands.  Each card of the split pair receives one and only one card.  (This last rule was mistakenly omitted from~\cite[Problem 21.19]{Et10}.  It is needed to avoid the possibility of running out of cards.\footnote{While it is not possible to run out of cards in a one-player vs.\ dealer single-deck game, it \textit{is} possible that all eight cards are needed.  For example, consider a pair of treys against a dealer deuce.  Player splits, getting a trey on each trey.  Dealer's downcard is also a deuce, and he draws an ace, then another ace, exhausting the deck for a total of six and a double push.  This is just one of several scenarios in which the full deck is needed.})  A two-card 7 after a split is not regarded as a natural and is therefore not entitled to a 3-to-2 payoff; in addition, it pushes a dealer 7 comprising three or more cards.  

As we have already assumed, the dealer checks for a natural when his upcard is an ace or a trey.  This is sometimes stated by saying that an untied dealer natural wins original bets only---additional bets due to doubling or splitting, if they could be made, would be pushed.

After each player has stood, busted, doubled down, or split pairs, the dealer acts according to a set of mandatory rules.  The dealer stands on hands of 6 and 7, including soft totals, and otherwise hits.

If the dealer busts, all remaining players are paid even money.  If the dealer stands, his total is compared to that of each remaining player.  If the player's total (which does not exceed 7) exceeds the dealer's total, the player is paid even money.  If the dealer's total (which does not exceed 7) exceeds the player's total, the player loses his bet.  If the player's total and the dealer's total are equal (and do not exceed 7), the hand is declared a push.

\section{Blackjack vs.\ grayjack vs.\ snackjack}\label{comparisons}

We pause to compare blackjack, grayjack, and snackjack.  First, it is important to clarify the specific blackjack rules that we assume.  The assumed set of rules was at one time standard on the Las Vegas Strip, so we consider it the benchmark against which other sets of rules can be measured.  In the notation of the blackjack literature, we assume S17 (dealer stands on soft 17), DOA (double down any first two cards), NDAS (no double after splits), SPA1 (split aces once, receiving only one card per ace), SPL3 (split non-ace pairs up to three times [up to four hands]), 3:2 (untied player natural pays 3 to 2), OBO (untied dealer natural wins original bets only), and NS (no surrender).

As for grayjack, we assume similarly S11, DOA, NDAS, SPA1 (in the case of multiple decks), SPL1 (no resplitting), 3:2, OBO, and NS.

Finally, snackjack rules can be summarized by S6, DOA, NDAS, SPP1 (split pairs once, receiving only one card per paircard), 3:2, OBO, NS, and NI (no insurance).

The more restrictive pair-splitting rules in grayjack and snackjack ensure against running out of cards in the single-deck (one player vs.\ dealer) games.\footnote{To see that SPL2 could result in an incomplete game in single-deck grayjack, consider a player $6,6$ vs.~a dealer ace.  Player splits, draws another 6, and resplits.  First hand is $6,2,2,3$, second hand is $6,4,5$, and third hand is $6,4,5$.  Dealer's hand is $\text{ace},3,6$, exhausting the deck before completing the hand.  Note that player violated basic strategy only when splitting 6s.}  We maintain these rules in the multiple-deck games even if running out of cards is no longer an issue.

Table~\ref{comparisons-1deck} (single deck) and Table~\ref{comparisons-multideck} (multiple deck) compare various statistics for blackjack, grayjack, and snackjack.  The aim is to justify our claim that blackjack and grayjack analyses require a computer, whereas a comparable analysis of snackjack, while tedious, does not.

\begin{table}[htb]
\caption{\label{comparisons-1deck}Single-deck comparisons of blackjack and its toy models.}
\catcode`@=\active \def@{\hphantom{0}}
\catcode`#=\active \def#{\hphantom{$\,^1$}}
\tabcolsep=.1cm
\begin{center}
\begin{tabular}{cccc}
\hline
\noalign{\smallskip}
statistic & blackjack & grayjack & snackjack \\
\noalign{\smallskip}
\hline
\noalign{\smallskip}
number of cards & 52 & 13 & 8 \\
\noalign{\medskip}
number of unordered & \\
two-card player hands & 55 & 20 & 6 \\
\noalign{\medskip}
number of unordered & \\
unbusted player hands & #$2{,}008\,^1$ & 87 & 14 \\
\noalign{\medskip}
number of comp.-dep.\ basic  & \\
strategy decision points &  #$19{,}620\,^2$ & 430   & 32 \\
\noalign{\medskip}
number of ordered & \\
dealer drawing sequences &  #$48{,}532\,^3$ & 498   & 17 \\
\noalign{\medskip}
number of unordered & \\
dealer drawing sequences &  #$2{,}741\,^4$ & 93   & 11 \\
\noalign{\medskip}
mimic-the-dealer & \\
strategy expectation     & #$-0.0568456\,^5$ & $-0.0584311\,^6$  & $+0.0952381\,^6$ \\
\noalign{\medskip}
composition-dep.\ basic & \\
strategy expectation     & #$+0.000412516\,^7$ & $+0.0218749\,^6$ & $+0.192857\,^8$ \\
\noalign{\smallskip}
\hline
\noalign{\smallskip}
\multicolumn{4}{l}{$^1$\,\cite[p.~275]{Ep13} or~\cite[p.~655]{Et10}. $^2$\,\cite[p.~655]{Et10}.  $^3$\,\cite[p.~275]{Ep13},~\cite[pp.~9--11, 649]{Et10}, or} \\
\multicolumn{4}{l}{\cite[p.~158]{G99}.  $^4$\,\cite[pp.~646, 648]{Et10}. $^5$\,\cite[p.~647]{Et10}.  $^6$\,Computed by the authors.}\\
\multicolumn{4}{l}{$^7$\,\cite[p.~661]{Et10}.  $^8$\,Section~\ref{BS results}.}\\
\end{tabular}
\end{center}
\end{table}

\begin{table}[htb]
\caption{\label{comparisons-multideck}Multiple-deck comparisons of blackjack and its toy models.}
\catcode`@=\active \def@{\hphantom{0}}
\catcode`#=\active \def#{\hphantom{$\,^1$}}
\tabcolsep=.1cm
\begin{center}
\begin{tabular}{cccc}
\hline
\noalign{\smallskip}
statistic & six-deck    & 24-deck  & 39-deck \\
         & blackjack & grayjack & snackjack \\
\noalign{\smallskip}
\hline
\noalign{\smallskip}
number of cards & 312 & 312 & 312 \\
\noalign{\medskip}
number of unordered & \\
two-card player hands & 55 & 21 & 6 \\
\noalign{\medskip}
number of unordered & \\
unbusted player hands & #$3{,}072\,^1$ & 291 & 27 \\
\noalign{\medskip}
number of comp.-dep.\ basic & \\
strategy decision points & $30{,}720$  & $1{,}746$   & 81 \\
\noalign{\medskip}
number of ordered & \\
dealer drawing sequences & #$54{,}433\,^2$  &  $1{,}121$   & 21 \\
\noalign{\medskip}
number of unordered & \\
dealer drawing sequences &  $3{,}357$ &  257  & 15 \\
\noalign{\medskip}
mimic-the-dealer & \\
strategy expectation     & $-0.0567565\,^3$ &  $-0.0628381\,^3$  &  $+0.0720903\,^3$  \\
\noalign{\medskip}
composition-dep.\ basic & \\
strategy expectation & $-0.00544565\,^4$ & $-0.0189084\,^3$  & $+0.139309\,^5$ \\
\noalign{\smallskip}
\hline
\noalign{\smallskip}
\multicolumn{4}{l}{$^1$\,\cite[p.~172]{G99}. $^2$\,\cite[p.~158]{G99}. $^3$\,Computed by the authors.  $^4$\,Computed by}\\
 \multicolumn{4}{l}{Marc Estafanous.  $^5$\,Table~\ref{EV(d)}.}
\end{tabular}
\end{center} 
\end{table}

Let us briefly explain these statistics.  The number of unordered two-card player hands in blackjack is well known to be $\binom{10}{2}+\binom{10}{1}=55$.  Similar calculations apply to grayjack and snackjack, except that a pair of aces is impossible in single-deck grayjack.  In single-deck blackjack the number of unordered unbusted player hands (of any size) is $2{,}008$.  This is simply the sum over $2\le n\le21$ of the number of partitions of the integer $n$ into two or more parts with no part greater than 10 and no part having multiplicity greater than 4.  It is the number of hands that must be analyzed for composition-dependent basic strategy.  

The number of composition-dependent basic strategy decision points is the number of unordered unbusted player hands multiplied by the number of possible dealer upcards, excluding those cases that require more cards than are available.  For example, in single-deck snackjack, $1,2,2$ (an ace and two deuces) vs.\ 2 is ruled out because it requires three deuces, more than are in the deck.  Epstein~\cite[p.~291]{Ep13} reported 33 decision points because he included a spurious one, namely $1,2,2$ vs.\ 1.  Indeed, all basic strategy expectations are conditioned on the dealer not having a natural, but in this case, only treys remain, so the dealer's downcard must be a trey.  In effect, we are conditioning on an event of probability 0, so this case must be excluded.  

The number of ordered dealer drawing sequences is readily computed by direct enumeration.  For example, sequence number $24{,}896$ (in reverse lexicographical order) of the $48{,}532$ such sequences in single-deck blackjack is $3,2,2,2,2,1,1,1,1,3$.  Without regard to order, this sequence would be listed as $(4,4,2,0,0,0,0,0,0,0)$ (i.e., 4 aces, 4 twos, and 2 threes), with total 18 and multiplicity 15 (i.e., 15 permutations of $3,2,2,2,2,1,1,1,1,3$ appear in the ordered list).  The ordered dealer drawing sequences are used to compute conditional expectations when standing.  The unordered dealer drawing sequences are used to evaluate the player's expectation under the mimic-the-dealer strategy~\cite[p.~647]{Et10}, which depends on $\P(\text{both player and dealer bust})$.  (The double bust is hypothetical; it assumes that the dealer deals out his hand even after the player has busted.)  In both blackjack and grayjack, the dealer advantage of acting last (because the dealer wins double busts) dominates the player advantage of a 3-to-2 payoff for an untied natural.  In snackjack, the opposite is true because double busts are rare (probability 2/105 in single deck) and winning player naturals are quite common (probability 8/35 in single deck).  

Finally, the player's expectation under composition-dependent basic strategy has been computed in blackjack, grayjack, and snackjack.  Of course it is substantially larger than that for the mimic-the-dealer strategy.  In single-deck blackjack this expectation is positive, barely, which may explain why the assumed set of rules is obsolete.  In single-deck grayjack it is about +2.19\%, well below Epstein's~\cite[p.~291]{Ep13} estimate of +7.5\%, but the positive expectation nevertheless ``mitigates its suitability as a casino game,'' as Epstein noted.  Perhaps 24-deck grayjack ($-1.89$\%) would be viable as a casino game, however; see Appendix~A.  Snackjack (+19.3\% for single deck, +13.9\% for 39 decks) would certainly not be.  But that is not our concern.  Instead, we want to gain insight into blackjack by studying snackjack. 

Certainly, it would be possible to make a simple rules change that would give the advantage to the house and make snackjack a potential casino game.  There are probably many ways to do this, but an especially simple approach would be to impose the rule, ``A player natural pays even money [instead of 3 to 2], with the exception that it loses to a dealer natural [instead of pushing].''  The result is $+3.10$\% for single deck, $-0.0959$\% for double deck, $-0.713$\% for triple deck, and $-1.73$\% for 39 decks.  We do not pursue this, however.  Instead, when we want the game to be slightly disadvantageous for the purpose of our card-counting analysis, we impose a suitable commission, specifically 1/7 of the amount initially bet in the 39-deck game, resulting in a net expectation of $-0.355$\%.

\section{Snackjack basic strategy methodology}\label{BS methodology}

The term ``basic strategy'' has several interpretations.  See Schlesinger~\cite[Appendix A]{S18} for a thorough discussion of the issues.  We will interpret it as composition-dependent basic strategy, since total-dependent or partially total-dependent basic strategy is an unnecessary compromise in this simple game.  Because of our restrictive rules on splitting, we need not concern ourselves with which cards are used for decisions about split hands.  We follow the approach originated for blackjack by Manson et al.~\cite{MBG75} and used by Griffin~\cite[p.~172]{G99} and Ethier~\cite[Section 21.2]{Et10}, and we use the notation of the latter source.  A completely different approach was taken by Werthamer~\cite[Section 7.2.1]{We18}, who wrote (p.~74), ``\dots\ no [previous] study describes its methodology in detail \dots,'' regrettably overlooking~\cite{Et10}, which was published eight years earlier.

A similar approach applies to grayjack, but here we must clarify how splits are treated.  Our convention is that the player makes use only of the cards in the hand he is currently playing, and of course the dealer's upcard. 

Returning to snackjack, an arbitrary pack is described by ${\bm n}=(n_1,n_2,n_3)$, meaning that it comprises $n_1$ aces, $n_2$ deuces, and $n_3$ treys, with 
\begin{equation*}
|{\bm n}|:=n_1+n_2+n_3
\end{equation*}
being the number of cards.  An unordered player hand is denoted by ${\bm l}=(l_1,l_2,l_3)$ if it comprises $l_1\le n_1$ aces, $l_2\le n_2$ deuces, and $l_3\le n_3$ treys.  The number of cards in the hand is
\begin{equation*}
|{\bm l}|:=l_1+l_2+l_3,
\end{equation*}
and the hand's total is 
\begin{equation*}
T({\bm l}):=\begin{cases}l_1+2l_2+3l_3+3&\text{if $l_1\ge1$ and $l_1+2l_2+3l_3\le4$},\\
\noalign{\smallskip}
l_1+2l_2+3l_3&\text{otherwise},\end{cases}
\end{equation*}
with the two cases corresponding to soft and hard totals.  For the hand to be unbusted, ${\bm l}$ must satisfy $T({\bm l})\le7$.

Let ${\bm X}$ denote the player's hand, let $Y$ denote the player's next card, if any, and let $U$ denote the dealer's upcard, $D$ his downcard, and $S$ his final total.  Finally, let $G_\std$, $G_\hit$, $G_\dbl$, and $G_\spl$ denote the player's profit from standing, hitting, doubling, and splitting, assuming an initial one-unit bet.

Here and in what follows, we often denote an ace not by A but by 1.

We denote the events on which we will condition by
\begin{equation*}
A({\bm l},u):=\begin{cases}\{{\bm X}={\bm l},\; U=1,\; D\ne3\}&\text{if $u=1$},\\
\{{\bm X}={\bm l},\; U=2\}&\text{if $u=2$},\\
\{{\bm X}={\bm l},\; U=3,\; D\ne1\}&\text{if $u=3$},\end{cases}
\end{equation*}
and we define the conditional expectations associated with each player hand, dealer upcard, and strategy:
\begin{align*}
E_\std({\bm l},u)&:=\E[G_\std\mid A({\bm l},u)],\\  
E_\hit({\bm l},u)&:=\E[G_\hit\mid A({\bm l},u)],\\  \
E_\dbl({\bm l},u)&:=\E[G_\dbl\mid A({\bm l},u)]\quad (|{\bm l}|=2), \\ 
E_\spl({\bm l},u)&:=\E[G_\spl\mid A({\bm l},u)]\quad ({\bm l}=2{\bm e}_i,\;i\in\{1,2,3\}),
\end{align*}
where ${\bm e}_1:=(1,0,0)$, ${\bm e}_2:=(0,1,0)$, and ${\bm e}_3:=(0,0,1)$.  Temporarily, we define the maximal stand/hit conditional expectation for each player hand and dealer upcard by
\begin{equation}\label{Emax}
E_{\max}^*({\bm l},u):=\max\{E_\std({\bm l},u),E_\hit({\bm l},u)\}.
\end{equation}

We specify more precisely the set of player hands and dealer upcards we will consider.  We denote the set of all unordered unbusted player hands of two or more cards by
\begin{equation*}
\mathscr{L}:=\{{\bm l}\le{\bm n}: |{\bm l}|\ge2,\; T({\bm l})\le 7\}
\end{equation*}
and the set of all pairs of such hands and dealer upcards by
\begin{equation*}
\mathscr{M}:=\{({\bm l},u)\in \mathscr{L}\times \{1,2,3\}: l_u\le n_u-1\}.
\end{equation*}
The cardinality of $\mathscr{L}$ is the sum over $2\le n\le7$ of the number of partitions of the integer $n$ into two or more parts with no part greater than 3 and 1s having multiplicity at most $n_1$, 2s having multiplicity at most $n_2$, and 3s having multiplicity at most $n_3$.  

The basic relations connecting the conditional expectations defined above  include, for all $({\bm l},u)\in\mathscr{M}$,
\begin{align}\label{Estd}
E_\std({\bm l},u)&=\P(S<T({\bm l}){\rm\ or\ }S>7\mid A({\bm l},u))\nonumber\\
&\qquad\;{}-\P(T({\bm l})<S\le7\mid A({\bm l},u)),\\ \label{Ehit}
E_\hit({\bm l},u)&=\sum_{1\le k\le3:\;({\bm l}+{\bm e}_k,u)\in\mathscr{M}}p(k\mid {\bm l},u)E_{\max}^*({\bm l}+{\bm e}_k,u)\nonumber\\
&\qquad{}+\sum_{1\le k\le3:\;{\bm l}+{\bm e}_k\notin\mathscr{L}}p(k\mid {\bm l},u)(-1),\\  \label{Edbl}
E_\dbl({\bm l},u)&=2\sum_{1\le k\le3:\;({\bm l}+{\bm e}_k,u)\in\mathscr{M}}p(k\mid {\bm l},u)E_\std({\bm l}+{\bm e}_k,u) \nonumber\\
&\qquad{}+2\sum_{1\le k\le3:\;{\bm l}+{\bm e}_k\notin\mathscr{L}}p(k\mid {\bm l},u)(-1)\qquad(|{\bm l}|=2),\\ \label{Espl}
E_\spl(2{\bm e}_i,u)&=2\sum_{1\le k\le3} p(k\mid 2{\bm e}_i,u)E_\std({\bm e}_i+{\bm e}_k,u\mid{\bm e}_i)\qquad (i=1,2,3),
\end{align}
where
\begin{equation*}
p(k\mid{\bm l},u):=\P(Y=k\mid A({\bm l},u)).
\end{equation*}
The probabilities $p(k\mid{\bm l},u)$ are derived from Bayes' law for $u=1$ and $u=3$:
\begin{align}\label{p(k,1)}
p(k\mid{\bm l},1)&=\frac{n_k-l_k-\delta_{1,k}}{|{\bm n}|-|{\bm l}|-1}\bigg(\frac{1-(n_3-l_3-\delta_{3,k})/(|{\bm n}|-|{\bm l}|-2)}{1-(n_3-l_3)/(|{\bm n}|-|{\bm l}|-1)}\bigg),\\ \label{p(k,2)}
p(k\mid{\bm l},2)&=\frac{n_k-l_k-\delta_{2,k}}{|{\bm n}|-|{\bm l}|-1},\\ \label{p(k,3)}
p(k\mid{\bm l},3)&=\frac{n_k-l_k-\delta_{3,k}}{|{\bm n}|-|{\bm l}|-1}\bigg(\frac{1-(n_1-l_1-\delta_{1,k})/(|{\bm n}|-|{\bm l}|-2)}{1-(n_1-l_1)/(|{\bm n}|-|{\bm l}|-1)}\bigg),
\end{align}
where $\delta_{u,k}$ is the Kronecker delta.
Equation \eqref{Espl} comes from~\cite[Eq.~(21.53)]{Et10} and requires a slight extension of our notation.  We define $E_{\std}({\bm l},u\mid {\bm m})$ analogously to $E_{\std}({\bm l},u)$, but with ${\bm m}=(m_1,m_2,m_3)$ indicating that the initial pack is depleted by removing $m_1$ aces, $m_2$ deuces, and $m_3$ treys (in addition to the cards in the player's hand and the dealer's upcard).  Thus, $E_{\std}({\bm l},u)=E_{\std}({\bm l},u\mid {\bm 0})$.

The quantities \eqref{Estd} are computed directly, while those in \eqref{Ehit} are obtained recursively.  They are recursive in the player's hard total
\begin{equation*}
T_{{\rm hard}}({\bm l}):=l_1+2l_2+3l_3.
\end{equation*}
The recursion is initialized with
\begin{equation}\label{initial}
E_\hit({\bm l},u)=-1,\qquad ({\bm l},u)\in\mathscr{M},\; T_{{\rm hard}}({\bm l})=7.
\end{equation}
There is one exception to \eqref{Estd} because an untied player natural is paid 3 to 2:
\begin{equation}\label{natural}
E_\std({\bm e}_1+{\bm e}_3,u)=\frac32,\qquad u=1,2,3.
\end{equation}

We begin by computing $E_\std({\bm l},u)$ for all $({\bm l},u)\in\mathscr{M}$ using \eqref{Estd} (except for  \eqref{natural}).   The number of ordered dealer drawing sequences that must be analyzed for each such ${\bm l}$ is at most $21$.  Then we go back and compute $E_{\max}^*({\bm l},u)$ of \eqref{Emax} for $T_{{\rm hard}}({\bm l})=7,6,5,4,3,2$ (in that order) and all $u$ using \eqref{Estd}, \eqref{Ehit}, \eqref{p(k,1)}--\eqref{p(k,3)}, and \eqref{initial}.  Finally, we compute $E_\dbl({\bm l},u)$ using \eqref{Estd}, \eqref{Edbl}, and \eqref{p(k,1)}--\eqref{p(k,3)}, and $E_\spl({\bm l},u)$ using \eqref{Estd}, \eqref{Espl}, and \eqref{p(k,1)}--\eqref{p(k,3)}.  We can finally evaluate
\begin{equation*}
E_{\max}({\bm l},u):=\begin{cases}\max\{E_\std({\bm l},u),E_\hit({\bm l},u),E_\dbl({\bm l},u),E_\spl({\bm l},u)\}&\text{if ${\bm l}=2{\bm e}_i$},\\
\max\{E_\std({\bm l},u),E_\hit({\bm l},u),E_\dbl({\bm l},u)\}&\text{if ${\bm l}={\bm e}_i+{\bm e}_j$},\\
\max\{E_\std({\bm l},u),E_\hit({\bm l},u)\}&\text{if $|{\bm l}|\ge3$},\end{cases}
\end{equation*}
where $i\in\{1,2,3\}$ in the first line and $i,j\in\{1,2,3\}$ with $i<j$ in the second.

We can also evaluate the player's overall expectation $E$ using the optimal strategy thus derived.  It is simply a matter of conditioning on the player's initial two-card hand and the dealer's upcard.  Now the 18 events $A({\bm l},u)$ for $|{\bm l}|=2$ and $u=1,2,3$ do not partition the sample space, but if we include the 12 events
\begin{align*}
B({\bm e}_i+{\bm e}_j,1)&:=\{{\bm X}={\bm e}_i+{\bm e}_j,\; U=1,\; D=3\}\\
B({\bm e}_i+{\bm e}_j,3)&:=\{{\bm X}={\bm e}_i+{\bm e}_j,\; U=3,\; D=1\}
\end{align*}
as well, where $1\le i\le j\le 3$, then we do have a partition, and conditioning gives the desired result, namely
\begin{align*}
E&=\sum_{u=1}^3\;\mathop{\sum\sum}_{1\le i\le j\le 3}\P(A({\bm e}_i+{\bm e}_j,u))E_{\max}({\bm e}_i+{\bm e}_j,u)\nonumber\\
&\qquad{}+\mathop{\sum\sum}_{1\le i\le j\le 3:\; (i,j)\ne(1,3)}\P(B({\bm e}_i+{\bm e}_j,1))(-1)\nonumber\\
&\qquad{}+\mathop{\sum\sum}_{1\le i\le j\le 3:\; (i,j)\ne(1,3)}\P(B({\bm e}_i+{\bm e}_j,3))(-1) \nonumber\\
&\qquad{}+\P(B({\bm e}_1+{\bm e}_3,1))(0)+\P(B({\bm e}_1+{\bm e}_3,3))(0) \nonumber\\
&=\sum_{u=1}^3\;\mathop{\sum\sum}_{1\le i\le j\le 3}\P(A({\bm e}_i+{\bm e}_j,u))E_{\max}({\bm e}_i+{\bm e}_j,u)\nonumber\\
&\qquad{}-\frac{\binom{n_1}{1}\binom{n_3}{1}}{\binom{|{\bm n}|}{2}}
\bigg(1-\frac{\binom{n_1-1}{1}\binom{n_3-1}{1}}{\binom{|{\bm n}|-2}{2}}\bigg).
\end{align*}
The second equality uses the fact that the union of the events $B({\bm e}_i+{\bm e}_j,u)$ ($1\le i\le j\le 3,\; (i,j)\ne(1,3),\;u\in\{1,3\}$) is the event that the dealer has a natural and the player does not.
Finally, we observe that, for $1\le i<j\le 3$ or $1\le i\le 3$,
\begin{align*}
\P(A({\bm e}_i+{\bm e}_j,1))&=\frac{\binom{n_i}{1}\binom{n_j}{1}}{\binom{|{\bm n}|}{2}}\,\frac{n_1-\delta_{1,i}-\delta_{1,j}}{|{\bm n}|-2}
\bigg(1-\frac{n_3-\delta_{3,i}-\delta_{3,j}}{|{\bm n}|-3}\bigg),\\
\P(A(2{\bm e}_i,1))&=\frac{\binom{n_i}{2}}{\binom{|{\bm n}|}{2}}\,\frac{n_1-2\delta_{1,i}}{|{\bm n}|-2}
\bigg(1-\frac{n_3-2\delta_{3,i}}{|{\bm n}|-3}\bigg),\\
\P(A({\bm e}_i+{\bm e}_j,2))&=\frac{\binom{n_i}{1}\binom{n_j}{1}}{\binom{|{\bm n}|}{2}}\,\frac{n_2-\delta_{2,i}-\delta_{2,j}}{|{\bm n}|-2},\\
\P(A(2{\bm e}_i,2))&=\frac{\binom{n_i}{2}}{\binom{|{\bm n}|}{2}}\,\frac{n_2-2\delta_{2,i}}{|{\bm n}|-2},\\
\P(A({\bm e}_i+{\bm e}_j,3))&=\frac{\binom{n_i}{1}\binom{n_j}{1}}{\binom{|{\bm n}|}{2}}\,\frac{n_3-\delta_{3,i}-\delta_{3,j}}{|{\bm n}|-2}\bigg(1-\frac{n_1-\delta_{1,i}-\delta_{1,j}}{|{\bm n}|-3}\bigg),\\
\P(A(2{\bm e}_i,3))&=\frac{\binom{n_i}{2}}{\binom{|{\bm n}|}{2}}\,\frac{n_3-2\delta_{3,i}}{|{\bm n}|-2}
\bigg(1-\frac{n_1-2\delta_{1,i}}{|{\bm n}|-3}\bigg),
\end{align*}
and the derivation is complete.

\section{Snackjack basic strategy results}\label{BS results}

In the case of a single deck, $(n_1,n_2,n_3)=(2,2,4)$, so by direct enumeration, $|\mathscr{L}|=14$ and $|\mathscr{M}|=32$ after we exclude $((1,2,0),1)$ from $\mathscr{M}$, as explained in Section~\ref{comparisons}.  The 87 conditional expectations (32 stand, 32 hit, 16 double, 7 split) needed for composition-dependent basic strategy are shown in Table~\ref{SJ-BS-1deck}.  The inner product of the last two columns, divided by 420, is the player's expectation under basic strategy, $27/140\approx0.192857$.

\begin{table}[htb]
\caption{\label{SJ-BS-1deck}Derivation of composition-dependent basic strategy for single-deck snackjack.  For computational convenience, rows are arranged in descending order of the hard total (htot).}
\catcode`@=\active \def@{\hphantom{0}}
\catcode`#=\active \def#{\hphantom{$\;^1$}}
\tabcolsep=.08cm
\begin{center}
\begin{tabular}{cccccccccccc}
\noalign{\smallskip}
\hline
\noalign{\smallskip}
no & nos of & htot & tot & up & $E_{\rm std}$ & $E_{\rm hit}$ & $E_{\rm dbl}$ & $E_{\rm spl}$ & bs & $E_{\rm max}$ & $420\times{}$ \\
& 1s, 2s, 3s & &&&&&&&&& probab \\
\noalign{\smallskip}
\hline
\noalign{\smallskip}
@1 & $(0,2,1)$   & 7 & h7 & 1 & $1$     & $-1$    & na     & na     & S     &         &    \\
@2 & $(0,2,1)$   & 7 & h7 & 3 & $1$     & $-1$    & na     & na     & S     &         &    \\
@3 & $(1,0,2)$   & 7 & h7 & 1 & $1$     & $-1$    & na     & na     & S     &         &    \\
@4 & $(1,0,2)$   & 7 & h7 & 2 & $2/3$   & $-1$    & na     & na     & S     &         &    \\
@5 & $(1,0,2)$   & 7 & h7 & 3 & $7/9$   & $-1$    & na     & na     & S     &         &    \\
@6 & $(2,1,1)$   & 7 & h7 & 2 & $1$     & $-1$    & na     & na     & S     &         &    \\
@7 & $(2,1,1)$   & 7 & h7 & 3 & $1$     & $-1$    & na     & na     & S     &         &    \\
@8 & $(0,0,2)$   & 6 & h6 & 1 & $-2/9$  & $-2/3$  & $-4/3$ & $-4/9$ & S     & $-2/9$  & @18 \\
@9 & $(0,0,2)$   & 6 & h6 & 2 & $-1/30$ & $-1/3$  & $-2/3$ & $1/5$  & Spl   & $1/5$   & @30 \\
10 & $(0,0,2)$   & 6 & h6 & 3 & $0$     & $-1/9$  & $-2/9$ & $2/9$  & Spl   & $2/9$   & @18 \\
11 & $(1,1,1)$   & 6 & h6 & 1 & $0$     & $-1$    & na     & na     & S     &         &    \\
12 & $(1,1,1)$   & 6 & h6 & 2 & $1/2$   & $-1/2$  & na     & na     & S     &         &    \\
13 & $(1,1,1)$   & 6 & h6 & 3 & $2/9$   & $-1/3$  & na     & na     & S     &         &    \\
14 & $(2,2,0)$   & 6 & h6 & 3 & $0$     & $-1$    & na     & na     & S     &         &    \\
15 & $(0,1,1)$   & 5 & h5 & 1 & $-1/2$  & $-5/8$  & $-5/4$ & na     & S     & $-1/2$  & @16 \\
16 & $(0,1,1)$   & 5 & h5 & 2 & $-2/5$  & $-2/5$  & $-4/5$ & na     & S/H   & $-2/5$  & @20 \\
17 & $(0,1,1)$   & 5 & h5 & 3 & $-2/3$  & $-1/18$ & $-1/9$ & na     & H     & $-1/18$ & @36 \\
18 & $(1,2,0)$   & 5 & h5 & 3 & $-1$    & $-2/3$  & na     & na     & H     &         &    \\
19 & $(2,0,1)$   & 5 & h5 & 2 & $0$     & $-1/2$  & na     & na     & S     &         &    \\
20 & $(2,0,1)$   & 5 & h5 & 3 & $-1/3$  & $0$     & na     & na     & H     &         &    \\
21 & $(0,2,0)$   & 4 & h4 & 1 & $1$     & $1$     & $2$    & $2$    & D/Spl & $2$     & @@1 \\
22 & $(0,2,0)$   & 4 & h4 & 3 & $-1$    & $1/6$   & $0$    & $-1$   & H     & $1/6$   & @@6 \\
23 & $(1,0,1)$   & 4 & s7 & 1 & $3/2$   & $3/4$   & $3/2$  & na     & S/D   & $3/2$   & @@8 \\
24 & $(1,0,1)$   & 4 & s7 & 2 & $3/2$   & $1/2$   & $1$    & na     & S     & $3/2$   & @40 \\
25 & $(1,0,1)$   & 4 & s7 & 3 & $3/2$   & $3/8$   & $7/12$ & na     & S     & $3/2$   & @48 \\
26 & $(2,1,0)$   & 4 & s7 & 2 & $1$     & $1$     & na     & na     & S/H   &         &    \\
27 & $(2,1,0)$   & 4 & s7 & 3 & $1$     & $3/4$   & na     & na     & S     &         &    \\
28 & $(1,1,0)$   & 3 & s6 & 1 & $0$     & $0$     & $0$    & na     & S/H/D & $0$     & @@2 \\
29 & $(1,1,0)$   & 3 & s6 & 2 & $3/5$   & $3/5$   & $6/5$  & na     & D     & $6/5$   & @10 \\
30 & $(1,1,0)$   & 3 & s6 & 3 & $3/16$  & $1/4$   & $3/8$  & na     & D     & $3/8$   & @32 \\
31 & $(2,0,0)$   & 2 & s5 & 2 & $1/5$   & $1/5$   & $2/5$  & $6/5$  & Spl   & $6/5$   & @@5 \\
32 & $(2,0,0)$   & 2 & s5 & 3 & $-2/5$  & $2/5$   & $2/5$  & $6/5$  & Spl   & $6/5$   & @10 \\
\noalign{\smallskip}
\hline
\noalign{\smallskip}
\multicolumn{10}{l}{dealer has natural, player does not} & $-1$ & @96 \\
\multicolumn{10}{l}{both player and dealer have naturals} & $0$ & @24 \\
\noalign{\smallskip}
\hline
\noalign{\smallskip}\multicolumn{11}{l}{total} & 420 \\
\end{tabular}
\end{center}
\end{table}
\afterpage{\clearpage}

To clarify our method for determining composition-dependent basic strategy, we provide several examples of how the conditional expectations in Table~\ref{SJ-BS-1deck} were computed.  Let us denote by $Y$ the player's next card and by $\bm Z=(Z_1,Z_2,\ldots)$ the dealer's hand beginning with $Z_1=U$ and $Z_2=D$.  $Z_3$ and $Z_4$ would be the third and fourth cards in the dealer's hand if needed.  No more than four cards are ever needed.

\begin{figure}[h]
\setlength{\unitlength}{0.65cm}
\begin{picture}(19,8)(.8,0)
\put(1,4.6){player's}
\put(1,3.9){$\text{hand}=3,3$}
\put(1,3){$\text{upcard}=1$}
\put(1,2.1){$\text{remainder}$}
\put(1,1.4){$=(1,2,2)$}
\put(4,3){\line(1,1){2}}
\put(4,3){\line(1,-1){2}}
\put(4.2,4.2){1/3}
\put(5.2,2){2/3}
\put(6.2,5){$1,1=\text{s}5$}
\put(6.2,4.1){$\text{remainder}$}
\put(6.2,3.4){$=(0,2,2)$}
\put(6.2,0.8){$1,2=\text{s}6$}
\put(9,5){\line(1,1){2}}
\put(9,5){\line(1,-1){2}}
\put(9.2,6.2){2/4}
\put(10.2,4){2/4}
\put(11.15,6.9){$1,1,2=\text{s}7$}
\put(11.15,2.9){$1,1,3=\text{h}5$}
\put(11.15,2){$\text{remainder}$}
\put(11.15,1.3){$=(0,2,1)$}
\put(14,2.9){\line(1,1){2}}
\put(14,2.9){\line(1,-1){2}}
\put(14.2,4.1){2/3}
\put(15.2,1.9){1/3}
\put(16.2,4.8){$1,1,3,2=\text{h}7$}
\put(16.2,0.7){$1,1,3,3=\text{h}8$}
\end{picture}
\caption{\label{tree}The tree diagram used to evaluate $E_\std(2\bm e_3,1)$. Notice that we are conditioning on the dealer not having a natural (i.e., the dealer's downcard is not a trey).}
\end{figure}
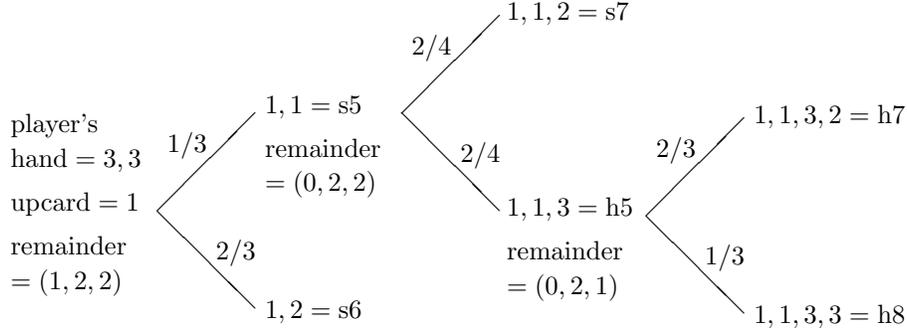

First, the conditional expectation when standing with $3,3$ vs.~1 can be evaluated with a tree diagram.  See Figure~\ref{tree}.  More formally,
\begin{align*}
&\!\!\!E_\std(2\bm e_3,1)\\
&=\E[G_\std\mid{\bm X}=2\bm e_3,\,U=1,\,D\ne3)\\
&=\P(S=6\mid {\bm X}=2\bm e_3,\,U=1,\,D\ne3)(0)\\
&\qquad{}+\P(S=7\mid {\bm X}=2\bm e_3,\,U=1,\,D\ne3)(-1)\\
&\qquad{}+\P(S=8\mid {\bm X}=2\bm e_3,\,U=1,\,D\ne3)(1)\\
&=\P(\bm Z=(1,2)\mid{\bm X}=2\bm e_3,\,U=1,\,D\ne3)(0)\\
&\qquad{}+\P(\bm Z=(1,1,2)\text{ or }(1,1,3,2)\mid{\bm X}=2\bm e_3,\,U=1,\,D\ne3)(-1)\\
&\qquad{}+\P(\bm Z=(1,1,3,3)\mid{\bm X}=2\bm e_3,\,U=1,\,D\ne3)(1)\\
&=\frac23(0)+\bigg(\frac13\,\frac24+\frac13\,\frac24\,\frac23\bigg)(-1)+\bigg(\frac13\,\frac24\,\frac13\bigg)(1)=-\frac29.
\end{align*} 

Second, the conditional expectation when hitting with $1,2$ vs.~3 is
\begin{align*}
&\!\!\!E_\hit(\bm e_1+\bm e_2,3)\\
&=\E[G_\hit\mid{\bm X}=\bm e_1+\bm e_2,\,U=3,\,D\ne1)\\
&=\P(Y=1\mid{\bm X}=\bm e_1+\bm e_2,\,U=3,\,D\ne1)E_\std((2,1,0),3)\\
&\qquad{}+\P(Y=2\mid{\bm X}=\bm e_1+\bm e_2,\,U=3,\,D\ne1)E_\hit((1,2,0),3)\\
&\qquad{}+\P(Y=3\mid{\bm X}=\bm e_1+\bm e_2,\,U=3,\,D\ne1)E_\std((1,1,1),3)\\
&=\frac15\,\frac{1-0}{1-\frac15}(1)+\frac15\,\frac{1-\frac14}{1-\frac15}\bigg(\!\!-\frac23\bigg)+\frac35\,\frac{1-\frac14}{1-\frac15}\,\bigg(\frac29\bigg)=\frac14,
\end{align*}
where \eqref{p(k,3)} was used to evaluate the conditional probabilities and we have used the facts that $E_\std((2,1,0),3)$, $E_\hit((1,2,0),3)$, and $E_\std((1,1,1),3)$ have already been computed, and are larger than $E_\hit((2,1,0),3)$, $E_\std((1,2,0),3)$, and $E_\hit((1,1,1),3)$, respectively.

Third, the conditional expectation when doubling with $1,3$ vs.~1 is
\begin{align*}
&\!\!\!E_\dbl(\bm e_1+\bm e_3,1)\\
&=\E[G_\dbl\mid{\bm X}=\bm e_1+\bm e_3,\,U=1,\,D\ne3)\\
&=2\{\P(Y=2\mid{\bm X}=\bm e_1+\bm e_3,\,U=1,\,D\ne3)E_\std((1,1,1),1)\\
&\qquad{}+\P(Y=3\mid{\bm X}=\bm e_1+\bm e_3,\,U=1,\,D\ne3)E_\std((1,0,2),1)\}\\
&=2\bigg[\frac25\,\frac{1-\frac34}{1-\frac35}(0)+\frac35\,\frac{1-\frac24}{1-\frac35}(1)\bigg]=\frac32,
\end{align*}
using \eqref{p(k,1)}.

Finally, the conditional expectation when splitting with $3,3$ vs.~2 is
\begin{align*}
&\!\!\!E_\spl(2\bm e_3,2)\\
&=2\sum_{j=1}^3\P(Y=j\mid\bm X=2\bm e_3,\, U=2)E_\std(\bm e_3+\bm e_j,2\mid \bm e_3)\\
&=2\bigg[\frac25 E_\std(\bm e_3+\bm e_1,2\mid\bm e_3)+\frac15 E_\std(\bm e_3+\bm e_2,2\mid\bm e_3)+\frac25 E_\std(2\bm e_3,2\mid\bm e_3)\bigg]\\
&=2\,\frac25[\P_{-\bm e_3}(\bm Z=(2,2,3)\text{ or }(2,3,2)\mid \bm X=\bm e_3+\bm e_1,\,U=2)(0)\\
&\qquad\quad\;{}+\P_{-\bm e_3}(\bm Z=(2,1),(2,2,1,3),(2,3,1),\text{ or }(2,3,3)\mid \\
&\qquad\qquad\qquad\qquad\qquad\qquad\qquad\qquad\qquad\qquad\bm X=\bm e_3+\bm e_1,\,U=2)(1)]\\
&\quad{}+2\,\frac15[\P_{-\bm e_3}(\bm Z=(2,1)\text{ or }(2,3,1)\mid\bm X=\bm e_3+\bm e_2,\,U=2)(-1)\\
&\qquad\qquad{}+\P_{-\bm e_3}(\bm Z=(2,3,3)\mid\bm X=\bm e_3+\bm e_2,\,U=2)(1)]\\
&\quad{}+2\,\frac25[\P_{-\bm e_3}(\bm Z=(2,1),(2,2,1,1),\text{ or }(2,3,1)\mid\bm X=2\bm e_3,\,U=2)(0)\\
&\qquad\qquad{}+\P_{-\bm e_3}(\bm Z=(2,2,3)\text{ or }(2,3,2)\mid\bm X=2\bm e_3,\,U=2)(-1)\\
&\qquad\qquad{}+\P_{-\bm e_3}(\bm Z=(2,2,1,3)\mid\bm X=2\bm e_3,\,U=2)(1)]\\
&=2\,\frac25\bigg[\bigg(\frac14\,\frac23+\frac24\,\frac13\bigg)(0)+\bigg(\frac14+\frac14\,\frac13\,\frac22+\frac24\,\frac13+\frac24\,\frac13\bigg)(1)\bigg]\\
&\qquad{}+2\,\frac15\bigg[\bigg(\frac24+\frac24\,\frac23\bigg)(-1)+\bigg(\frac24\,\frac13\bigg)(1)\bigg]\\
&\qquad{}+2\,\frac25\bigg[\bigg(\frac24+\frac14\,\frac23\,\frac12+\frac14\,\frac23\bigg)(0)+\bigg(\frac14\,\frac13+\frac14\,\frac13\bigg)(-1)+\bigg(\frac14\,\frac23\,\frac12\bigg)(1)\bigg]\\
&=2\bigg[\frac25\,\frac23+\frac15\,\bigg(\!\!-\frac23\bigg)+\frac25\bigg(\!\!-\frac{1}{12}\bigg)\bigg]=\frac15,
\end{align*}
where the subscript $-\bm e_3$ means that the deck has been depleted by one trey.

In the case of two decks, $(n_1,n_2,n_3)=(4,4,8)$, so $|\mathscr{L}|=23$ and $|\mathscr{M}|=66$.  

In the case of three decks, $(n_1,n_2,n_3)=(6,6,12)$, so $|\mathscr{L}|=26$ and $|\mathscr{M}|=77$. 

In the case of $d$ decks with $d\ge4$, $(n_1,n_2,n_3)=(2d,2d,4d)$, so $|\mathscr{L}|=27$ and $|\mathscr{M}|=81$.  Of course, some of these 27 hands are never seen by the basic strategist.  For example, $(3,0,0)$, $(4,0,0)$, $(5,0,0)$, $(6,0,0)$, and $(7,0,0)$ are never encountered because the basic strategist splits $(2,0,0)$.

In Table~\ref{SJ-BS-d-deck-condensed} we present basic strategy for $d$ decks, where $d$ is a positive integer.  For $d\ge9$, composition-dependent basic strategy does not depend on $d$.  It may be surprising that basic strategy has very little dependence on the dealer's upcard, but that is the nature of snackjack.  See Appendix~B for a more complete description of composition-dependent basic strategy.  Overall player expectation, as a function of the number of decks, is shown in Table~\ref{EV(d)}.  As in blackjack~\cite[p.~177]{G99},~\cite[p.~394]{S18},~\cite[p.~11]{We18}, it is a decreasing function of $d$.

\begin{table}[htb]
\caption{\label{SJ-BS-d-deck-condensed}Composition-dependent basic strategy for $d$-deck snackjack, $d$ a positive integer.  For $d=1$ there are five decision points where composition-dependent basic strategy is nonunique; for $d=2$ there is one.  We have excluded exceptions that do not occur to the basic strategist.  For example, with $1,2,2$ (hard 5) vs.~2 it is correct to stand if $2\le d\le6$, but the basic strategist doubles $1,2$ vs.~2 and $2,2$ vs.~2, so this exception never arises.  For more-complete tables, see Appendix~B.}
\catcode`@=\active \def@{\hphantom{0}}
\catcode`#=\active \def#{\hphantom{$\;^1$}}
\tabcolsep=.08cm
\begin{center}
\begin{tabular}{cccc}
\hline
\noalign{\smallskip}
player & \multicolumn{3}{c}{dealer upcard} \\
\noalign{\smallskip}
total & $1$ & $2$ & $3$ \\
\noalign{\smallskip}
\hline
\noalign{\smallskip}
@hard 7@&  S & S & S \\
\noalign{\smallskip}
@hard 6@&  S & S & S \\
\noalign{\smallskip}   
@hard 5@&  #H$\,^1$ & #H# & H \\
\noalign{\smallskip}
\hline
\noalign{\smallskip}       
@soft 7@&  S & S & S \\
\noalign{\smallskip} 
@soft 6@&  #H# & #D# & D \\ 
\noalign{\smallskip}
\hline
\noalign{\smallskip}         
$(3,3)$ &  #Spl$\,^2$ & Spl & #Spl# \\
\noalign{\smallskip} 
$(2,2)$ &  D & D & #D$\,^3$ \\
\noalign{\smallskip} 
$(1,1)$ &  Spl & #Spl# & Spl \\
\noalign{\smallskip}
\hline
\noalign{\smallskip}
\multicolumn{4}{l}{$^1\,$S if $d=1$. $^2\,$S if $d\le8$.}\\
\multicolumn{4}{l}{$^3\,$H if $d=1$.}\\
\end{tabular}
\end{center}
\end{table}

\begin{table}[htb]
\caption{\label{EV(d)}Player expectation at $d$-deck snackjack under composition-dependent basic strategy, as a function of $d$.}
\catcode`@=\active \def@{\hphantom{0}}
\catcode`#=\active \def#{\hphantom{$\;^1$}}
\tabcolsep=.2cm
\begin{center}
\begin{tabular}{cccccccc}
\hline
\noalign{\smallskip}
$d$ & expectation && $d$ & expectation && $d$ & expectation \\
\noalign{\smallskip}
\hline
\noalign{\smallskip}
1 & 0.192857 &@@& @7 & 0.144558 &@@& 13 & 0.141548 \\
2 & 0.163144 && @8 & 0.143639 && 26 & 0.139871 \\
3 & 0.154360 && @9 & 0.143031 && 39 & 0.139309 \\
4 & 0.150073 && 10 & 0.142550 && 52 & 0.139028 \\
5 & 0.147500 && 11 & 0.142156 &&    & \\
6 & 0.145784 && 12 & 0.141827 && $\infty$ & 0.138184 \\
\noalign{\smallskip}
\hline
\end{tabular}
\end{center}
\end{table}

\section{Potential gain from bet variation}\label{FTCC}

The fundamental theorem of card counting~\cite{EL05,TW73} (see~\cite[Section 11.3]{Et10} for a textbook treatment) tells us that the player's conditional expectation under basic strategy, given the $n$ cards seen so far, is a random variable with mean that is constant in $n$, mean positive part that is nondecreasing in $n$, and standard deviation that is increasing in $n$.  

To illustrate in the simplest possible situation, we consider the toy game of red-and-black mentioned in Section~\ref{intro}, which could just as well be \textit{odd-and-even}.  An advantage of the latter formulation is that the cards can be numbered from 1 to $N$ ($N$ is the size of the deck, assumed even), and then
$$
Z_n:=\frac{1}{N-n}\sum_{i=1}^n (-1)^{X_i}
$$
gives the exact player conditional expectation of a one-unit even-money bet that the next card dealt is odd, given that the first $n$ cards, $X_1,X_2,\ldots,X_n$, have been seen.  It is easy to verify that $\E[Z_n]=0$,
\begin{align}\label{E[Z_n^+]}
\E[(Z_n)^+]&=\frac{1}{N-n}\sum_{k=0}^{\lfloor n/2\rfloor}(n-2k)\frac{\binom{N/2}{k}\binom{N/2}{n-k}}{\binom{N}{n}},\\  
\text{SD}(Z_n)&=\sqrt{\frac{n}{(N-n)(N-1)}},\label{SD(Z_n)}
\end{align}
for $n=1,2,\ldots,N-1$, where $a^+:=\max(a,0)$.  The expectation \eqref{E[Z_n^+]} is nondecreasing in $n$, and the standard deviation \eqref{SD(Z_n)} is increasing in $n$, both as a result of the FTCC~\cite{EL05}.  (We cannot express \eqref{E[Z_n^+]} in closed form, but we can show analytically that $\E[(Z_n)^+]\le\frac12\,\text{SD}(Z_n)$.) 

Let us consider a shoe comprising 39 decks, or 312 cards, at snackjack.  Basic strategy is the strategy of Table~\ref{SJ-BS-d-deck-condensed} without the footnotes.  To summarize it, the player mimics the dealer except when he has a soft 6 or a pair.  He hits a soft 6 against an ace and otherwise doubles.  He splits a pair of aces and a pair of treys and doubles a pair of deuces.  A single round with one player can be completed with certainty if at least eight cards remain.   The mean profit, given that $n_1$ aces, $n_2$ deuces, and $n_3$ treys remain, is
\begin{equation}\label{E(n_1,n_2,n_3)}
E(n_1,n_2,n_3)=\frac{P(n_1,n_2,n_3)}{(n_1+n_2+n_3)_8},
\end{equation}
where $P(n_1,n_2,n_3)$ is a polynomial of degree 8 in $n_1$, $n_2$, and $n_3$ with 147 terms (see Appendix~C), and $(N)_8:=N(N-1)\cdots(N-7)$.  For example,    

\begin{align*}
&E(2d,2d,4d)\\
&\quad{}=\frac{-630 + 4{,}017 d - 2{,}673 d^2 - 32{,}132 d^3 + 92{,}560 d^4 - 97{,}144 d^5 + 36{,}224 d^6}{(8d-1)_3(8d-5)_3},
\end{align*} 
which yields the entries in Table~\ref{EV(d)} for $d\ge9$ because basic strategy optimized for $d$ decks coincides with 39-deck basic strategy provided $d\ge9$.  As a check of \eqref{E(n_1,n_2,n_3)}, we can confirm that $E(n_1,0,0)=-2$ (player splits, gets two soft 5s, dealer wins both with soft 6), $E(0,n_2,0)=0$ (player doubles, gets hard 6, dealer pushes with hard 6), and $E(0,0,n_3)=0$ (player splits, gets two hard 6s, dealer pushes both with hard 6).

Because of the simplicity of snackjack, we can compute the means and variances arising in the fundamental theorem of card counting.  The analogous computations at blackjack would be prohibitively time-consuming.  To justify this claim, we need to do some counting.  

In 39-deck snackjack, if $n$ cards have been seen, the numbers $M_1$, $M_2$, and $M_3$ of aces, deuces, and treys among them are such that $(M_1,M_2,M_3)$ has the multivariate hypergeometric distribution
$$
\P(M_1=m_1,\,M_2=m_2,\,M_3=m_3)=\frac{\binom{78}{m_1}\binom{78}{m_2}\binom{156}{m_3}}{\binom{312}{n}},
$$
where $0\le m_1\le78$, $0\le m_2\le78$, $0\le m_3\le156$, and $m_1+m_2+m_3=n$.  The number $s(n)$ of distinct values of $(M_1,M_2,M_3)$ such that $M_1+M_2+M_3=n$ satisfies $s(n)=s(312-n)$ and is given by
\begin{align}\label{s(n)}
s(n)&=\sum_{k=0}^2(-1)^k\binom{2}{k}\bigg[\binom{n+2-79k}{2}I(n\ge79k)\nonumber\\
&\qquad\qquad\qquad\qquad{}-\binom{n+2-79k-157}{2}I(n\ge79k+157)\bigg];
\end{align}
see below for details.  In particular, $\max_n s(n)=s(156)=6{,}241$ and $\sum_n s(n)=(79)^2 157=979{,}837$.  That is, there are fewer than one million distinguishable subsets of the 39-deck snackjack shoe.

In 24-deck grayjack, if $n$ cards have been seen, the numbers $M_1$, $M_2$, \dots, $M_6$ of aces, 2s, \dots, 6s among them are such that $(M_1,M_2,\ldots,M_6)$ has the multivariate hypergeometric distribution
$$
\P(M_1=m_1,\,M_2=m_2,\ldots,M_6=m_6)=\frac{\binom{24}{m_1}\big[\prod_{i=2}^5\binom{48}{m_i}\big]\binom{96}{m_6}}{\binom{312}{n}},
$$
where $0\le m_1\le24$, $0\le m_i\le48$ for $2\le i\le5$, $0\le m_6\le96$, and $m_1+m_2+\cdots+m_6=n$.  The number $g(n)$ of distinct values of $(M_1,M_2,\ldots,M_6)$ such that $M_1+M_2+\cdots+M_6=n$ satisfies $g(n)=g(312-n)$ and is given by
\begin{align}\label{g(n)}
g(n)&=\sum_{k=0}^4(-1)^k\binom{4}{k}\bigg[\binom{n+5-49k}{5}I(n\ge49k)\nonumber\\
&\qquad\qquad\qquad\qquad{}-\binom{n+5-25-49k}{5}I(n\ge25+49k)\nonumber\\
&\qquad\qquad\qquad\qquad{}-\binom{n+5-49k-97}{5}I(n\ge49k+97)\\
&\qquad\qquad\qquad\qquad{}+\binom{n+5-25-49k-97}{5}I(n\ge25+49k+97)\bigg].\nonumber
\end{align}
In particular, $\max_n g(n)=g(156)=130{,}046{,}539$ and $\sum_n g(n)=25(49)^4 97=13{,}979{,}642{,}425$.  That is, there are about 14 billion distinguishable subsets of the 24-deck grayjack shoe.

In six-deck blackjack, if $n$ cards have been seen, the numbers $M_1$, $M_2$, \dots, $M_{10}$ of aces, 2s, \dots, tens among them are such that $(M_1,M_2,\ldots,M_{10})$ has the multivariate hypergeometric distribution
$$
\P(M_1=m_1,\,M_2=m_2,\ldots,M_{10}=m_{10})=\frac{\big[\prod_{i=1}^9\binom{24}{m_i}\big]\binom{96}{m_{10}}}{\binom{312}{n}},
$$
where $0\le m_i\le24$ for $1\le i\le9$, $0\le m_{10}\le96$, and $m_1+m_2+\cdots+m_{10}=n$.  The number $b(n)$ of distinct values of $(M_1,M_2,\ldots,M_{10})$ such that $M_1+M_2+\cdots+M_{10}=n$  satisfies $b(n)=b(312-n)$ and is given by
\begin{align}\label{b(n)}
b(n)&=\sum_{k=0}^9(-1)^k\binom{9}{k}\bigg[\binom{n+9-25k}{9}I(n\ge25k)\nonumber\\
&\qquad\qquad\qquad\qquad{}-\binom{n+9-25k-97}{9}I(n\ge25k+97)\bigg].
\end{align}
Therefore, $\max_n b(n)=b(156)=3{,}726{,}284{,}230{,}655$ and $\sum_n b(n)=(25)^9 97=370{,}025{,}634{,}765{,}625$.  That is, there are more than 370 trillion distinguishable subsets of the six-deck blackjack shoe. 

Denoting by $b_1(n)$ the analogous quantity in single-deck blackjack, we have $b_1(n)=b_1(52-n)$, $\max_n b_1(n)=b_1(26)=1{,}868{,}755$, and $\sum_n b_1(n)=5^9 17=33{,}203{,}125$.  See Griffin~\cite[p.~159]{G99} and Thorp~\cite[p.~126]{T00}.

To clarify how \eqref{s(n)}--\eqref{b(n)} were derived, we elaborate on \eqref{s(n)}.  Let
\begin{align*}
A&:=\{(m_1,m_2,m_3):m_1\ge0,\,m_2\ge0,\,m_3\ge0,\,m_1+m_2+m_3=n\},\\
B_1&:=\{(m_1,m_2,m_3):m_1\ge79,\,m_2\ge0,\,m_3\ge0,\,m_1+m_2+m_3=n\},\\
B_2&:=\{(m_1,m_2,m_3):m_1\ge0,\,m_2\ge79,\,m_3\ge0,\,m_1+m_2+m_3=n\},\\
B_3&:=\{(m_1,m_2,m_3):m_1\ge0,\,m_2\ge0,\,m_3\ge157,\,m_1+m_2+m_3=n\}.
\end{align*}
Then $|A|=\binom{n+2}{2}$, $|B_1|=\binom{n+2-79}{2}I(n\ge79)$, $|B_1\cap B_3|=\binom{n+2-79-157}{2}I(n\ge79+157)$,
and so on.  By inclusion-exclusion,
\begin{align*}
&\!\!\!\!\!|A-(B_1\cup B_2\cup B_3)|\\
&=|A|-|B_1|-|B_2|-|B_3|+|B_1\cap B_2|+|B_1\cap B_3|+|B_2\cap B_3|\\
&=|A|-|B_3|-2(|B_1|-|B_1\cap B_3|)+|B_1\cap B_2|,
\end{align*}
where we have used $B_1\cap B_2\cap B_3=\varnothing$, $|B_1|=|B_2|$, and $|B_1\cap B_3|=|B_2\cap B_3|$, and the result follows.

Returning to snackjack, let $Z_n$ denote the player's conditional expectation, given that $n$ cards have been seen.  Then 
\begin{align}\label{FTCC-mean}
\E[Z_n]&=\sum_{m_1+m_2+m_3=n}\frac{\binom{78}{m_1}\binom{78}{m_2}\binom{156}{m_3}}{\binom{312}{n}}\,E(78-m_1,78-m_2,156-m_3)\nonumber\\
&=\E[Z_0]=E(78,78,156)=\frac{220{,}204{,}549{,}189}{1{,}580{,}689{,}046{,}285}=:\mu,
\end{align}
where $0\le m_1\le78$, $0\le m_2\le78$, and $0\le m_3\le156$ in the sum.  The second equality is a consequence of the martingale property of $\{Z_n\}$, for which see Ethier and Levin~\cite{EL05}.

The expected positive part of the difference between the player's conditional expectation and a positive number $\nu$ is
\begin{align}\label{FTCC-meanpos}
&\E[(Z_n-\nu)^+]\\
&\;{}=\sum_{m_1+m_2+m_3=n}\frac{\binom{78}{m_1}\binom{78}{m_2}\binom{156}{m_3}}{\binom{312}{n}}[E(78-m_1,78-m_2,156-m_3)-\nu]^+,\nonumber
\end{align}
where the sum is constrained as in \eqref{FTCC-mean}.  To interpret this, suppose that players are required to pay a commission $\nu$ per unit bet initially on each hand.  (Doubling and splitting do not require any additional commission.)
Then this is the player's expected profit, assuming $n$ cards have been seen and assuming he bets one unit if and only if his net conditional expectation (taking the commission into account) is nonnegative.

Finally, the variance of the player's conditional expectation is
\begin{align}\label{FTCC-variance}
\text{Var}(Z_n)&=\sum_{m_1+m_2+m_3=n}\frac{\binom{78}{m_1}\binom{78}{m_2}\binom{156}{m_3}}{\binom{312}{n}}[E(78-m_1,78-m_2,156-m_3)-\mu]^2,
\end{align}
and again the sum is constrained as in \eqref{FTCC-mean}.

The quantities \eqref{FTCC-mean}--\eqref{FTCC-variance} can be computed for $n=1,2,3,\ldots,304\;(=312-8)$, and Figure~\ref{meanpos,sdev} displays the graph of $f(n):=\E[(Z_n-\nu)^+]$ with $\nu=1/7$, as well as the graph of the standard deviations $g(n):=\text{SD}(Z_n)$.  The two curves have similar shapes, and are also very similar to the graphs of \eqref{E[Z_n^+]} and \eqref{SD(Z_n)} with $N=312$.  They increase gradually over the first 2/3 of the shoe and more rapidly over the final 1/6.  The increase in the slope is gradual throughout, unlike with the famous ``hockey stick graph'' of climate science.

Notice that, for $n=1,2,3,\ldots,304$, each of the quantities in \eqref{FTCC-mean}--\eqref{FTCC-variance} requires up to $6{,}241$ evaluations of the rational function $E(n_1,n_2,n_3)$, which is computationally routine.  The corresponding quantities in blackjack would require up to 3.7 trillion evaluations of the basic strategy expectation (which itself is too complicated to be usefully expressed as a rational function; see Table~\ref{comparisons-multideck}), and would be computationally prohibitive.

\begin{figure}[htb]
\centering
\includegraphics[width=2.25in]{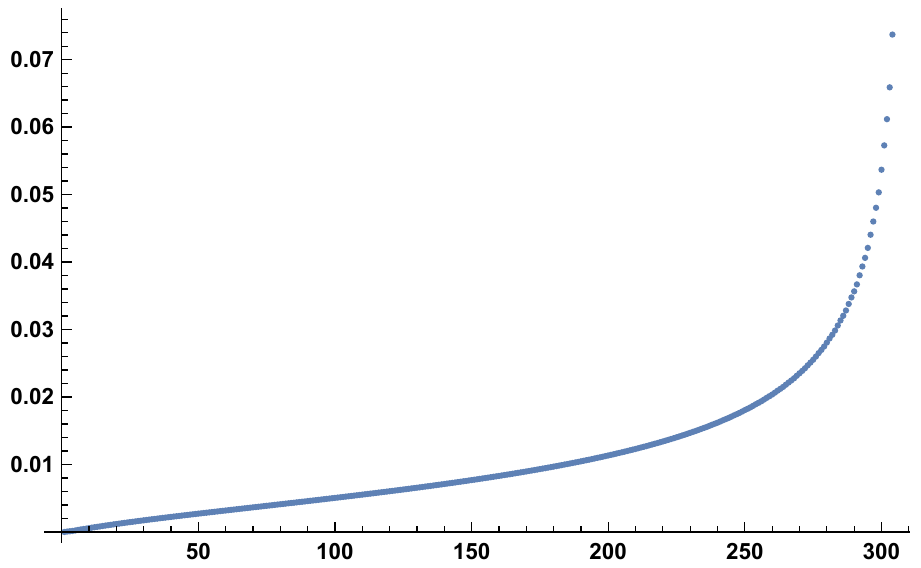}\quad\includegraphics[width=2.25in]{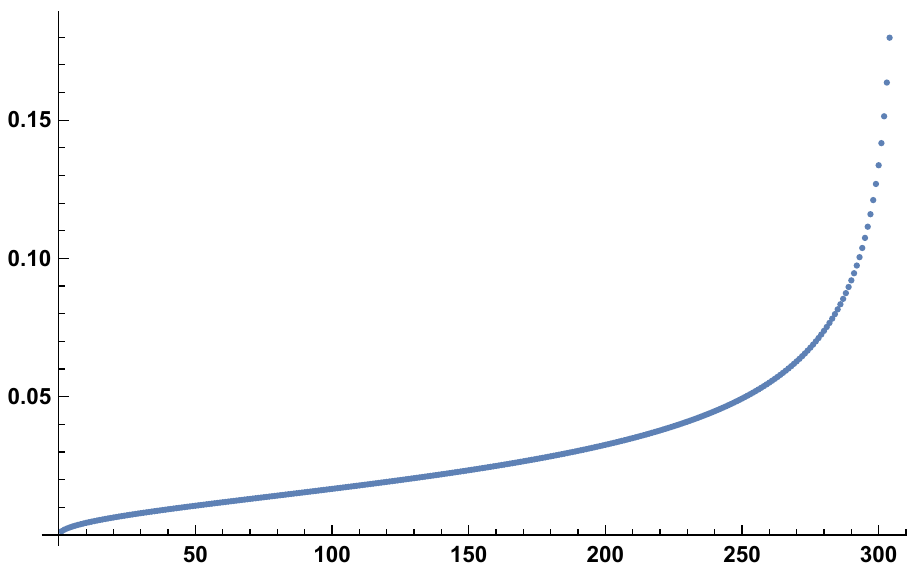}
\caption{\label{meanpos,sdev}On the left is the graph of $f(n):=\E[(Z_n-\nu)^+]$ with $\nu=1/7$, $1\le n\le304$, for 39-deck snackjack.  On the right is the graph of $g(n):=\text{SD}(Z_n)$, $1\le n\le304$, for the same game.}
\end{figure}

\section{Card counting and bet variation}\label{bet variation}

It is well known in blackjack~\cite[Chap.~4]{G99} that when the point values of a card-counting system are highly correlated with the effects of removal, a high betting efficiency is achieved but not necessarily a high strategic efficiency.  In snackjack, we continue to treat the case of a 39-deck, 312-card, shoe.  

Let us denote by $\mu({\bm m})$, where ${\bm m}=(m_1,m_2,m_3)$, the expected profit from an initial one-unit snackjack wager, assuming composition-dependent basic strategy (optimized for the 39-deck shoe), when the 39-deck shoe is depleted by $m_1$ aces, $m_2$ deuces, and $m_3$ treys.  Using \eqref{E(n_1,n_2,n_3)}, 
$$
\mu(\bm m)=E(78-m_1,78-m_2,156-m_3).
$$
We can then evaluate the \textit{effects of removal} on the expected profit from an initial one-unit snackjack wager, assuming composition-dependent basic strategy:
\begin{equation}\label{EoR}
\text{EoR}(i):=\mu({\bm e}_i)-\mu({\bm 0}),\qquad i=1,2,3.
\end{equation}
The numbers \eqref{EoR}, multiplied by $311$, are 
$$
E_1=-\frac{849{,}581{,}527}{1{,}793{,}859{,}330},\;\;E_2=\frac{3{,}539{,}587{,}453}{5{,}082{,}601{,}435},\;\;E_3=-\frac{6{,}794{,}638{,}759}{60{,}991{,}217{,}220},
$$
with decimal equivalents listed in Table~\ref{EoRs}.  A simple probabilistic argument shows that
\begin{equation}\label{EoR constraint}
E_1+E_2+2E_3=0.
\end{equation}
We provide the exact fractions above to allow confirmation that \eqref{EoR constraint} holds exactly, not just to a certain number of decimal places.

\begin{table}[htb]
\caption{\label{EoRs}Effects of removal, multiplied by 311, for an initial one-unit bet in 39-deck snackjack, assuming composition-dependent basic strategy for the 39-deck shoe.  Results are rounded to six decimal places.  Also included are two card-counting systems, one of level one, the other of level six.\medskip}
\catcode`#=\active \def#{\hphantom{$-$}}
\catcode`@=\active \def@{\hphantom{0}}
\renewcommand{\arraystretch}{1.}
\begin{center}
\begin{tabular}{cccccc}
\hline
\noalign{\smallskip}
 card  && $E_i:=$ & level one & level six \\
value $i$ && $311\,{\rm EoR}(i)$ & system & system \\
\noalign{\smallskip}\hline
\noalign{\smallskip}
  @1 && $-0.473605$ & $-1$ & $-4$ \\
  @2 && #$0.696413$ &   #1 &   #6 \\
  @3 && $-0.111404$ &   #0 & $-1$ \\
\noalign{\smallskip}\hline
\noalign{\smallskip}
\multicolumn{3}{l}{correlation $\rho$} & @$0.965597$ & @$0.999921$ \\
\noalign{\smallskip}\hline
\noalign{\smallskip}
\multicolumn{3}{l}{regression coefficient $\gamma$} & @$0.585009$ & @$0.116587$ \\
\noalign{\smallskip}\hline
\end{tabular}
\end{center}
\end{table}

In the blackjack literature (e.g., Schlesinger~\cite[pp.~503--504, 522]{S18}), it is conventional to evaluate the effects of removal for the single-deck game and then use a conversion factor to handle the multiple-deck games.  We could follow this precedent with $6\frac12$-deck, 52-card, snackjack playing the role of single-deck blackjack, but we prefer to work directly with the 39-deck, 312-card, game.

Recall that, in a \textit{balanced} card-counting system, the sum of the point values over the entire pack is 0.  For the system $(J_1,J_2,J_3)$, this means that
$$
J_1+J_2+2J_3=0.
$$
Table~\ref{EoRs} lists two balanced card-counting systems, the best level-one system and the best level-six system, the level being defined by $\max(|J_1|,|J_2|,|J_3|)$.  In each case we indicate the correlation $\rho$ with the effects of removal, and the relevant regression coefficient $\gamma$ defined below.  Based on~\cite[Eqs.\ (11.76), (11.95), (21.69), and (21.70)]{Et10}, an estimate of $Z_n$ (the player's conditional expectation under basic strategy, given that $n$ cards have been seen) is
\begin{equation}\label{Z_n hat}
\widehat Z_n:=\mu+\frac{1}{312-n}\sum_{j=1}^n E_{X_j},
\end{equation}
where $E_i:=311\,\text{EoR}(i)$ and $X_1,X_2,\ldots,X_{312}$ is the sequence of card values in the order in which they are exposed,  which in turn is approximated by
\begin{equation}\label{Z_n^*}
Z_n^*:=\mu+\frac{\gamma}{52}\bigg(\frac{52}{312-n}\sum_{j=1}^n J_{X_j}\bigg)=\mu+\frac{\gamma}{52}\text{TC}_n,
\end{equation}
where $(J_1,J_2,J_3)$ is one of the two card-counting systems listed in Table~\ref{EoRs}, and 
$$
\gamma:=\frac{E_1 J_1+E_2 J_2+2E_3 J_3}{J_1^2+J_2^2+2J_3^2}
$$
is the regression coefficient that minimizes the sum of squares $(E_1-\gamma J_1)^2+(E_2-\gamma J_2)^2+2(E_3-\gamma J_3)^2$.  We find that
\begin{align*}
\rho&\approx0.965597,\quad\gamma=\frac{35{,}680{,}410{,}677}{60{,}991{,}217{,}220},\quad\text{if }(J_1,J_2,J_3)=(-1,1,0),\\
\rho&\approx0.999921,\quad\gamma=\frac{63{,}997{,}110{,}301}{548{,}920{,}954{,}980},\quad\text{if }(J_1,J_2,J_3)=(-4,6,-1).
\end{align*}
Finally, $\text{TC}_n$ is the \textit{true count}, which is the \textit{running count} (the sum of the point values of the cards seen so far) divided by the number of unseen 52-card packs, namely $(312-n)/52$.  We use 52 instead of  8 here because it may be easier to estimate the number of unseen 52-card packs than the number of unseen 8-card decks.
      
The first question we would like to address is, how accurate is card counting?  There are several ways to answer this question, but a first step would be to compare $Z_n$ with its approximations $\widehat Z_n$ and $Z_n^*$.  More specifically, we compare the $L^1$ distances between $Z_n$ and its approximations. So we evaluate
\begin{align*}
\|Z_n-\widehat Z_n\|_1&=\sum_{m_1+m_2+m_3=n}\frac{\binom{78}{m_1}\binom{78}{m_2}\binom{156}{m_3}}{\binom{312}{n}}\bigg|E(78-m_1,78-m_2,156-m_3)\\
&\qquad\qquad\qquad\qquad{}-\bigg(\mu+\frac{1}{312-n}(m_1 E_1+m_2 E_2+m_3 E_3)\bigg)\bigg|,
\end{align*}
where $0\le m_1\le78$, $0\le m_2\le78$, and $0\le m_3\le156$ in the sum, and
\begin{align*}
\|Z_n-Z_n^*\|_1&=\sum_{m_1+m_2+m_3=n}\frac{\binom{78}{m_1}\binom{78}{m_2}\binom{156}{m_3}}{\binom{312}{n}}\bigg|E(78-m_1,78-m_2,156-m_3)\\
&\qquad\qquad\qquad\qquad{}-\bigg(\mu+\frac{\gamma}{312-n}(m_1 J_1+m_2 J_2+m_3 J_3)\bigg)\bigg|,
\end{align*}
where the sum is constrained in the same way,
with partial results appearing in Table~\ref{L1-distance}.  By definition, $\widehat Z_1=Z_1$.  We can regard $\|Z_n-\widehat Z_n\|_1$ as a measurement of the lack of linearity of the player's conditional expectation under basic strategy when $n$ cards have been seen.  It increases gradually as cards are dealt and then more sharply near the end of the shoe.  Replacing the EoRs of $\widehat Z_n$ by the level-six point count has only a small effect, whereas the use of the rather crude level-one point count has a rather substantial effect.  

\begin{table}[htb]
\caption{\label{L1-distance}$L^1$ distances between $Z_n$ (exact player conditional expectation under 39-deck composition-dependent basic strategy), $\widehat Z_n$ (approximate player conditional expectation based on EoRs), and $Z_n^*$ (approximate player conditional expectation based on a card-counting system).\medskip}
\catcode`#=\active \def#{\hphantom{$-$}}
\catcode`@=\active \def@{\hphantom{0}}
\tabcolsep=2mm
\renewcommand{\arraystretch}{1.}
\begin{center}
\begin{tabular}{cccccc}
\hline
\noalign{\smallskip}
    &       (a)                &      (b)          & \% incr. &    (c)            & \% incr.\\
$n$ & $\|Z_n-\widehat Z_n\|_1$ & $\|Z_n-Z_n^*\|_1$ & of (b)   & $\|Z_n-Z_n^*\|_1$ & of (c) \\
@seen@ &                       &  for $(-4,6,-1)$  & over (a) & for $(-1,1,0)$    & over (a)    \\                               
\noalign{\smallskip}\hline
\noalign{\smallskip}
@@1 & 0           & 0.00001667 &   --  & 0.0003582 &  --  \\
@@2 & 0.000006083 & 0.00001886 & 210.1@ & 0.0003614 & 5842. \\
@@3 & 0.00001110@ & 0.00002510 & 126.1@ & 0.0005390 & 4755. \\
@@4 & 0.00001631@ & 0.00002843 & 74.27@ & 0.0005442 & 3236. \\
\noalign{\medskip}
@26 & 0.0001120@@ & 0.0001243@ & 11.05@ & 0.001515@ & 1253. \\
@52 & 0.0002480@@ & 0.0002607@ & 5.121@ & 0.002257@ & 810.1 \\
@78 & 0.0004141@@ & 0.0004275@ & 3.245@ & 0.002919@ & 605.0 \\
104 & 0.0006225@@ & 0.0006364@ & 2.233@ & 0.003579@ & 474.9 \\
130 & 0.0008911@@ & 0.0009063@ & 1.709@ & 0.004282@ & 380.5 \\
156 & 0.001249@@@ & 0.001265@@ & 1.332@ & 0.005072@ & 306.1 \\
182 & 0.001754@@@ & 0.001771@@ & 0.9314 & 0.006013@ & 242.8 \\
208 & 0.002514@@@ & 0.002536@@ & 0.8654 & 0.007220@ & 187.1 \\
234 & 0.003788@@@ & 0.003821@@ & 0.8801 & 0.008960@ & 136.5 \\
260 & 0.006394@@@ & 0.006433@@ & 0.6069 & 0.01197@@ & 87.14 \\
286 & 0.01457@@@@ & 0.01466@@@ & 0.6154 & 0.02013@@ & 38.22 \\
\noalign{\medskip}
301 & 0.04268@@@@ & 0.04286@@@ & 0.4175 & 0.04682@@ & 9.704 \\
302 & 0.04954@@@@ & 0.04973@@@ & 0.3965 & 0.05239@@ & 5.753 \\
303 & 0.05917@@@@ & 0.05912@@@ &  --    & 0.06177@@ & 4.409 \\
304 & 0.06990@@@@ & 0.07004@@@ & 0.1882 & 0.07335@@ & 4.929 \\
\noalign{\smallskip}\hline
\end{tabular}
\end{center}
\end{table}

Next, we return to a previously computed quantity.  We supposed that players are required to pay a commission on each hand equal to $\nu=1/7$ of the initial amount bet, which would make snackjack a subfair game for the basic strategist.  Then $\E[(Z_n-\nu)^+]$ is the player's expected profit, assuming $n$ cards have been seen and assuming he bets one unit if and only if his net conditional expectation (taking the commission into account) is nonnegative.  The only problem is how does the player know whether his net conditional expectation is nonnegative?  Unless he has an electronic device programmed to evaluate $E(78-M_1,78-M_2,156-M_3)$ (which would be illegal in Nevada), he does not.  The best he can do is estimate his conditional expectation using card counting.  The \textit{betting efficiency} of a card-counting system could then be defined in terms of how close to the ideal $\E[(Z_n-\nu)^+]$ one could come in practice.  This would be
\begin{align*}
&\E[(Z_n-\nu)\bm1\{Z_n^*-\nu\ge0\}]\\
&\;{}=\E[(Z_n-\nu)\bm1\{\mu+(\gamma/52)\text{TC}_n\ge\nu\}]\\
&\;{}=\sum_{m_1+m_2+m_3=n}\frac{\binom{78}{m_1}\binom{78}{m_2}\binom{156}{m_3}}{\binom{312}{n}}[E(78-m_1,78-m_2,156-m_3)-\nu]\\
\noalign{\vglue-3mm}
&\qquad\qquad\qquad\qquad\qquad\qquad\qquad{}\cdot \bm1\{\text{TC}_n\ge52(\nu-\mu)/\gamma\},\quad
\end{align*}
where $0\le m_1\le78$, $0\le m_2\le78$, and $0\le m_3\le156$ in the sum.  Thus, the ratio
$$
\text{BE}_n=\frac{\E[(Z_n-\nu)\bm1\{Z_n^*-\nu\ge0\}]}{\E[(Z_n-\nu)^+]}
$$
is the \textit{betting efficiency} when $n$ cards have been seen.  We note that the threshold for the true count to suggest a positive expectation is $52(\nu-\mu)/\gamma\approx0.315367$ in the level-one system.  Partial results are shown in Table~\ref{efficiency}.  $\text{BE}_n$ is undefined if $n=1$ and is 1 (i.e., 100\%) for $n=2,3,4$ for the level-one system and for $n=2,3,\ldots,17$ except $n=10$ and $n=15$ for the level-six system.  

\begin{table}[htb]
\caption{\label{efficiency}The betting efficiency of two card-counting systems at 39-deck snackjack as a function of the number of cards seen.}
\catcode`#=\active \def#{\hphantom{$-$}}
\catcode`@=\active \def@{\hphantom{0}}
\tabcolsep=1.3mm
\renewcommand{\arraystretch}{1.}
\begin{center}
\begin{tabular}{ccccccc}
\hline
\noalign{\smallskip}
$n$ & $\text{BE}_n$ of & $\text{BE}_n$ of && $n$ & $\text{BE}_n$ of & $\text{BE}_n$ of\\
seen &  $(-4,6,-1)$    & $(-1,1,0)$ && seen & $(-4,6,-1)$    & $(-1,1,0)$       \\                       
\noalign{\smallskip}\hline
\noalign{\smallskip}
@26 & 0.999989 & 0.939475 && 234 & 0.992589 & 0.957591 \\
@52 & 0.998986 & 0.948355 && 260 & 0.991316 & 0.956100 \\
@78 & 0.999497 & 0.950998 && 286 & 0.986838 & 0.942744 \\
104 & 0.998716 & 0.951662 && \\
130 & 0.998168 & 0.951270 && 301 & 0.950243 & 0.895330 \\
156 & 0.998069 & 0.954986 && 302 & 0.953991 & 0.892247 \\
182 & 0.997312 & 0.956942 && 303 & 0.954517 & 0.889094 \\
208 & 0.995350 & 0.957795 && 304 & 0.892340 & 0.862632 \\
\noalign{\smallskip}\hline
\end{tabular}
\end{center}
\end{table}

It is useful to have a single number that can be called the \textit{betting efficiency} of a card-counting system.  For this we use an average of the quantities $\text{BE}_n$.  Since it is likely that the last one-quarter of the shoe is not dealt, we exclude decisions based on 234 or more cards.  This leads to
$$
\text{BE}:=\frac{1}{232}\sum_{n=2}^{233}\text{BE}_n.
$$
We find that $\text{BE}\approx0.9982$ for the level-six system $(-4,6,-1)$, and $\text{BE}\approx0.9508$ for the level-one system $(-1,1,0)$.  These numbers are not far from the correlations, 0.9999 and 0.9656, between the EoRs and the numbers of the point count.  Griffin~\cite[Chapter~4]{G99} used this correlation as a proxy for betting efficiency, unable to compute for blackjack numbers analogous to those in Table~\ref{efficiency} other than by computer simulation.

\begin{figure}[htb]
\centering
\includegraphics[width=2.25in]{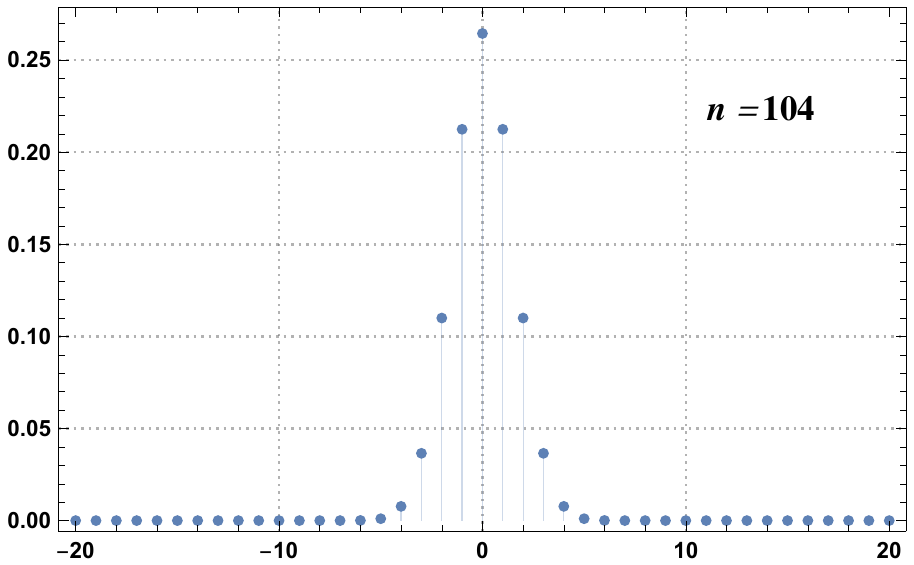}\quad
\includegraphics[width=2.25in]{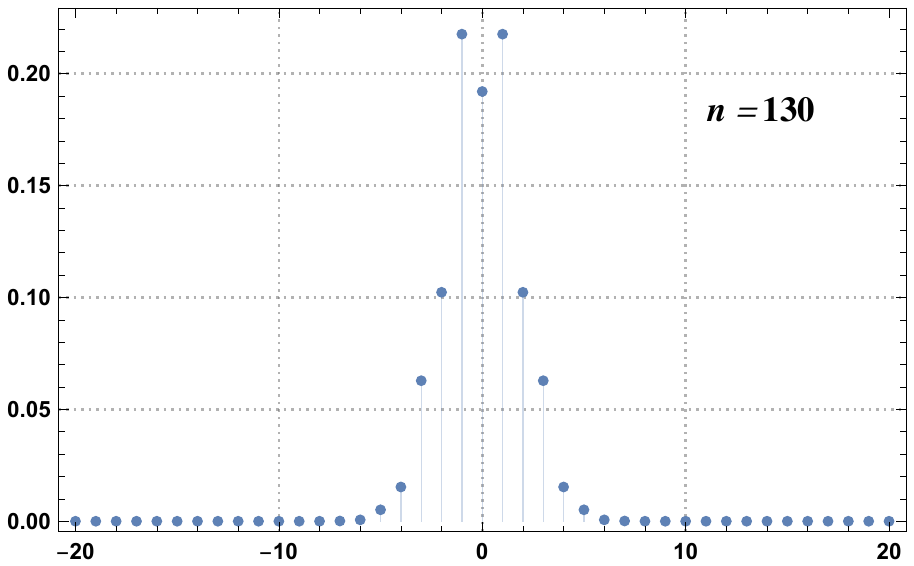}\medskip\par
\includegraphics[width=2.25in]{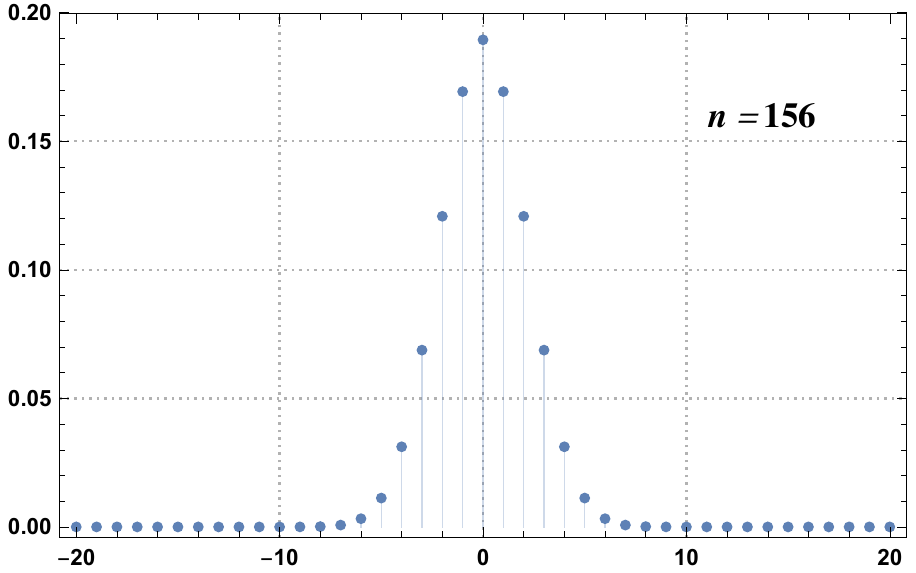}\quad
\includegraphics[width=2.25in]{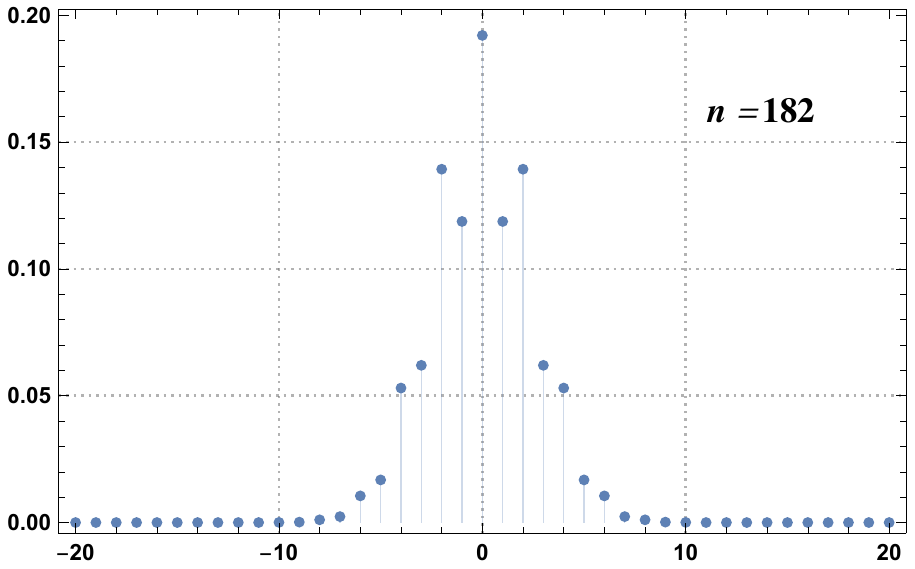}\medskip\par
\includegraphics[width=2.25in]{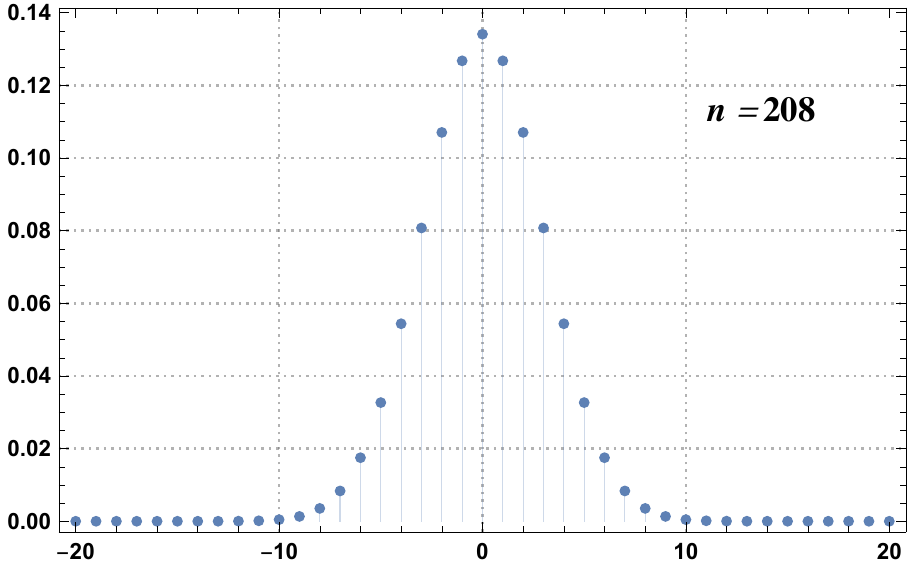}\quad
\includegraphics[width=2.25in]{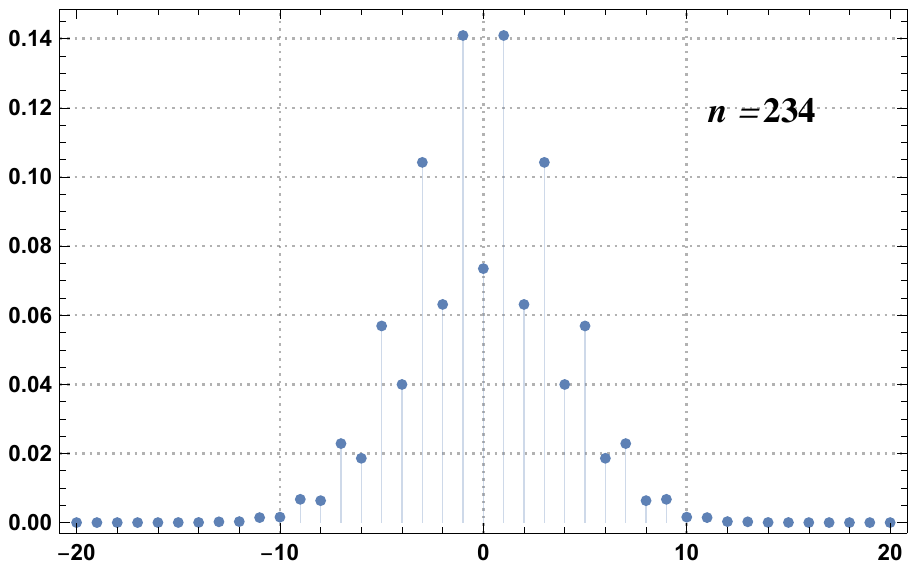}\medskip\par
\includegraphics[width=2.25in]{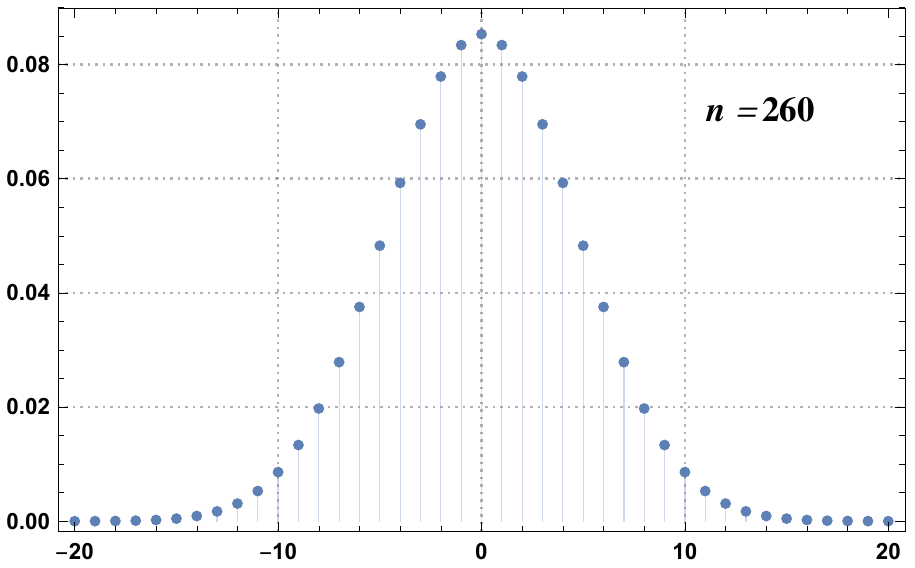}\quad
\includegraphics[width=2.25in]{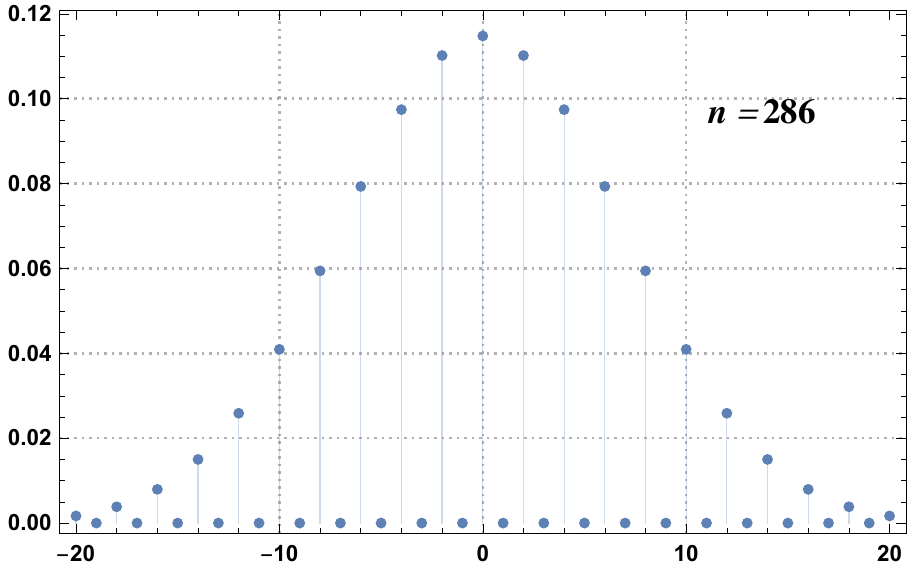}
\caption{\label{count-dist}At 39-deck snackjack, the distribution of $\text{TC}_n$ rounded to the nearest integer (assuming the deuces-minus-aces count), with $n$ being the number of cards seen.  Notice that the distribution is normal-like for $n$ an even multiple of 26 (left column) but not for $n$ an odd multiple of 26 (right column).}
\end{figure}

Now let us examine the level-one counting system, which we call the \textit{deuces-minus-aces} system, in more detail.  It is snackjack's analogue of the Hi-Lo system at blackjack.  The true count, when $n$ cards have been seen, including $M_1$ aces, $M_2$ deuces, and $M_3$ treys, is given by
$$
\text{TC}_n:=\frac{52(M_2-M_1)}{312-n},
$$
and the \textit{rounded true count} is $\text{TC}_n$ rounded to the nearest integer, denoted by $[\text{TC}_n]$.  More precisely, if $k-1/2<\text{TC}_n<k+1/2$, we define $[\text{TC}_n]:=k$  and, if $\text{TC}_n=k+1/2$, then $[\text{TC}_n]=k$ with probability 1/2 and $[\text{TC}_n]=k+1$ with probability 1/2.  This symmetric rounding ensures that the distribution of $[\text{TC}_n]$ is symmetric about 0.  Indeed,
\begin{align*}
\P([\text{TC}_n]=k)&=\sum_{m_1+m_2+m_3=n}\frac{\binom{78}{m_1}\binom{78}{m_2}\binom{156}{m_3}}{\binom{312}{n}}\\
&\qquad\qquad\qquad\qquad\cdot\bigg[\bm1\bigg\{k-\frac12<\frac{52(m_2-m_1)}{312-n}<k+\frac12\bigg\}\\
&\qquad\qquad\qquad\qquad\quad{}+\frac12\,\bm1\bigg\{\frac{52(m_2-m_1)}{312-n}=k-\frac12\text{ or }k+\frac12\bigg\}\bigg],
\end{align*}
where the sum is constrained by $0\le m_1\le78$, $0\le m_2\le78$, and $0\le m_3\le156$.  Figure~\ref{count-dist} plots the graph of $[\text{TC}_n]$ for $n=26m$, $m=4,5,\ldots,11$.  

\afterpage{\clearpage}

Next, we evaluate the conditional expectation at snackjack (assuming the commission of $\nu=1/7$), given the rounded true count.  This is
\begin{align*}
&\!\!\!\!\!\E[Z_n-\nu\mid[\text{TC}_n]=k]\\
&=\sum_{m_1+m_2+m_3=n}\frac{\binom{78}{m_1}\binom{78}{m_2}\binom{156}{m_3}}{\binom{312}{n}}[E(78-m_1,78-m_2,156-m_3)-\nu]\\
&\qquad\qquad\qquad\cdot\bigg[\bm1\bigg\{k-\frac12<\frac{52(m_2-m_1)}{312-n}<k+\frac12\bigg\}\\
&\qquad\qquad\qquad\quad{}+\frac12\,\bm1\bigg\{\frac{52(m_2-m_1)}{312-n}=k-\frac12\text{ or }k+\frac12\bigg\}\bigg]\bigg/\P([\text{TC}_n]=k),
\end{align*}
with the same constraints on the sum, and results are tabulated in Table~\ref{condl-EVs} for $n=78$, 156, and 234.  We find that the rounded true count is a good estimate of the player's expectation in percentage terms. 

We now consider a betting strategy similar to the one assumed by Schlesinger \cite[Chapter 5]{S18}, but a little simpler.  We assume that the player bets 
$$
\max(1,\min([\text{TC}_n],6)).
$$
That is, the player bets the rounded true count, but never less than one unit or more than six units.  Thus, this betting strategy has a 6 to 1 spread.  It could be argued that the bettor should walk away if the true count falls below some threshold, but we assume that he continues to play and bet one unit, perhaps to disguise his status as a card counter.

\begin{figure}[htb]
\centering
\includegraphics[width=3.in]{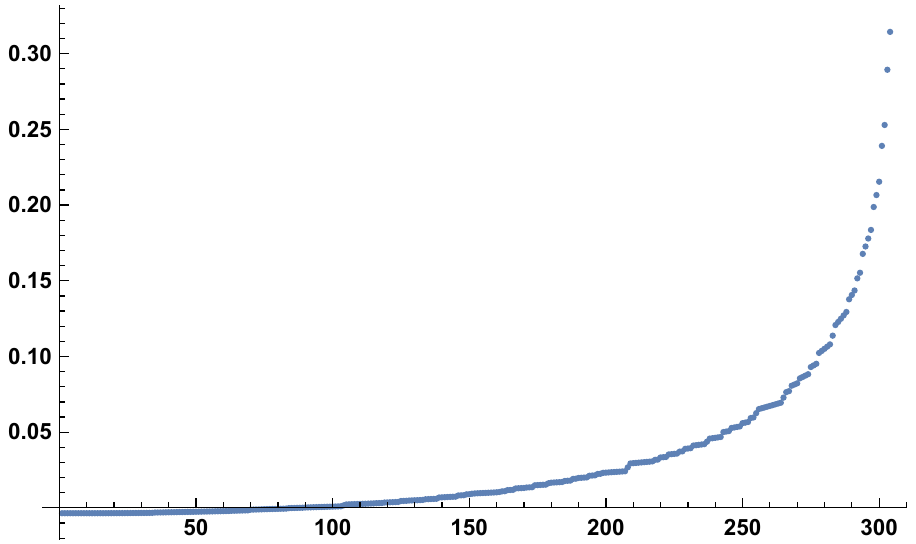}
\caption{\label{gain-bet-var}Snackjack expectation as a function of the number $n$ of cards seen, assuming a commission of $\nu=1/7$ per unit initially bet and bets equal to the rounded true count, but always at least one unit and at most six units.}
\end{figure}

We can then evaluate the player's expected profit at each level of penetration.  The formula is
\begin{align*}
&\!\!\!\!\!\E[\max(1,\min([\text{TC}_n],6))(Z_n-\nu)]\\
&=\sum_{m_1+m_2+m_3=n}\frac{\binom{78}{m_1}\binom{78}{m_2}\binom{156}{m_3}}{\binom{312}{n}}[E(78-m_1,78-m_2,156-m_3)-\nu]\\
&\qquad\qquad\cdot\sum_k\max(1,\min(k,6))\bigg[\bm1\bigg\{k-\frac12<\frac{52(m_2-m_1)}{312-n}<k+\frac12\bigg\}\\
&\qquad\qquad\qquad\qquad\qquad\qquad\quad{}+\frac12\,\bm1\bigg\{\frac{52(m_2-m_1)}{312-n}=k-\frac12\text{ or }k+\frac12\bigg\}\bigg],
\end{align*}
with the same constraints on the outer sum, and the results are plotted in Figure~\ref{gain-bet-var}.

\begin{table}[htb]
\caption{\label{condl-EVs}The conditional expectation at snackjack (assuming a commission of $\nu=1/7$ per unit initially bet), given the (rounded) true count.}
\catcode`@=\active \def@{\hphantom{0}}
\catcode`#=\active \def#{\hphantom{$-$}}
\tabcolsep=.2cm
\begin{center}
\begin{tabular}{ccccccc}
&\multicolumn{2}{c}{$n=78$}&\multicolumn{2}{c}{$n=156$}&\multicolumn{2}{c}{$n=234$}\\
\noalign{\smallskip}
\hline
\noalign{\smallskip}
$[\text{TC}_n]$ & cond'l ex & probab & cond'l ex & probab & cond'l ex & probab \\
\noalign{\smallskip}
\hline
\noalign{\smallskip}
$-6$ & $-0.0688$@ & 0.00000221 & $-0.0724$@ &  0.00320 & $-0.0742$@ &  0.0186 \\
$-5$ & $-0.0585$@ & 0.0000640@ & $-0.0606$@ &  0.0112@ & $-0.0614$@ &  0.0569 \\
$-4$ & $-0.0463$@ & 0.00195@@@ & $-0.0489$@ &  0.0312@ & $-0.0497$@ &  0.0400 \\
$-3$ & $-0.0361$@ & 0.0146@@@@ & $-0.0373$@ &  0.0688@ & $-0.0372$@ &  0.104@ \\
$-2$ & $-0.0246$@ & 0.0984@@@@ & $-0.0257$@ &  0.121@@ & $-0.0254$@ &  0.0631 \\
$-1$ & $-0.0143$@ & 0.207@@@@@ & $-0.0143$@ &  0.169@@ & $-0.0134$@ &  0.141@ \\
#0   & $-0.00352$ & 0.355@@@@@ & $-0.00311$ &  0.189@@ & $-0.00183$ &  0.0735 \\
#1   & #$0.00723$ & 0.207@@@@@ & #$0.00791$ &  0.169@@ & #$0.00950$ &  0.141@ \\
#2   & #$0.0175$@ & 0.0984@@@@ & #$0.0187$@ &  0.121@@ & #$0.0207$@ &  0.0631 \\
#3   & #$0.0288$@ & 0.0146@@@@ & #$0.0292$@ &  0.0688@ & #$0.0312$@ &  0.104@ \\
#4   & #$0.0387$@ & 0.00195@@@ & #$0.0395$@ &  0.0312@ & #$0.0418$@ &  0.0400 \\
#5   & #$0.0506$@ & 0.0000640@ & #$0.0494$@ &  0.0112@ & #$0.0513$@ &  0.0569 \\
#6   & #$0.0606$@ & 0.00000221 & #$0.0590$@ &  0.00320 & #$0.0611$@ &  0.0186\\
\noalign{\smallskip}
\hline
\end{tabular}
\end{center}
\end{table} 

The average expected value over the first 3/4 of the shoe ($0\le n\le 233$) is 0.00779463.
The average over the first 5/6 of the shoe ($0\le n\le 259$) is 0.0123218.

\section{Card counting and strategy variation}\label{strategy variation}

In this section we show that card counting can be used to determine when a departure from basic strategy is called for.  A table analogous to Table~\ref{SJ-BS-1deck} for 39-deck snackjack could be generated, but it would have 81 rows.  For simplicity, we consider only two-card player hands, which leads to the 18-row Table~\ref{SJ-BS-39deck}.  The omitted rows are quite similar to the included rows with the same hard or soft total and the same dealer upcard.  It appears that $(0,0,2)$, $(0,1,1)$, and $(1,1,0)$ have the greatest potential for profitable strategy variation, but this must be quantified in a systematic way.

\begin{table}[htb]
\caption{\label{SJ-BS-39deck}The analogue of Table~\ref{SJ-BS-1deck} for 39-deck snackjack, but restricted to two-card player hands for simplicity.}
\catcode`@=\active \def@{\hphantom{0}}
\catcode`#=\active \def#{\hphantom{$-$}}
\tabcolsep=.12cm
\begin{center}
\begin{tabular}{ccccccccc}
\noalign{\smallskip}
\hline
\noalign{\smallskip}
nos of & htot & tot & up & $E_{\rm std}$ & $E_{\rm hit}$ & $E_{\rm dbl}$ & $E_{\rm spl}$ & bs \\
1s, 2s, 3s & &&&&&& \\
\noalign{\smallskip}
\hline
\noalign{\smallskip}
$(0,0,2)$   & 6 & h6 & 1 & $-0.065126$ & $-0.550044$ & $-1.100088$ & $-0.040216$ & Spl \\
$(0,0,2)$   & 6 & h6 & 2 & #$0.014416$ & $-0.562389$ & $-1.124778$ & #$0.164520$ & Spl \\
$(0,0,2)$   & 6 & h6 & 3 & #$0.083319$ & $-0.514954$ & $-1.029908$ & #$0.209534$ & Spl \\
\noalign{\medskip}
$(0,1,1)$   & 5 & h5 & 1 & $-0.747557$ & $-0.316759$ & $-0.633519$ & na & H  \\
$(0,1,1)$   & 5 & h5 & 2 & $-0.435486$ & $-0.314908$ & $-0.629816$ & na & H  \\
$(0,1,1)$   & 5 & h5 & 3 & $-0.666667$ & $-0.247564$ & $-0.495128$ & na & H  \\
\noalign{\medskip}
$(0,2,0)$   & 4 & h4 & 1 & $-0.742607$ & #$0.317544$ & #$0.426296$ & $-1.145240$ & D \\
$(0,2,0)$   & 4 & h4 & 2 & $-0.429399$ & #$0.299155$ & #$0.544761$ & $-0.631158$ & D \\
$(0,2,0)$   & 4 & h4 & 3 & $-0.668859$ & #$0.417413$ & #$0.625288$ & $-0.956068$ & D \\
\noalign{\medskip}
$(1,0,1)$   & 4 & s7 & 1 & #$1.5@@@@@$ & #$0.315713$ & #$0.417682$ & na & S  \\
$(1,0,1)$   & 4 & s7 & 2 & #$1.5@@@@@$ & #$0.294149$ & #$0.530176$ & na & S  \\
$(1,0,1)$   & 4 & s7 & 3 & #$1.5@@@@@$ & #$0.415415$ & #$0.622725$ & na & S  \\
\noalign{\medskip}
$(1,1,0)$   & 3 & s6 & 1 & $-0.060261$ & #$0.088154$ & $-0.032482$ & na & H  \\
$(1,1,0)$   & 3 & s6 & 2 & #$0.021141$ & #$0.114036$ & #$0.175894$ & na & D  \\
$(1,1,0)$   & 3 & s6 & 3 & #$0.085129$ & #$0.208650$ & #$0.211340$ & na & D  \\
\noalign{\medskip}
$(2,0,0)$   & 2 & s5 & 1 & $-0.749293$ & #$0.064931$ & $-0.381205$ & #$0.430233$ & Spl \\
$(2,0,0)$   & 2 & s5 & 2 & $-0.429377$ & #$0.049446$ & $-0.062008$ & #$0.540126$ & Spl \\
$(2,0,0)$   & 2 & s5 & 3 & $-0.663062$ & #$0.157834$ & $-0.159281$ & #$0.634053$ & Spl \\
\noalign{\smallskip}
\hline
\end{tabular}
\end{center}
\end{table}

Let us first treat the case of $(0,0,2)$ (i.e., a pair of treys), which the 39-deck basic strategist splits against any dealer upcard, but with which standing may be preferable in some situations.  

We denote $E_{\std}((0,0,2),u)$ of Section~\ref{BS methodology} by $E_{\std,(n_1,n_2,n_3)}((0,0,2),u)$
with $(n_1,n_2,n_3)$ indicating the post-deal shoe composition (i.e., the hand's two 3s and dealer upcard $u$ are excluded from the unseen shoe).  A similar interpretation applies to $E_{\spl,(n_1,n_2,n_3)}((0,0,2),u)$.  

The difference $E_{\std,(n_1,n_2,n_3)}((0,0,2),u)-E_{\spl,(n_1,n_2,n_3)}((0,0,2),u)$, which represents the expected gain by departing from basic strategy, is well defined under the following conditions (see Appendix~D for explicit formulas): $n_1,n_2,n_3\ge0$ and
\begin{itemize}
\item if $u=1$, then $n_1+n_2+n_3\ge4$, $n_1+n_2\ge1$;
\item if $u=2$, then $n_1+n_2+n_3\ge4$;
\item if $u=3$, then $n_1+n_2+n_3\ge3$, $n_2+n_3\ge1$.
\end{itemize}
We assume in fact that $n_1+n_2+n_3\ge5$ for all $u$ because a new hand should never be dealt with fewer than eight cards remaining.

With $(0,0,2)$ vs.~1, the proportion of shoe compositions that call for a departure from basic strategy is 
$439{,}742/954{,}925\approx0.460499$.  With $(0,0,2)$ vs.~2, the proportion is $271{,}854/955{,}075\approx0.284642$.   With $(0,0,2)$ vs.~3, the proportion is $358{,}973/961{,}005\approx0.373539$.  

With $(0,0,2)$ vs.~1, the probability that a departure from basic strategy is called for when $n$ cards have been seen (before the hand is dealt) is 
\begin{equation}\label{weightedave1}
\frac{\sum_{(m_1,m_2,m_3)\in\Gamma_n((0,0,2),1)}\binom{78}{m_1}\binom{78}{m_2}\binom{156}{m_3}\alpha_1(m_1,m_2,m_3)}{\sum_{(m_1,m_2,m_3)\in\Gamma_n((0,0,2),1)}\binom{78}{m_1}\binom{78}{m_2}\binom{156}{m_3}}
\end{equation}
for $n=1,2,\ldots,304$, where

\begin{align*}
\Gamma_n((0,0,2),1)&:=\{(m_1,m_2,m_3): 0\le m_1\le77,\, 0\le m_2\le78,\, 0\le m_3\le154,\nonumber\\
&\qquad\qquad\qquad\qquad\quad{}m_1+m_2\le154,\,m_1+m_2+m_3=n\}
\end{align*}
and
\begin{equation}\label{alpha_1}
\alpha_1(m_1,m_2,m_3):=\begin{cases}1&\text{if }E_{\std,(77-m_1,78-m_2,154-m_3)}((0,0,2),1)\\
& \quad{}>E_{\spl,(77-m_1,78-m_2,154-m_3)}((0,0,2),1),\\
0&\text{otherwise}.\end{cases}
\end{equation}
(The condition $m_1+m_2\le154$ ensures that there are enough 1s and 2s remaining to allow the dealer's downcard to be other than a 3.)

The corresponding probability for $(0,0,2)$ vs.~2 is 
\begin{equation}\label{weightedave2}
\frac{\sum_{(m_1,m_2,m_3)\in\Gamma_n((0,0,2),2)}\binom{78}{m_1}\binom{78}{m_2}\binom{156}{m_3}\alpha_2(m_1,m_2,m_3)}{\sum_{(m_1,m_2,m_3)\in\Gamma_n((0,0,2),2)}\binom{78}{m_1}\binom{78}{m_2}\binom{156}{m_3}}
\end{equation}
for $n=1,2,\ldots,304$, where
\begin{align*}
\Gamma_n((0,0,2),2)&:=\{(m_1,m_2,m_3): 0\le m_1\le78,\, 0\le m_2\le77,\, 0\le m_3\le154,\\
&\qquad\qquad\qquad\qquad\quad{}m_1+m_2+m_3=n\}
\end{align*}
and
\begin{equation}\label{alpha_2}
\alpha_2(m_1,m_2,m_3):=\begin{cases}1&\text{if }E_{\std,(78-m_1,77-m_2,154-m_3)}((0,0,2),2)\\
& \quad{}>E_{\spl,(78-m_1,77-m_2,154-m_3)}((0,0,2),2),\\
0&\text{otherwise}.\end{cases}
\end{equation}

\begin{figure}[htb]
\centering
\includegraphics[width=2.25in]{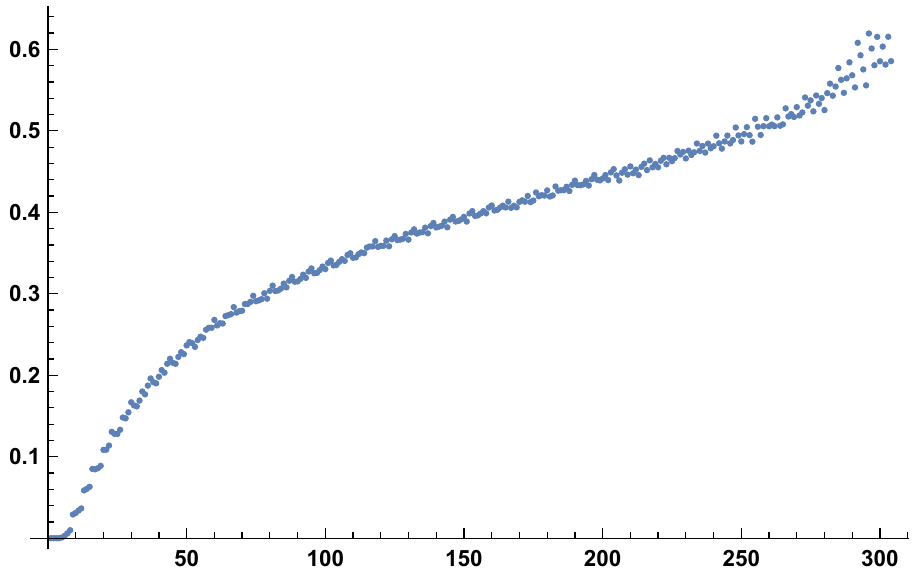}\quad
\includegraphics[width=2.25in]{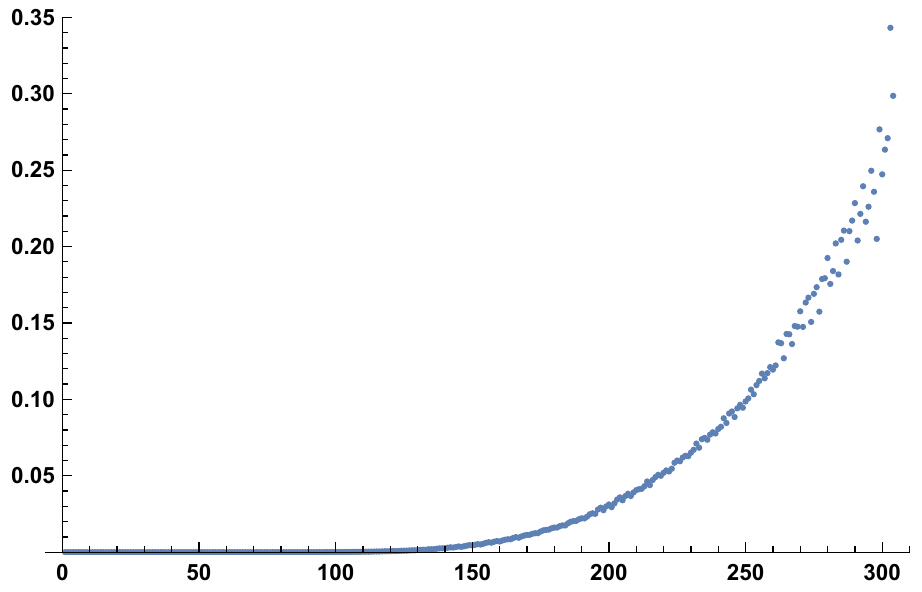}\vglue3mm
\includegraphics[width=2.25in]{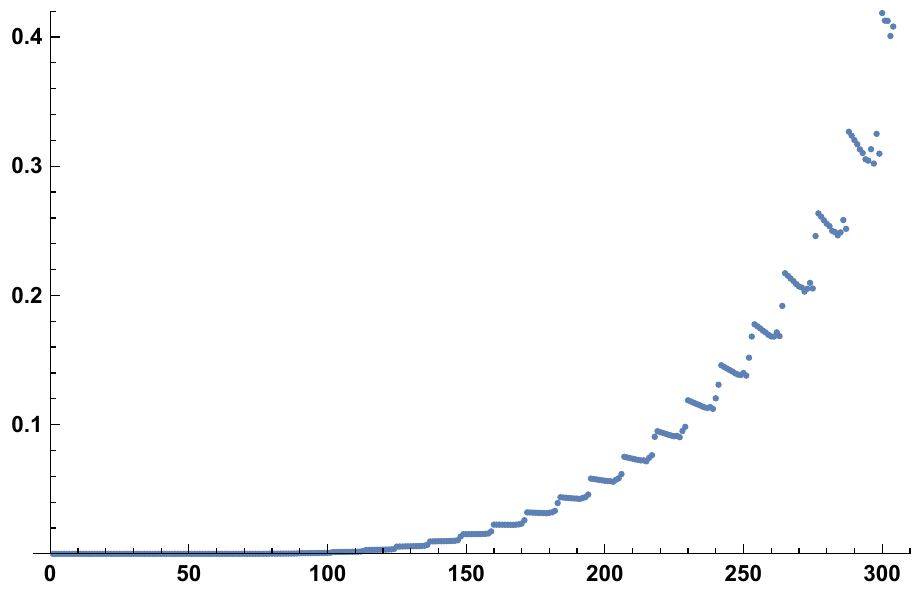}
\caption{\label{prob-33v123}With $(0,0,2)$ vs.~1 (top left), $(0,0,2)$ vs.~2 (top right), and $(0,0,2)$ vs.~3 (bottom), the probability that a departure from basic strategy is called for when $n$ cards have been seen (before the hand is dealt), as a function of $n$, $1\le n\le304$.}
\end{figure}

The corresponding probability for $(0,0,2)$ vs.~3 is 
\begin{equation}\label{weightedave3}
\frac{\sum_{(m_1,m_2,m_3)\in\Gamma_n((0,0,2),3)}\binom{78}{m_1}\binom{78}{m_2}\binom{156}{m_3}\alpha_3(m_1,m_2,m_3)}{\sum_{(m_1,m_2,m_3)\in\Gamma_n((0,0,2),3)}\binom{78}{m_1}\binom{78}{m_2}\binom{156}{m_3}}
\end{equation}
for $n=1,2,\ldots,304$, where
\begin{align*}
\Gamma_n((0,0,2),3)&:=\{(m_1,m_2,m_3): 0\le m_1\le78,\, 0\le m_2\le78,\, 0\le m_3\le153,\\
&\qquad\qquad\qquad\qquad\quad{}m_2+m_3\le230,\,m_1+m_2+m_3=n\}
\end{align*}
and
\begin{equation}\label{alpha_3}
\alpha_3(m_1,m_2,m_3):=\begin{cases}1&\text{if }E_{\std,(78-m_1,78-m_2,153-m_3)}((0,0,2),3)\\
& \quad{}>E_{\spl,(78-m_1,78-m_2,153-m_3)}((0,0,2),3),\\
0&\text{otherwise}.\end{cases}
\end{equation}
The expressions \eqref{weightedave1}, \eqref{weightedave2}, and \eqref{weightedave3} are graphed in Figure~\ref{prob-33v123}.

\begin{figure}[htb]
\centering
\includegraphics[width=2.25in]{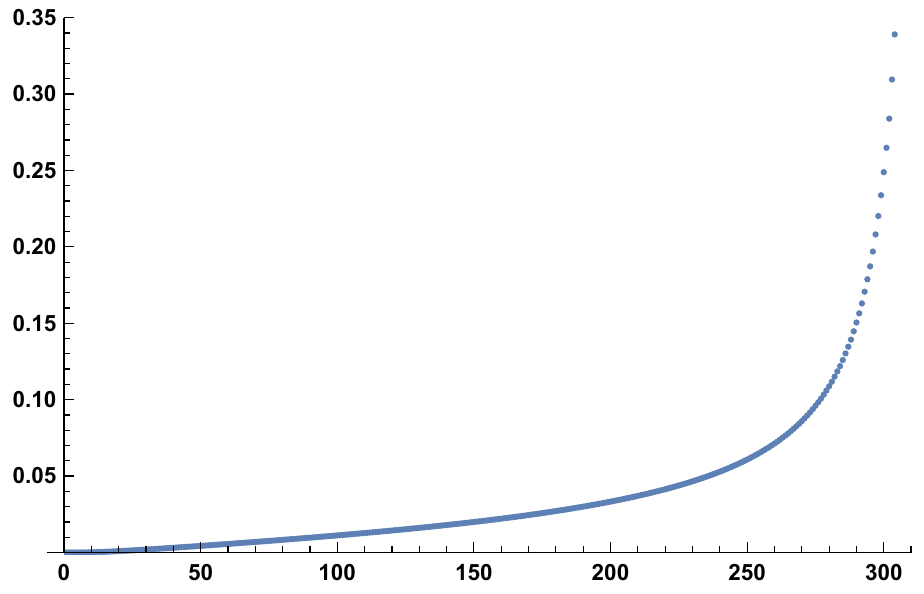}\quad
\includegraphics[width=2.25in]{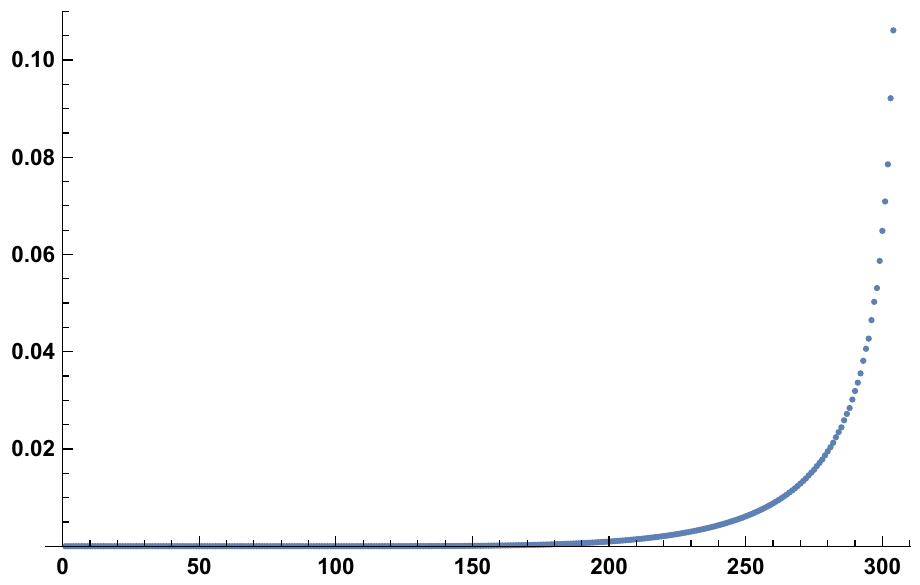}\vglue3mm
\includegraphics[width=2.25in]{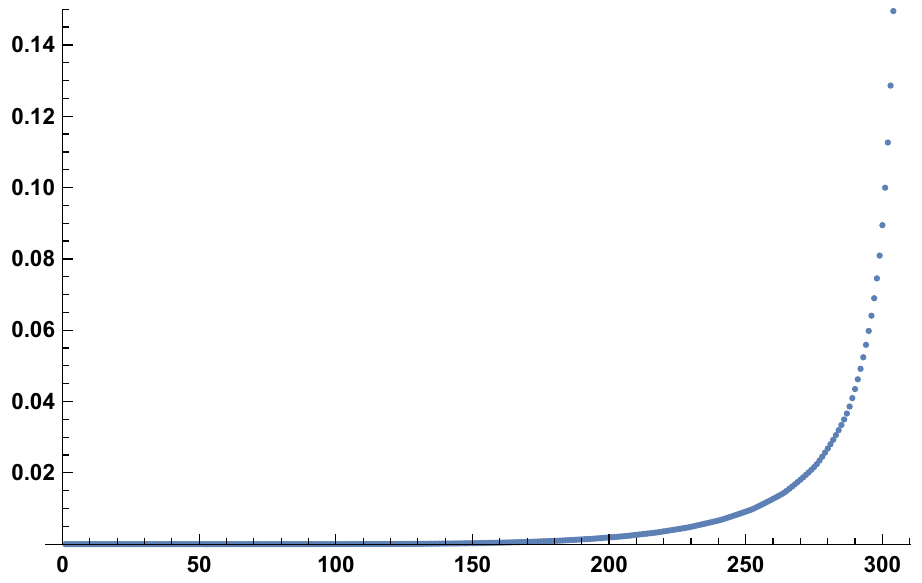}
\caption{\label{ev-33v123}With $(0,0,2)$ vs.~1 (left top), $(0,0,2)$ vs.~2 (right top), and $(0,0,2)$ vs.~3 (bottom), the additional expectation that a departure from basic strategy provides when $n$ cards have been seen (before the hand is dealt), as a function of $n$, $1\le n\le304$.  Notice that the vertical scales differ considerably.}
\end{figure}

More important than the probability that a departure from basic strategy is called for is the additional expectation that such a departure provides.  With $(0,0,2)$ vs.~1, 2, or 3 this is given by \eqref{weightedave1}, \eqref{weightedave2}, or \eqref{weightedave3} but with 
\begin{align}\label{alpha_1-alt}
\alpha_1(m_1,m_2,m_3)&:=[E_{\std,(77-m_1,78-m_2,154-m_3)}((0,0,2),1)\nonumber\\
& \qquad{}-E_{\spl,(77-m_1,78-m_2,154-m_3)}((0,0,2),1)]^+,
\end{align}
\begin{align}\label{alpha_2-alt}
\alpha_2(m_1,m_2,m_3)&:=[E_{\std,(78-m_1,77-m_2,154-m_3)}((0,0,2),2)\nonumber\\
& \qquad{}-E_{\spl,(78-m_1,77-m_2,154-m_3)}((0,0,2),2)]^+,
\end{align}
or
\begin{align}\label{alpha_3-alt}
\alpha_3(m_1,m_2,m_3)&:=[E_{\std,(78-m_1,78-m_2,153-m_3)}((0,0,2),3)\nonumber\\
& \qquad{}-E_{\spl,(78-m_1,78-m_2,153-m_3)}((0,0,2),3)]^+
\end{align}
in place of \eqref{alpha_1}, \eqref{alpha_2}, or \eqref{alpha_3}.  The expressions \eqref{weightedave1},  \eqref{weightedave2}, and \eqref{weightedave3}, using \eqref{alpha_1-alt}, \eqref{alpha_2-alt}, and \eqref{alpha_3-alt}, are graphed in Figure~\ref{ev-33v123}.

As with bet variation, the only way to recognize potentially profitable departures from basic strategy is with card counting.  First, we analyze the case $(0,0,2)$ vs.~1, which is complicated by the assumption that the dealer does not have a natural.  The effects of removal are
\begin{align}\label{EoR-33v1}
\text{EoR}(i)&:=E_{\std,(77,78,154)-\bm e_i}((0,0,2),1)-E_{\spl,(77,78,154)-\bm e_i}((0,0,2),1)\nonumber\\
&\qquad{}-[E_{\std,(77,78,154)}((0,0,2),1)-E_{\spl,(77,78,154)}((0,0,2),1)]
\end{align}
for $i=1,2,3$.
The numbers \eqref{EoR-33v1}, multiplied by 308, are
\begin{equation*}
E_1=\frac{78{,}498{,}676}{49{,}345{,}645},\quad E_2=-\frac{895{,}474{,}426}{444{,}110{,}805},\quad E_3=\frac{33{,}220{,}264}{148{,}036{,}935},
\end{equation*}
with decimal equivalents 1.59079, $-2.01633$, and 0.224405.  The analogue of \eqref{EoR constraint} is 
\begin{equation}\label{EoRs weighted}
w_1\,E_1+w_2\,E_2+w_3\,E_3=0
\end{equation}
with weights
\begin{equation}\label{weights-33v1}
w_1=\frac{77}{308}\,\frac{154}{155},\quad w_2=\frac{78}{308}\,\frac{154}{155},\quad w_3=\frac{154}{308};
\end{equation}
see Epstein~\cite[p.~244]{Ep67}.
The correlation between the effects of removal and the deuces-minus-aces counting system $(J_1,J_2,J_3)=(-1,1,0)$ is
\begin{equation*}
\rho=\frac{w_1\,E_1 J_1+w_2\,E_2 J_2+w_3\,E_3 J_3}{\sigma_E\,\sigma_J}\approx-0.985649,
\end{equation*}
where $\sigma_E^2:=w_1\,E_1^2+w_2\,E_2^2+w_3\,E_3^2$ and $\sigma_J^2:=w_1\,J_1^2+w_2\,J_2^2+w_3\,J_3^2-(w_1\,J_1+w_2\,J_2+w_3\,J_3)^2$, and the regression coefficient is
$$
\gamma=\frac{w_1\,E_1 J_1+w_2\,E_2 J_2+w_3\,E_3 J_3}{w_1\,J_1^2+w_2\,J_2^2+w_3\,J_3^2}=-\frac{41{,}415{,}529{,}232}{22{,}945{,}724{,}925}\approx-1.80493,
$$
which is the $\gamma$ that minimizes the sum of squares $w_1(E_1-\gamma J_1)^2+w_2(E_2-\gamma J_2)^2+w_3(E_3-\gamma J_3)^2$.  

The analogues of \eqref{Z_n hat} and \eqref{Z_n^*} can be found by observing from \eqref{EoRs weighted} and \eqref{weights-33v1} that
$$
77\,E_1+78\,E_2+154\,E_3=-E_3.
$$
Hence
\begin{align*}
\widehat Z_n&=\frac{1}{309-n}\sum_{j=n+1}^{309}(\mu-E_{X_j})\nonumber\\
&=\mu+\frac{1}{309-n}\bigg(E_3+\sum_{j=1}^n E_{X_j}\bigg)
\end{align*}
and 
\begin{align}\label{Z_n^*-33v1}
Z_n^*&=\mu+\frac{\gamma}{52}\bigg(\frac{52}{309-n}\bigg[J_3+\sum_{j=1}^n J_{X_j}\bigg]\bigg)\nonumber\\
&=\mu+\frac{\gamma}{52}\bigg(\frac{52}{309-n}\sum_{j=1}^n J_{X_j}\bigg)=\mu+\frac{\gamma}{52}\text{TC}_n^*,
\end{align}
where
\begin{align*}
\mu&=E_{\std,(77,78,154)}((0,0,2),1)-E_{\spl,(77,78,154)}((0,0,2),1)\\
&=-\frac{60{,}451}{2{,}426{,}835}\approx-0.0249094.
\end{align*}
This allows the card counter to know (approximately) when it is advantageous to depart from basic strategy when holding $(0,0,2)$ vs.~1.  Indeed, $Z_n^*>0$ is equivalent to 
$$
\text{TC}_n^*:=\frac{52(m_2-m_1)}{309-n}<-\frac{52\mu}{\gamma}
$$
(the inequality is reversed because $\gamma<0$ or, equivalently, $\rho<0$).  The fraction 
$$
-\frac{52\mu}{\gamma}=-\frac{7{,}430{,}334{,}665}{10{,}353{,}882{,}308}\approx-0.717638
$$
is the \textit{index number} for this departure.  In this case, if the \textit{adjusted true count} $\text{TC}_n^*$ is less than this index number, standing on $(0,0,2)$ vs.~1 is called for instead of splitting.  We say ``adjusted'' because the player's treys and the dealer's ace are excluded from the count.

In practice, we would round the index number to $-1$, and this play would occur relatively often (the rounded adjusted true count would have to be at $-1$ or less).  Of course we would be betting only one unit, and we can infer an upper bound on the profit potential from the first panel in Figure~\ref{ev-33v123}.  In fact, we can compute it precisely using \eqref{weightedave1} with
\begin{align*}
\alpha_1(m_1,m_2,m_3)&:=[E_{\std,(77-m_1,78-m_2,154-m_3)}((0,0,2),1)\\
&\qquad{}-E_{\spl,(77-m_1,78-m_2,154-m_3)}((0,0,2),1)]\\
&\qquad\qquad{}\cdot\bm1\bigg\{\bigg[\frac{52(m_2-m_1)}{309-n}\bigg]\le-1\bigg\}.
\end{align*}
The average of these expectations over $1\le n\le 233$ is approximately 0.0162143.

Next, we analyze the simpler case of $(0,0,2)$ vs.~2.  The effects of removal are
\begin{align}\label{EoR-33v2}
\text{EoR}(i)&:=E_{\std,(78,77,154)-\bm e_i}((0,0,2),2)-E_{\spl,(78,77,154)-\bm e_i}((0,0,2),2)\nonumber\\
&\qquad{}-[E_{\std,(78,77,154)}((0,0,2),2)-E_{\spl,(78,77,154)}((0,0,2),2)]
\end{align}
for $i=1,2,3$.
The numbers \eqref{EoR-33v2}, multiplied by 308, are
$$
E_1=\frac{318{,}420{,}487}{295{,}118{,}793},\quad E_2=-\frac{2{,}651{,}203{,}088}{1{,}475{,}593{,}965},\quad E_3=\frac{519{,}211{,}999}{1{,}475{,}593{,}965},
$$
with decimal equivalents 1.07896, $-1.79670$, and 0.351866.  The analogue of \eqref{EoR constraint} is \eqref{EoRs weighted} with weights
$$
w_1=\frac{78}{309},\quad w_2=\frac{77}{309},\quad w_3=\frac{154}{309}.
$$
The correlation between the effects of removal and the deuces-minus-aces counting system $(J_1,J_2,J_3)=(-1,1,0)$ is
$\rho\approx-0.943999$, and the regression coefficient is
$$
\gamma=-\frac{328{,}326{,}627{,}706}{228{,}717{,}064{,}575}\approx-1.43551.
$$

The analogues of \eqref{Z_n hat} and \eqref{Z_n^*} are
$$
\widehat Z_n=\mu+\frac{1}{309-n}\sum_{j=1}^n E_{X_j}
$$
and 
$$
Z_n^*=\mu+\frac{\gamma}{52}\bigg(\frac{52}{309-n}\sum_{j=1}^n J_{X_j}\bigg)=\mu+\frac{\gamma}{52}\text{TC}_n^*,
$$
where
\begin{align*}
\mu&=E_{\std,(78,77,154)}((0,0,2),2)-E_{\spl,(78,77,154)}((0,0,2),2)\\
&=-\frac{1{,}452{,}413}{9{,}676{,}026}\approx-0.150104.
\end{align*}
This analysis is similar to~\cite[pp.~668--670]{Et10} for $6,\text{T}$ vs.~9 in blackjack.

This allows the card counter to know (approximately) when it is advantageous to depart from basic strategy when holding $(0,0,2)$ vs.~2.  Indeed, $Z_n^*>0$ is equivalent to 
$$
\text{TC}_n^*:=\frac{52(m_2-m_1)}{309-n}<-\frac{52\mu}{\gamma}.
$$
The fraction 
$$
-\frac{52\mu}{\gamma}=-\frac{892{,}616{,}719{,}475}{164{,}163{,}313{,}853}\approx-5.43737
$$
is the index number for this departure from basic strategy.  In this case, if the adjusted true count $\text{TC}_n^*$ is less than this index number, standing on $(0,0,2)$ vs.~2 is called for instead of splitting. 

In practice, the index number would be rounded to $-6$.  This play would seldom occur, for the rounded adjusted true count would have to be at $-6$ or less (see Table~\ref{condl-EVs}), and when it did occur, the bet size would be one unit.  So an upper bound on the value of this departure from basic strategy can be inferred from the second panel in Figure~\ref{ev-33v123}.

Finally, we analyze the case $(0,0,2)$ vs.~3, which also involves the assumption that the dealer does not have a natural.  The effects of removal are
\begin{align}\label{EoR 33v3}
\text{EoR}(i)&:=E_{\std,(78,78,153)-\bm e_i}((0,0,2),3)-E_{\spl,(78,78,153)-\bm e_i}((0,0,2),3)\nonumber\\
&\qquad{}-[E_{\std,(78,78,153)}((0,0,2),3)-E_{\spl,(78,78,153)}((0,0,2),3)]
\end{align}
for $i=1,2,3$.
The numbers \eqref{EoR 33v3}, multiplied by 308, are
$$
E_1=\frac{5{,}425{,}240}{3{,}616{,}767},\quad E_2=-\frac{209{,}017{,}702}{138{,}642{,}735},\quad E_3=\frac{24{,}804}{46{,}214{,}245},
$$
with decimal equivalents 1.50002, $-1.50760$, and 0.000536718.  The analogue of \eqref{EoR constraint} is 
\eqref{EoRs weighted} with weights
\begin{equation}\label{weights 33v3}
w_1=\frac{78}{308},\quad w_2=\frac{78}{308}\,\frac{230}{231},\quad w_3=\frac{153}{308}\,\frac{230}{231}.
\end{equation}
The correlation between the effects of removal and the deuces-minus-aces counting system $(J_1,J_2,J_3)=(-1,1,0)$ is
$\rho\approx-0.999998$, and the regression coefficient is
$$
\gamma=-\frac{835{,}778{,}884}{555{,}776{,}529}\approx-1.50380.
$$

The analogues of \eqref{Z_n hat} and \eqref{Z_n^*} can be found by observing from \eqref{EoRs weighted} and \eqref{weights 33v3} that
$$
78\,E_1+78\,E_2+153\,E_3=-\frac{78}{230}E_1.
$$
Hence
\begin{align*}
\widehat Z_n&=\frac{1}{309-n}\sum_{j=n+1}^{309}(\mu-E_{X_j})\nonumber\\
&=\mu+\frac{1}{309-n}\bigg(\frac{78}{230}E_1+\sum_{j=1}^n E_{X_j}\bigg)
\end{align*}
and 
\begin{align*}
Z_n^*&=\mu+\frac{\gamma}{52}\bigg(\frac{52}{309-n}\bigg[\frac{78}{230}J_1+\sum_{j=1}^n J_{X_j}\bigg]\bigg)\nonumber\\
&=\mu+\frac{\gamma}{52}\bigg(\frac{52}{309-n}\bigg[-\frac{78}{230}+\sum_{j=1}^n J_{X_j}\bigg]\bigg)\nonumber\\
&=\mu+\frac{\gamma}{52}\bigg[-\frac{52}{309-n}\,\frac{78}{230}+\text{TC}_n^*\bigg],
\end{align*}
where
\begin{align*}
\mu&=E_{\std,(78,78,153)}((0,0,2),3)-E_{\spl,(78,78,153)}((0,0,2),3)\\
&=-\frac{229{,}736}{1{,}820{,}203}\approx-0.126214.
\end{align*}
This allows the card counter to know (approximately) when it is advantageous to depart from basic strategy when holding $(0,0,2)$ vs.~3.  Indeed, $Z_n^*>0$ is equivalent to 
$$
\text{TC}_n^*:=\frac{52(m_2-m_1)}{309-n}<-\frac{52\mu}{\gamma}+\frac{52}{309-n}\,\frac{78}{230}.
$$
The fraction 
$$
-\frac{52\mu}{\gamma}=-\frac{70{,}217{,}200{,}248}{16{,}088{,}743{,}517}\approx-4.36437
$$
plus the fraction $(52)(39)/[(309-n)115]$ is the variable index number for this departure.  As Griffin~\cite[p.~181]{G99} noted for blackjack, ``different change of strategy parameters will be required at different levels of the deck when the dealer's up card is an ace.''  If the adjusted true count $\text{TC}_n^*$ is less than this index number, standing on $(0,0,2)$ vs.~3 is called for instead of splitting.

Even with the extra $n$-dependent term, the index would be rounded to $-5$ (if $n\le260$), so this departure would seldom occur, and the bet size would be one unit when it did occur.  Therefore, the profit potential is rather limited, and an upper bound can be inferred from the third panel in Figure~\ref{ev-33v123}.

It should be pointed out that the methodology we have used to analyze strategy variation differs slightly from that of the blackjack literature and in particular from that of Griffin~\cite[pp.~72--90]{G99} and Schlesinger~\cite[Appendix D]{S18}.  The distinction was described by Griffin~\cite[Appendix to Chapter 6]{G99} as follows:

\begin{quote}
The strategy tables presented here are not the very best we could come up with in a particular situation.  As mentioned in this chapter more accuracy can be obtained with the normal approximation if we work with a 51 rather than a 52 card deck.  One could even have separate tables of effects for different two card player hands such as $(\text{T},6)$ v T.  Obviously a compromise must be reached, and my motivation has been in the direction of simplicity of exposition and ready applicability to multiple deck play.
\end{quote}

There are two issues here.  First, in analyzing a particular strategic situation, it is conventional to assume that the player has an ``abstract'' total, and to even regard the dealer's upcard as ``abstract.''  We regard this convention as  an unnecessary simplification, especially in snackjack.  It also explains why the player's two cards and dealer's upcard are not included in our ``adjusted'' true count;  indeed, those three cards are not part of the 309-card shoe on which the analysis is based.  An advantage of the conventional approach is that the sum of the effects of removal is 0, rather than some weighted average, with the exception of the cases in which the upcard is an ace or a ten in blackjack (an ace or a trey in snackjack).  Here Griffin~\cite[p.~197]{G99} achieved an EoR sum of 0 with an ace up by multiplying the EoR for ten by 36/35.  With a ten up, the EoR for ace is multiplied by 48/47.  This also seems to be the approach of Schlesinger~\cite[Appendix D]{S18}, but it sacrifices accuracy for simplicity.  The second issue is that the strategic EoRs (as well as the betting EoRs) are typically computed for the single-deck game and then converted to the multiple-deck games with the aid of a conversion factor (Griffin~\cite[Chapter 6]{G99}, Schlesinger~\cite[Appendix~D]{S18}).  Here, to maximize accuracy, we compute the effects of removal directly for the game we are interested in, 39-deck snackjack.

\begin{figure}[htb]
\centering
\includegraphics[width=3in]{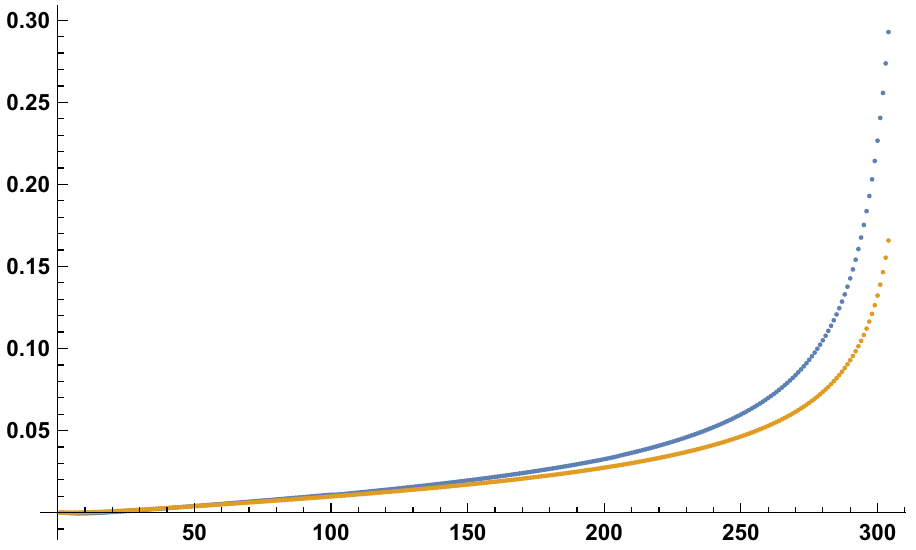}
\caption{\label{normal-approx}With $(0,0,2)$ vs.~1, the exact (blue) and approximate (orange) expected gain by departing from basic strategy when the rounded adjusted true count is at $-1$ or less, as a function of the number $n$ of cards seen (before the hand is dealt), $1\le n\le304$.}
\end{figure}

The computations that generated the first panel in Figure~\ref{ev-33v123} were exact, but similar computations cannot be done for six-deck blackjack (recall the 370 trillion distinguishable shoe compositions from Section~\ref{FTCC}).  Instead, approximate methods, developed by Griffin~\cite{G76}, are available, and it may be of some interest to see how accurate they are in the case of 39-deck snackjack.  With $Z_n^*$ defined by \eqref{Z_n^*-33v1} and $\sigma_J^2:=w_1+w_2-(w_2-w_1)^2$, a simple computation shows that $\E[Z_n^*]=\mu$ and 
\begin{equation}\label{SD(Z_n) in SJ}
\text{SD}(Z_n^*)=|\gamma|\sigma_J\sqrt{\frac{n}{(312-n)311}}.
\end{equation}
(The coefficient of the square root, $|\gamma|\sigma_J$, is sometimes written as $|\rho|\sigma_E$, but this is not quite the same thing.)  Notice that \eqref{SD(Z_n) in SJ} is proportional to \eqref{SD(Z_n)} with $N=312$.
By the normal approximation, $(Z_n^*-\mu)/\text{SD}(Z_n^*)$ is approximately $N(0,1)$.  Now in general, if $Z$ is $N(0,1)$, $\mu$ is real, and $\sigma>0$, then
$$
\E[(\mu+\sigma Z)^+]=\sigma\,\E[(Z+\mu/\sigma)^+]=\sigma\, \text{UNLLI}(-\mu/\sigma),
$$
where UNLLI stands for \textit{unit normal linear loss integral}~\cite[p.~87]{G99}, defined for real $x$ by
$$
\text{UNLLI}(x):=\E[(Z-x)^+]=\int_x^\infty (z-x)\phi(z)\,dz=\phi(x)-x(1-\Phi(x)),
$$
where $\phi$ and $\Phi$ denote the standard normal probability density function and cumulative distribution function.
We conclude from the normal approximation that, with $\sigma_n:=\text{SD}(Z_n^*)$,
\begin{equation}\label{UNLLI-approx}
\E[(Z_n^*)^+]\approx \E[(\mu+\sigma_n Z)^+]=\sigma_n\,\text{UNLLI}(-\mu/\sigma_n).
\end{equation}
Figure~\ref{normal-approx} shows that the quality of the approximation deteriorates over the course of the shoe.
If we average the approximate quantities over $1\le n\le 233$, we obtain $0.0139785$, which underestimates the exact value found above by about $13.8$\%.
  
In blackjack, the decision points with the greatest profit potential for varying basic strategy are the \textit{Illustrious 18} of Schlesinger~\cite[Chap.~5]{S18}.  Ten of the 18 involve \textit{stiffs}, player hands valued at hard 12--16, which when hit can be busted and when stood will lose unless the dealer busts.  
The only stiff in snackjack is a hard 5, so let us treat that case next.  We focus on the two-card hard 5, namely $2,3$. 

With $(0,1,1)$ vs.~1, the proportion of shoe compositions that call for a departure from basic strategy is 
$156{,}807/948{,}918\approx0.165248$.  With $(0,1,1)$ vs.~2, the proportion is $387{,}717/948{,}913\approx0.408591$.   With $(0,1,1)$ vs.~3, the proportion is $0/955{,}001=0$.  In particular, it is never correct to depart from basic strategy with $(0,1,1)$ vs.~3.  Thus, we discuss only the other two situations.  

\begin{figure}[htb]
\centering
\includegraphics[width=2.25in]{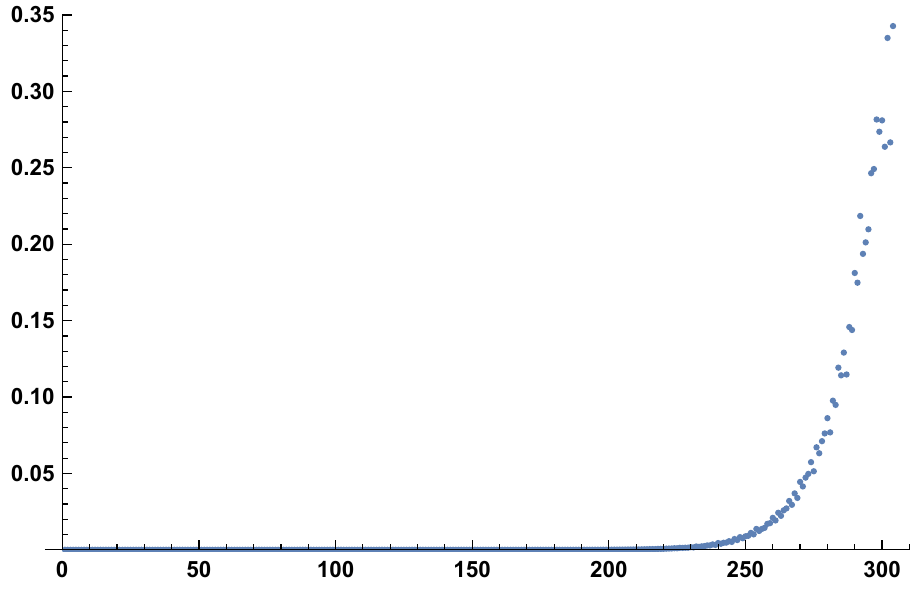}\quad
\includegraphics[width=2.25in]{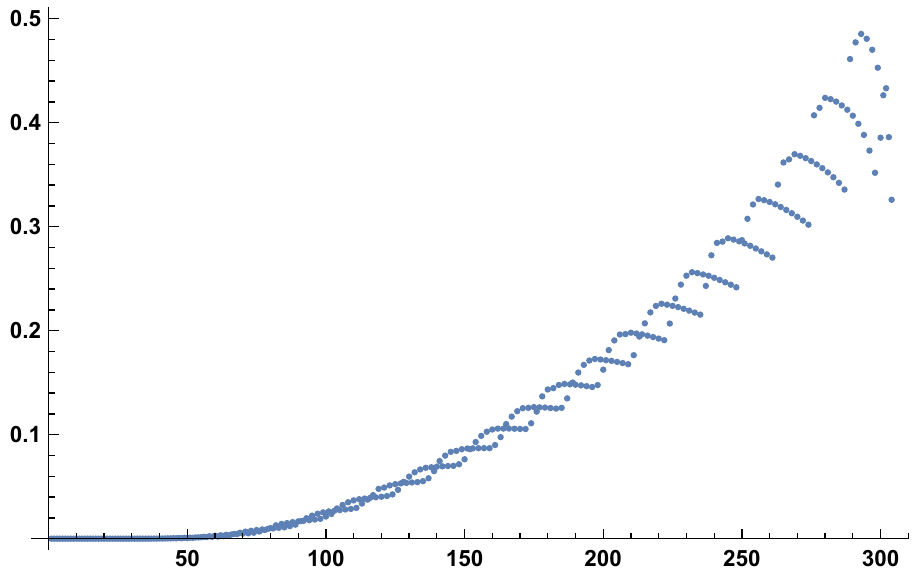}
\caption{\label{prob-23v1or2}With $(0,1,1)$ vs.~1 (left) and $(0,1,1)$ vs.~2 (right), the probability that a departure from basic strategy is called for when $n$ cards have been seen (before the hand is dealt), as a function of $n$, $1\le n\le304$.}
\end{figure}

With $(0,1,1)$ vs.~1, the probability that a departure from basic strategy is called for when $n$ cards have been seen (before the hand is dealt) can be evaluated as in \eqref{weightedave1}.  Similar computations can be done for $(0,1,1)$ vs.~2.  The two expressions are graphed in Figure~\ref{prob-23v1or2}.

\begin{figure}[htb]
\centering
\includegraphics[width=2.25in]{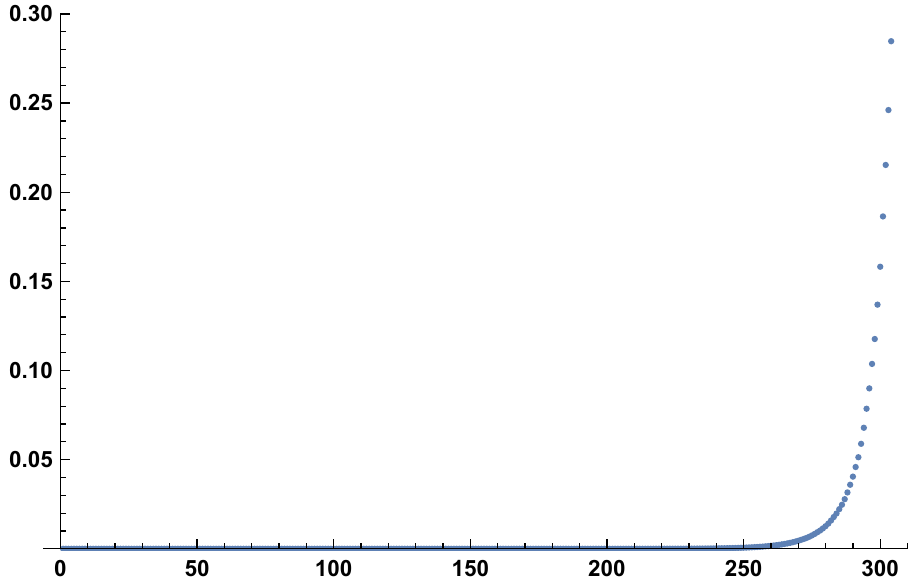}\quad
\includegraphics[width=2.25in]{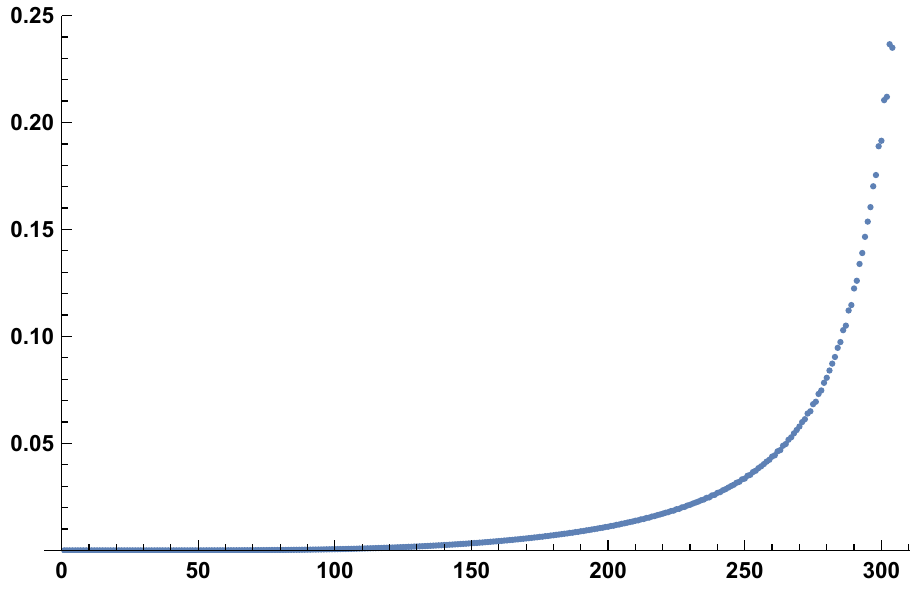}
\caption{\label{ev-23v1or2}With $(0,1,1)$ vs.~1 (left) and $(0,1,1)$ vs.~2 (right), the additional expectation that a departure from basic strategy provides when $n$ cards have been seen (before the hand is dealt), as a function of $n$, $1\le n\le304$.}
\end{figure}

The additional expectation that such a departure provides is computed as in \eqref{weightedave1} and \eqref{alpha_1-alt}.  The two expressions are graphed in  Figure~\ref{ev-23v1or2}.

We assume that the card counter employs the deuces-minus-aces count, $(J_1,J_2,J_3)=(-1,1,0)$.  This count has correlation 0.0325 (resp., 0.503006) with the effects of removal for $2,3$ vs.~2 (resp., $2,3$ vs.~1), implying that card counting is ineffective in this setting.  Indeed, potentially profitable strategy variation with $(0,1,1)$ vs.~2 will be recognized only by the player who keeps track of all three denominations, permitting the evaluation of the count $(K_1,K_2,K_3)=(1,1,-1)$.

We conclude this section by considering the hand $(1,1,0)$, for which basic strategy is to hit vs.~1 and double vs.~2 and 3.  The interesting thing about this hand is that all three strategies, standing, hitting, and doubling, are optimal for certain shoe compositions.  Thus, a single difference of two expectations for each upcard is insufficient to analyze departures from basic strategy.  In Appendix~D, we provide formulas for $E_{\std,(n_1,n_2,n_3)}((1,1,0),u)$, $E_{\hit,(n_1,n_2,n_3)}((1,1,0),u)$, and $E_{\dbl,(n_1,n_2,n_3)}((1,1,0),u)$ for $u=1,2,3$.

With $(1,1,0)$ vs.~1, the number of shoe compositions is $N_1:=942{,}755$, the proportion for which it is optimal to stand is $110{,}497/N_1\approx0.117206$, to hit (basic strategy) is $477{,}244/N_1\approx0.506223$, and to double is $355{,}014/N_1\approx0.376571$.
With $(1,1,0)$ vs.~2, the number of shoe compositions is $N_2:=942{,}907$, the proportion for which it is optimal to stand is $60{,}215/N_2\approx0.063861$, to hit is $335{,}256/N_2\approx0.355556$, and to double (basic strategy) is $547{,}436/N_2\approx0.580583$.
With $(1,1,0)$ vs.~3, the number of shoe compositions is $N_3:=948{,}996$, the proportion for which it is optimal to stand is $152{,}902/N_3\approx0.161120$, to hit is $363{,}196/N_3\approx0.382716$, and to double  (basic strategy) is $432{,}898/N_3\approx0.456164$.

\begin{figure}[htb]
\centering
\includegraphics[width=2.25in]{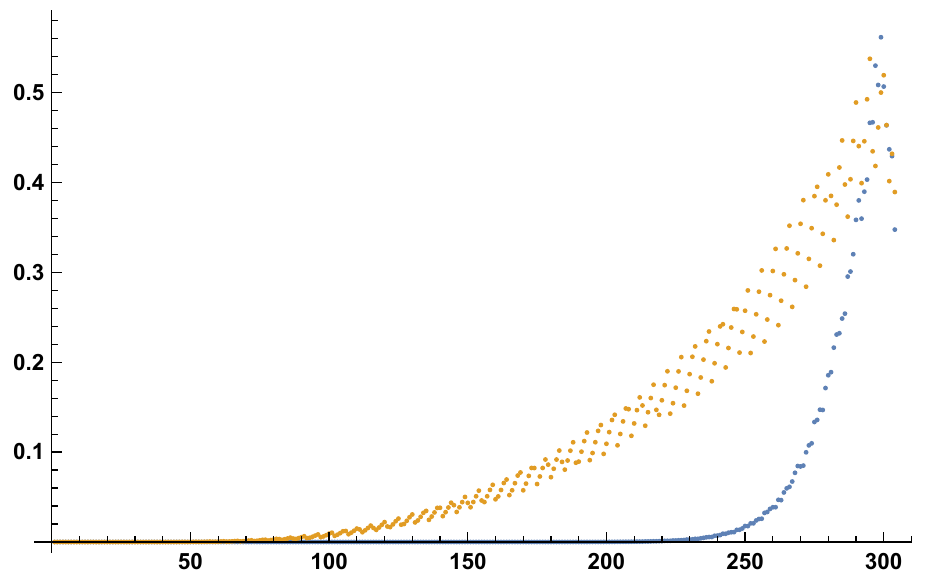}\quad
\includegraphics[width=2.25in]{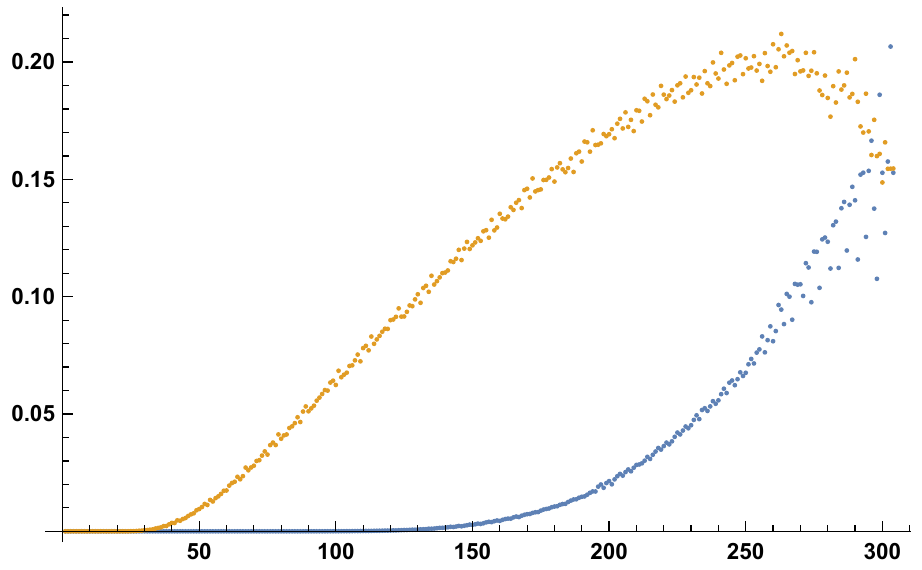}\vglue3mm
\includegraphics[width=2.25in]{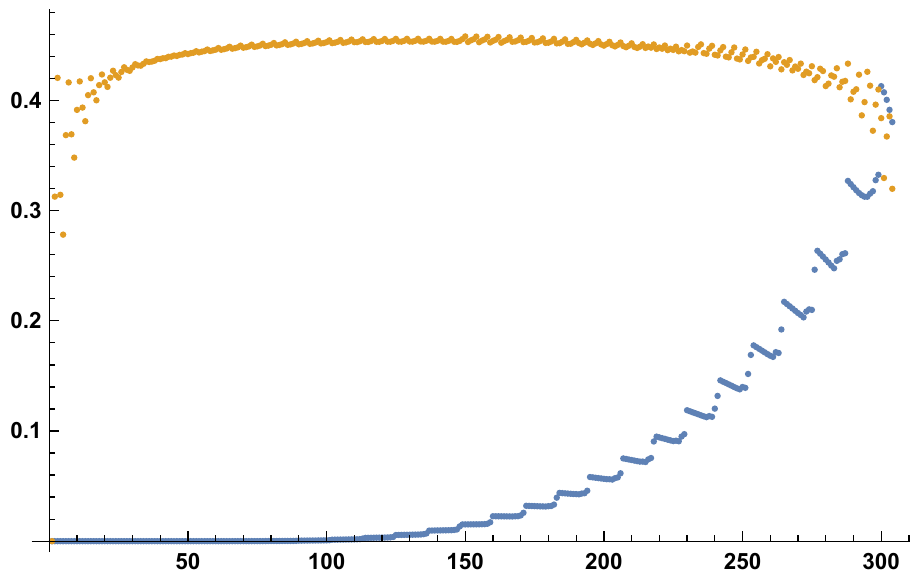}
\caption{\label{prob-12v123}With $(1,1,0)$ vs.~1 (top left), $(1,1,0)$ vs.~2 (top right), and $(1,1,0)$ vs.~3 (bottom), the probability that a departure from basic strategy is called for when $n$ cards have been seen (before the hand is dealt), as a function of $n$, $1\le n\le304$.  The blue graph is for standing, and the orange graph is for doubling vs.~1 and hitting vs.~2 or 3.}
\end{figure}

With $(1,1,0)$ vs.~1, the probability that a departure from basic strategy is called for when $n$ cards have been seen (before the hand is dealt) can be evaluated as in \eqref{weightedave1}.  We treat the cases of standing and doubling separately.  Similar computations can be done for $(1,1,0)$ vs.~2 and for $(1,1,0)$ vs.~3, with the cases of standing and hitting treated separately.  The six expressions are graphed in Figure~\ref{prob-12v123}.

\begin{figure}[htb]
\centering
\includegraphics[width=2.25in]{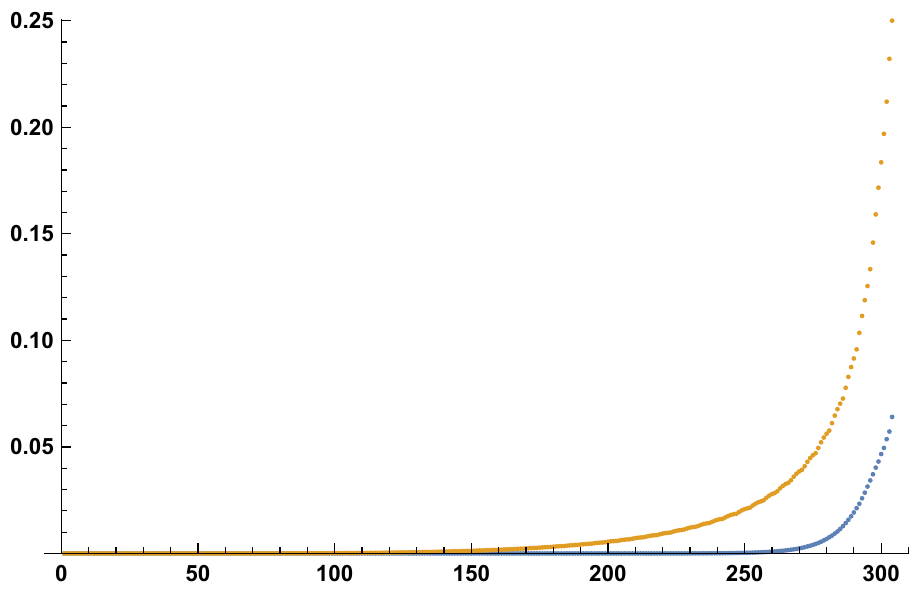}\quad
\includegraphics[width=2.25in]{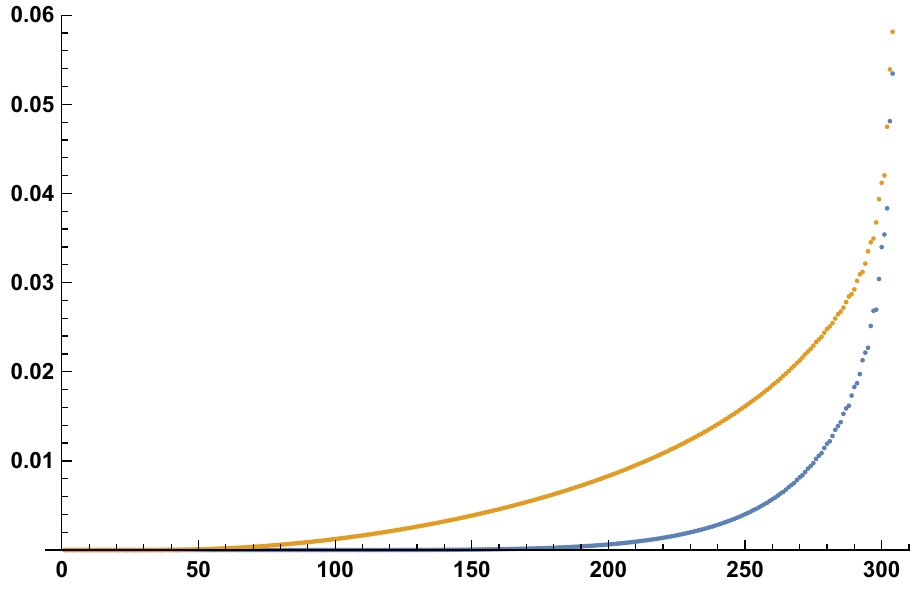}\vglue3mm
\includegraphics[width=2.25in]{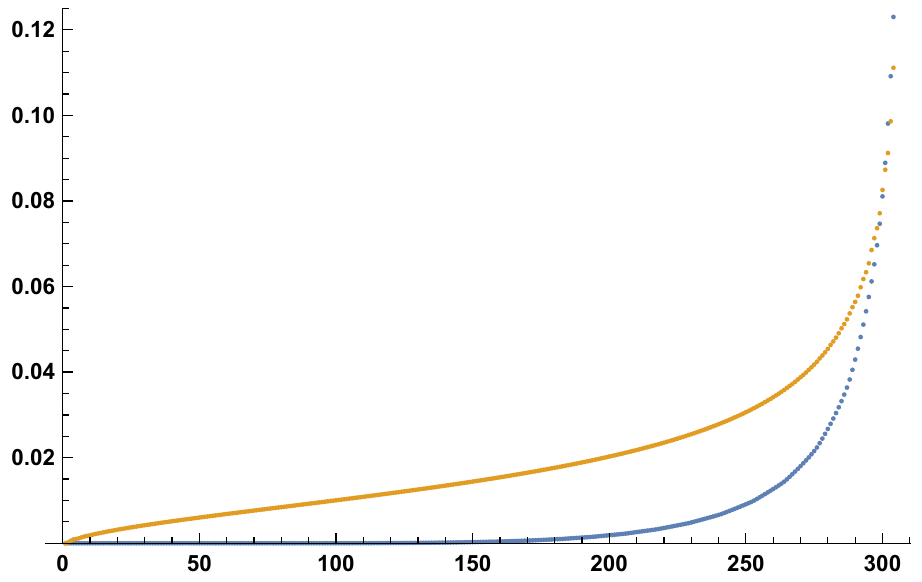}
\caption{\label{ev-12v123}With $(1,1,0)$ vs.~1 (left top), $(1,1,0)$ vs.~2 (right top), and $(1,1,0)$ vs.~3 (bottom), the additional expectation that a departure from basic strategy provides when $n$ cards have been seen (before the hand is dealt), as a function of $n$, $1\le n\le304$. The blue graph is for standing, and the orange graph is for doubling vs.~1 and hitting vs.~2 or 3.  Notice that the vertical scales differ considerably.}
\end{figure}

The additional expectation that such a departure provides is computed as in \eqref{weightedave1} and \eqref{alpha_1-alt}, with the cases of standing and doubling (or standing and hitting) treated separately.  The six expressions are graphed in Figure~\ref{ev-12v123}.

Later, in Table~\ref{strategy variation analysis}, we summarize the strategy variation analysis done for $(0,0,2)$, $(0,1,1)$, and $(1,1,0)$.  Since there are six cases for $(1,1,0)$, we do not work through the details for all of them, only for $(1,1,0)$ vs.~1 when doubling is the alternative to hitting, and for $(1,1,0)$ vs.~3 when hitting is the alternative to doubling.  In the first case, the effects of removal are
\begin{align}\label{EoR-12v1d}
\text{EoR}(i)&:=E_{\dbl,(76,77,156)-\bm e_i}((1,1,0),1)-E_{\hit,(76,77,156)-\bm e_i}((1,1,0),1)\nonumber\\
&\qquad{}-[E_{\dbl,(76,77,156)}((1,1,0),1)-E_{\hit,(76,77,156)}((1,1,0),1)]
\end{align}
for $i=1,2,3$.
The numbers \eqref{EoR-12v1d}, multiplied by 308, are
$$
E_1=-\frac{110{,}244{,}809{,}177}{140{,}671{,}381{,}176},\quad E_2=-\frac{104{,}264{,}078}{48{,}708{,}927},\quad E_3=-\frac{293{,}864{,}077}{438{,}380{,}343},
$$
with decimal equivalents $-0.783705$, 2.14055, and $-0.67034$.  The analogue of \eqref{EoR constraint} is 
\eqref{EoRs weighted} with weights
\begin{equation}\label{weights-12v1}
w_1=\frac{76}{308}\,\frac{152}{153},\quad w_2=\frac{77}{308}\,\frac{152}{153},\quad w_3=\frac{156}{308}.
\end{equation}
The correlation between the effects of removal and the deuces-minus-aces counting system $(J_1,J_2,J_3)=(-1,1,0)$ is
$\rho\approx0.836717$, and the regression coefficient is
$$
\gamma=\frac{415{,}321{,}501{,}405}{283{,}193{,}701{,}578}\approx1.46656.
$$

The analogues of \eqref{Z_n hat} and \eqref{Z_n^*} can be found by observing from \eqref{EoRs weighted} and \eqref{weights-12v1} that
$$
76\,E_1+77\,E_2+156\,E_3=-\frac{156}{152}\,E_3.
$$
Hence
\begin{align*}
\widehat Z_n&=\frac{1}{309-n}\sum_{j=n+1}^{309}(\mu-E_{X_j})\nonumber\\
&=\mu+\frac{1}{309-n}\bigg(\frac{156}{152}\,E_3+\sum_{j=1}^n E_{X_j}\bigg)
\end{align*}
and 
\begin{align*}
Z_n^*&=\mu+\frac{\gamma}{52}\bigg(\frac{52}{309-n}\bigg[\frac{156}{152}\,J_3+\sum_{j=1}^n J_{X_j}\bigg]\bigg)\nonumber\\
&=\mu+\frac{\gamma}{52}\bigg(\frac{52}{309-n}\sum_{j=1}^n J_{X_j}\bigg)=\mu+\frac{\gamma}{52}\text{TC}_n^*,
\end{align*}
where
\begin{align*}
\mu&=E_{\dbl,(76,77,156)}((1,1,0),1)-E_{\hit,(76,77,156)}((1,1,0),1)\\
&=-\frac{452{,}457{,}716}{3{,}750{,}587{,}379}\approx-0.120636.
\end{align*}
This allows the card counter to know (approximately) when it is advantageous to depart from basic strategy when holding $(1,1,0)$ vs.~1.  Indeed, $Z_n^*>0$ is equivalent to 
$$
\text{TC}_n^*:=\frac{52(m_2-m_1)}{309-n}>-\frac{52\mu}{\gamma}.
$$
The fraction 
$$
-\frac{52\mu}{\gamma}=\frac{136{,}790{,}636{,}362{,}848}{31{,}979{,}755{,}608{,}185}\approx4.27741
$$
is the index number for this departure.  If the adjusted true count is greater than this index number, doubling on $(1,1,0)$ vs.~1 is called for instead of hitting.

Here we would round the index to $5$, and this play would occur relatively seldom (the rounded adjusted true count would have to be at 5 or more).  On the other hand, we would be betting at least five units, so the profit potential is nontrivial.  Indeed, we can compute it precisely using \eqref{weightedave1} with
\begin{align*}
\alpha_1(m_1,m_2,m_3)&:=[E_{\dbl,(76-m_1,77-m_2,156-m_3)}((1,1,0),1)\\
&\qquad{}-E_{\hit,(76-m_1,77-m_2,156-m_3)}((1,1,0),1)]\\
&\qquad\qquad{}\cdot\bm1\bigg\{\bigg[\frac{52(m_2-m_1)}{309-n}\bigg]\ge5\bigg\}\\
&\qquad\qquad{}\cdot \min\bigg(\bigg[\frac{52(m_2-m_1)}{309-n}\bigg],6\bigg).
\end{align*}
The average of these expectations over $1\le n\le 233$ is about 0.00522862.

In the case $(1,1,0)$ vs.~3 when hitting is the alternative to doubling, the effects of removal are
\begin{align}\label{EoR-12v3h}
\text{EoR}(i)&:=E_{\hit,(77,77,155)-\bm e_i}((1,1,0),3)-E_{\dbl,(77,77,155)-\bm e_i}((1,1,0),3)\nonumber\\
&\qquad{}-[E_{\hit,(77,77,155)}((1,1,0),3)-E_{\dbl,(77,77,155)}((1,1,0),3)]
\end{align}
for $i=1,2,3$.
The numbers \eqref{EoR-12v3h}, multiplied by 308, are
$$
E_1=\frac{213{,}948{,}581}{332{,}366{,}796},\quad E_2=-\frac{1{,}261{,}124{,}737}{997{,}100{,}388},\quad E_3=\frac{170{,}145{,}283}{553{,}944{,}660},
$$
with decimal equivalents $0.643712$, $-1.26479$, and $0.307152$.  The analogue of \eqref{EoR constraint} is 
\eqref{EoRs weighted} with weights
\begin{equation}\label{weights-12v3}
w_1=\frac{77}{308},\quad w_2=\frac{77}{308}\,\frac{231}{232},\quad w_3=\frac{155}{308}\,\frac{231}{232}.
\end{equation}
The correlation between the effects of removal and the deuces-minus-aces counting system $(J_1,J_2,J_3)=(-1,1,0)$ is
$\rho\approx-0.908998$, and the regression coefficient is
$$
\gamma=-\frac{146{,}742{,}675{,}541}{153{,}885{,}826{,}548}\approx-0.953581.
$$

The analogues of \eqref{Z_n hat} and \eqref{Z_n^*} can be found by observing from \eqref{EoRs weighted} and \eqref{weights-12v3} that
$$
77\,E_1+77\,E_2+155\,E_3=-\frac{1}{3}\,E_1.
$$
Hence
\begin{align*}
\widehat Z_n&=\frac{1}{309-n}\sum_{j=n+1}^{309}(\mu-E_{X_j})\nonumber\\
&=\mu+\frac{1}{309-n}\bigg(\frac{1}{3}\,E_1+\sum_{j=1}^n E_{X_j}\bigg)
\end{align*}
and 
\begin{align*}
Z_n^*&=\mu+\frac{\gamma}{52}\bigg(\frac{52}{309-n}\bigg[\frac{1}{3}\,J_1+\sum_{j=1}^n J_{X_j}\bigg]\bigg)\nonumber\\
&=\mu+\frac{\gamma}{52}\bigg(\frac{52}{309-n}\bigg[-\frac13+\sum_{j=1}^n J_{X_j}\bigg]\bigg)\nonumber\\
&=\mu+\frac{\gamma}{52}\bigg(\!\!-\frac{52}{3(309-n)}+\text{TC}_n^*\bigg),
\end{align*}
where
\begin{align*}
\mu&=E_{\hit,(77,77,155)}((1,1,0),3)-E_{\dbl,(77,77,155)}((1,1,0),3)\\
&=-\frac{191}{70{,}992}\approx-0.00269044.
\end{align*}
This allows the card counter to know (approximately) when it is advantageous to depart from basic strategy when holding $(1,1,0)$ vs.~3.  Indeed, $Z_n^*>0$ is equivalent to 
$$
\text{TC}_n^*:=\frac{52(m_2-m_1)}{309-n}>-\frac{52\mu}{\gamma}+\frac{52}{3(309-n)}.
$$
The fraction 
$$
-\frac{52\mu}{\gamma}=-\frac{21{,}529{,}102{,}283}{146{,}742{,}675{,}541}\approx-0.146713
$$
plus the fraction $52/[3(309-n)]$
is the index number for this departure.  If the adjusted true count is less than this index number, hitting on $(1,1,0)$ vs.~3 is called for instead of doubling.

Here we would round the index to $-1$, ignoring the extra $n$-dependent term with perhaps a negligible effect, and this play would occur relatively often (the rounded adjusted true count would have to be at $-1$ or less).  Here we would be betting only one unit unit, but the profit potential is nontrivial.  Indeed, we can compute it precisely using \eqref{weightedave1} with
\begin{align*}
\alpha_1(m_1,m_2,m_3)&:=[E_{\hit,(77-m_1,77-m_2,155-m_3)}((1,1,0),3)\\
&\qquad{}-E_{\dbl,(77-m_1,77-m_2,155-m_3)}((1,1,0),3)]\\
&\qquad\qquad{}\cdot\bm1\bigg\{\bigg[\frac{52(m_2-m_1)}{309-n}\bigg]\le-1\bigg\}.
\end{align*}
The average of these expectations over $1\le n\le 233$ is approximately 0.0105474.

\begin{table}[htb]
\caption{\label{strategy variation analysis}Analysis of several departures from basic strategy.  $\rho$ is the correlation between the effects of removal for the strategic situation and the deuces-minus-aces count $(-1,1,0)$.  The inequality listed under ``departure criterion'' is what $[\text{TC}_n^*]$ must satisfy for a departure from basic strategy to be called for.  The average EV is the player's expectation averaged over $1\le n\le 233$ when betting $\max(1,\min(\text{[TC}_n^*],6))$.  An asterisk in the ``index'' column signifies a missing $n$-dependent term.  We omit rows for $(0,2,0)$ and $(2,0,0)$ vs.~1, 2, and 3, each of which would have a 0 in the last column.}
\catcode`@=\active \def@{\hphantom{0}}
\catcode`#=\active \def#{\hphantom{$-$}}
\tabcolsep=.12cm
\begin{center}
\begin{tabular}{cccccccr}
\noalign{\smallskip}
\hline
\noalign{\smallskip}
nos of     & up & bs & alt & corr.  & index & departure   & $10^6\times$  \\
1s, 2s, 3s &    &    &     & $\rho$ &       & criterion &  ave EV   \\
\noalign{\smallskip}
\hline
\noalign{\smallskip}
$(0,0,2)$  & 1 & Spl & S  & $-0.986$ & $-0.718$ & $\le-1$ & $16{,}214$ \\
$(0,0,2)$  & 2 & Spl & S  & $-0.944$ & $-5.44\phantom{^*}$ & $\le-6$ & 250 \\
$(0,0,2)$  & 3 & Spl & S  & $-1.000$ & $-4.36^*$ & $\le-5$  & 688 \\
\noalign{\medskip}
$(0,1,1)$  & 1 & H & S  & #$0.503$ & #25.4$\phantom{^*}$ & -- & 0 \\
$(0,1,1)$  & 2 & H & S  & #$0.033$ & #91.7$\phantom{^*}$ & -- & 0 \\
$(0,1,1)$  & 3 & H & S  & $-0.149$ & $-184.^*$ & -- & 0 \\
\noalign{\medskip}
$(1,1,0)$  & 1 & H & S & $-0.472$ & $-22.5\phantom{^*}$ & -- & 0 \\
$(1,1,0)$  & 1 & H & D & #$0.837$ & #$4.28\phantom{^*}$ & $\ge+5$   & $5{,}229$ \\
$(1,1,0)$  & 2 & D & S & $-0.940$ & $-5.63\phantom{^*}$ & $\le-6$ & 88 \\
$(1,1,0)$  & 2 & D & H & $-0.661$ & $-3.06\phantom{^*}$ & $\le-4$ & 650 \\
$(1,1,0)$  & 3 & D & S & $-1.000$ & $-4.36^*$ & $\le-5$ & 687 \\
$(1,1,0)$  & 3 & D & H & $-0.909$ & $\phantom{^*}$$-0.147^*$ & $\le-1$ & $10{,}547$ \\
\noalign{\smallskip}
\hline
\end{tabular}
\end{center}
\end{table}

The results of Table~\ref{strategy variation analysis} show that only three departures from basic strategy are ``illustrious.''  $(0,0,2)$ vs.~1 (stand instead of split) offers the greatest profit potential, then $(1,1,0)$ vs.~3 (hit instead of double), and finally $(1,1,0)$ vs.~1 (double instead of hit).  None of the others is close to these three.

This table also emphasizes another distinction between our approach and the blackjack literature.  Our strategy variations are formulated as departures from basic strategy.  Instead of saying ``stand instead of split if the rounded adjusted true count is $-1$ or less,'' it would be more conventional to say, ``split instead of stand if the rounded adjusted true count is 0 or more.''  This way, all inequalities point in the same direction $(\ge)$.

\section{What does snackjack tell us about blackjack?}\label{conclusions}

The simpler a toy model is, the fewer features it shares with the original.  Grayjack is closer to blackjack than is snackjack.  For example, the proportions of aces and tens in blackjack are 1/13 and 4/13.  These proportions are maintained in grayjack for aces and sixes, thereby making naturals about as frequent and allowing insurance.  In snackjack, the proportions of aces and treys are unavoidably rather different.  In blackjack the numbers of pat totals (17--21) and stiff totals (12--16) are the same, five each.  In grayjack the numbers (8--10 stiff and 11--13 pat) are also the same, three each.  But in snackjack the numbers (5 stiff and 6--7 pat) are different, again unavoidably.  This makes it more difficult to bust, mitigating the dealer's principal advantage, the double bust.  Ultimately, we felt that the benefits of having a hand-computable toy model of blackjack outweighed the drawbacks of a significant player advantage and a largely upcard-independent basic strategy.  Actually, more important than hand-computability are the explicit formulas available for basic strategy expectations with arbitrary shoe compositions.  This allows exact computation of quantities that can only be estimated at blackjack.

What then have we learned about blackjack from its computable toy model, snackjack?

\begin{itemize}
\item The derivation of basic strategy at blackjack is conceptually very simple, despite its computational complexity.
Basic strategy for blackjack is now so well known and understood that there is little insight to be gained by deriving basic strategy for snackjack or grayjack.  Nevertheless, perhaps surprising to some is the conceptual simplicity of the basic strategy derivation, as illustrated by the tree diagram in Figure~\ref{tree}.  The corresponding tree diagram for standing with a pair of tens vs.\ a playable ace in six-deck blackjack would have 8,496 terminal vertices~\cite[p.~158]{G99} instead of four, but \textit{conceptually} it is the same thing.

\item It is truly remarkable, as has been noted elsewhere~\cite[p.~17]{G99}, that single-deck blackjack (under classic Las Vegas Strip rules), which was played long before it was analyzed, turned out to be an essentially fair game, with a player advantage of about four hundredths of 1\%.  The present study emphasizes the sensitivity of basic strategy expectations to minor rules changes.  For example, the rule ``A player natural pays even money, with the exception that it loses to a dealer natural,'' reduces the player advantage at double-deck snackjack from $+16.3$\% to $-0.0959$\%.  The less extreme rules change in blackjack in which untied player naturals pay 6 to 5 instead of 3 to 2 has a smaller but still significant effect, as every advantage player will acknowledge.

\item A formula for basic strategy expectation in six-deck blackjack as a function of shoe composition, if found, would likely be highly impractical.  A polynomial in three variables of degree 8 or less has at most $\binom{8+3}{3}=165$ terms, so the 147 terms of the polynomial in \eqref{E(n_1,n_2,n_3)} is not surprising.  The analogous polynomial in blackjack would have at most $\binom{m+10}{10}$ terms, where $m$ is the maximum number of cards needed to complete a round.  In six-deck blackjack that number is at least 24 (e.g., two hands of $1,1,1,1,1,1,6,1,1,1,1,5$), even before considering splits, and $\binom{24+10}{10}=131{,}128{,}140$.  The actual number of terms in the blackjack basic strategy expectation polynomial would likely be somewhat smaller but still highly impractical.

\item With $Z_n$ being the player's conditional expectation when $n$ cards have been seen, $\E[(Z_n)^+]$ (or $\E[(Z_n-\nu)^+]$) can be computed directly for the game of red-and-black (see \eqref{E[Z_n^+]}) and for 39-deck snackjack (see \eqref{FTCC-meanpos}), assuming basic strategy.  As we have explained, such computations for six-deck blackjack are likely impossible, but perhaps computer simulation would give the best results.  Another potential approach would be to approximate $Z_n$ by its linearization $\widehat Z_n$ based on EoRs, and then use a normal approximation involving the UNLLI function, much as we did in \eqref{UNLLI-approx}.  There are some things in blackjack that simply cannot be known exactly.

\item Section~\ref{bet variation} contains several computations for 39-deck snackjack that cannot be replicated for six-deck blackjack.  If they could be, we would presumably reach the same conclusions as we do at snackjack.  Specifically, we computed the $L^1$ distances between $Z_n$ and its linearization $\widehat Z_n$ based on EoRs and its linearization $Z_n^*$ based on the chosen card-counting system.  We find from Table~\ref{L1-distance} that, for the first 2/3 of the shoe, the bulk of the error in approximating $Z_n$ by $Z_n^*$ is explained by the use of the level-one deuces-minus-aces point count in place of the EoRs; the nonlinearity effect is relatively inconsequential.  Another finding, based on limited evidence, was that the betting efficiency of a card-counting system is well approximated by the betting correlation, that is, the correlation between the EoRs and the numbers of the point count.  The latter is computable for blackjack, whereas the former is not (except by simulation).  There is of course a theoretical reason for this~\cite[p.~51]{G99}.

\item A surprise to us was the extent to which the distribution of the true count at snackjack departs from normality.  This is undoubtedly true at blackjack as well, but less easy to verify.  For snackjack, it is a consequence of Figure~\ref{count-dist}, which shows that the rounded true count fails to be discrete normal for some choices of $n$.  The $n=260$ case is what we expected, whereas the $n=234$ case illustrates what can happen.  A more complete analysis than that done for the figure shows that the distribution of the rounded true count is bimodal if and only if $105\le n\le 138$ or $209\le n\le255$.  Theory tells us that the true count is asymptotically normal but of course this lacks rigor because we never let $N$ (the number of cards in the shoe) tend to infinity; instead, $N$ is fixed at 312.  As Griffin~\cite[p.~38]{G99} put it, ``the proof of the pudding is in the eating.''

\item As we mentioned in Section~\ref{strategy variation}, it is conventional in blackjack to compute the effects of removal based on a 52-card deck, then multiply them by a conversion factor for multiple-deck applications.  The justification for this is based on the following observation.  Let $\text{EoR}_N(i)$ denote the effect of removal of card value $i$ from a deck of $N$ cards on basic strategy expectation.  Then it can be shown that $\lim_{N\to\infty}(N-1)\text{EoR}_N(i)$ exists, and therefore $\text{EoR}_{312}(i)$ is approximately equal to $(51/311)\text{EoR}_{52}(i)$, for example.  We can use snackjack to investigate how accurate we can expect this approximation to be for six-deck blackjack.  For snackjack, Table~\ref{SJ-eors} displays the relevant data.  The correlation between the $N=52$ EoRs and the $N=312$ EoRs is 0.999135.  The result is that the approximate effects of removal based on a 52-card ($6\frac12$-deck) pack, instead of a 312-card (39-deck) shoe, are considerably less accurate than our level-six point count but substantially more accurate than our level-one point count.  This can also confirmed in terms of $L^1$ distances, as in Table~\ref{L1-distance}.

\begin{table}[htb]
\caption{\label{SJ-eors}Effects of removal on snackjack's basic strategy expectation.  For simplicity we used \eqref{E(n_1,n_2,n_3)} to evaluate these numbers, even though basic strategy in the case $N=52$ differs slightly from the strategy implicit in \eqref{E(n_1,n_2,n_3)}.  The entries for $N=312$ coincide with those of Table~\ref{EoRs}.}
\catcode`#=\active \def#{\hphantom{$-$}}
\catcode`@=\active \def@{\hphantom{0}}
\tabcolsep=1.3mm
\renewcommand{\arraystretch}{1.}
\begin{center}
\begin{tabular}{cccc}
\hline
\noalign{\smallskip}
$N$ & $(N-1)\text{EoR}_N(1)$ & $(N-1)\text{EoR}_N(2)$ & $(N-1)\text{EoR}_N(3)$ \\
\noalign{\smallskip}\hline
\noalign{\smallskip}
52  & $-0.516148$ & 0.711619 & $-0.0977352$ \\
104 & $-0.490108$ & 0.702236 & $-0.106064$ \\
312 & $-0.473605$ & 0.696413 & $-0.111404$ \\
$\infty$ & $-0.465576$ & 0.693604 & $-0.114014$ \\
\noalign{\smallskip}\hline
\end{tabular}
\end{center}
\end{table}

\item In snackjack we have seen that some strategy variation decisions (such as standing instead of hitting a hard 5) are not well suited to the deuces-minus-aces count.  Similarly, and it is well known to experts, some strategy variation decisions in blackjack are not well suited to the Hi-Lo count (14 vs.~10 and 16 vs.~7 are two examples mentioned by Schlesinger~\cite[p.~57]{S18}).
\end{itemize}

In conclusion, this paper provides theoretical support for a conclusion for which abundant anecdotal evidence exists, namely that card counting works, despite its attempt to linearize a function that is clearly nonlinear.

\section*{Acknowledgments}

We thank Don Schlesinger for valuable comments on a preliminary draft of the manuscript.  The work of SNE was partially supported by a grant from the Simons Foundation (429675).  The work of JL was supported by a 2019 Yeungnam University Research Grant.

\section*{Appendix A}
Basic strategy for grayjack is complicated by many composition-dependent exceptions in the single-deck game (Table~\ref{GJ-BS-1-deck}), but has a rather simple form in the 24-deck game (Table~\ref{GJ-BS-24-deck}).

\begin{table}[htb]
\caption{\label{GJ-BS-1-deck}Composition-dependent basic strategy for single-deck grayjack.}
\catcode`@=\active \def@{\hphantom{0}}
\catcode`#=\active \def#{\hphantom{$\;^1$}}
\tabcolsep=.18cm
\begin{small}
\begin{center}
\begin{tabular}{cccccccc}
\hline
\noalign{\smallskip}
\multicolumn{2}{c}{player} & \multicolumn{6}{c}{dealer upcard}\\
\noalign{\smallskip}
cards&total& $1$  & $2$ & $3$ & $4$ & $5$ & $6$ \\
\noalign{\smallskip}\hline
\noalign{\smallskip}
&h11--h13&S&S&S&S&S&S\\
$4,6$&h10&H&S&S&H&H&H\\
$3,3,4$&h10&S&S&--&S&H&S\\
$2,4,4$&h10&S&S&S&--&H&H\\
$2,3,5$&h10&S&S&S&S&H&S\\
$2,2,6$&h10&H&--&S&S&S&S\\
$2,2,3,3$&h10&S/H&--&--&S&S&S\\
$1,4,5$&h10&--&S&S&S&H&H\\
$1,3,6$&h10&--&S&S&S&S&S\\
$1,2,3,4$&h10&--&S&S&S&S&S\\
$1,2,2,5$&h10&--&--&S&S&S&S\\
$4,5$&h9&H&S&S&S&H&H\\
$3,6$&h9&H&H&S&S&H&H\\
$2,3,4$&h9&S&S&S&S&H&H\\
$2,2,5$&h9&H&--&S&S&H&H\\
$1,4,4$&h9&--&S&S&--&H&H\\
$1,3,5$&h9&--&S&S&S&H&H\\
$1,2,6$&h9&--&H&S&S&H&H\\
$1,2,3,3$&h9&--&S&--&S&H&H\\
$1,2,2,4$&h9&--&--&S&S&H&H\\
$3,5$&h8&H&S&S&S&H&H\\
$2,6$&h8&H&H&H&H&H&H\\
$2,3,3$&h8&H&H&--&S&H&H\\
$2,2,4$&h8&H&--&S&S&H&H\\
$1,3,4$&h8&--&S&S&S&H&H\\
$1,2,5$&h8&--&H&S&S&H&H\\
$1,2,2,3$&h8&--&--&S&S&H&H\\
$3,4$&h7&D&D&D&D&D&D\\
$2,5$&h7&D&D&D&D&D&D\\
$2,2,3$&h7&H&--&H&H&H&H\\
$2,4$&h6&H&D&D&D&D&H\\
$2,3$&h5&H&H&D&D&H&H\\
\noalign{\smallskip}\hline
\noalign{\smallskip}
&s12--s13&--&S&S&S&S&S\\
$1,4$&s11&--&S&D&D&H&H\\
$1,2,2$&s11&--&--&S&S&S&S\\
$1,3$&s10&--&H&H&D&H&H\\
$1,2$&s9&--&H&D&D&D&H\\
\noalign{\smallskip}\hline
\noalign{\smallskip}
$6,6$&&S&S&S&S&S&S\\
$5,5$&&H&Spl&Spl&Spl&--&H\\
$4,4$&&Spl&Spl&Spl&--&Spl&Spl\\
$3,3$&&H&D&--&D&D&H\\
$2,2$&&Spl&--&Spl&Spl&Spl&Spl\\
\noalign{\smallskip}
\hline
\noalign{\smallskip}
\end{tabular}
\end{center}
\end{small}
\end{table}
\afterpage{\clearpage}

\begin{table}[htb]
\caption{\label{GJ-BS-24-deck}Composition-dependent basic strategy for 24-deck grayjack.  More generally, this table applies to $d$-deck grayjack if $21\le d\le 59$; for $d\ge60$ there is one change: The entry for $(2,2)$ vs.~2 is H instead of Spl.  (DH = double if allowed, hit if not.)}
\catcode`@=\active \def@{\hphantom{0}}
\catcode`#=\active \def#{\hphantom{$\;^1$}}
\tabcolsep=2mm
\begin{center}
\begin{tabular}{ccccccc}
\hline
\noalign{\smallskip}
player & \multicolumn{6}{c}{dealer upcard} \\
\noalign{\smallskip}
total & @$1$@ & $2$ & $3$ & $4$ & $5$ & $6$ \\
\noalign{\smallskip}
\hline
\noalign{\smallskip}
hard 11--13 &  S & S & S & S & S & S\\
\noalign{\smallskip}
hard 10&  H & S & S & S & H & H \\
\noalign{\smallskip}   
hard 9 &  H & H & S & S & H & H \\
\noalign{\smallskip}
hard 8 &  H & H & H & H & H & H \\
\noalign{\smallskip}
hard 7 &  H & DH & DH & DH & DH & DH \\
\noalign{\smallskip}
hard 6 &  H & H & DH & DH & H & H \\
\noalign{\smallskip}
hard 5 &  H & H & H & H & H & H \\
\noalign{\smallskip}
\hline
\noalign{\smallskip}       
soft 12--13&  S & S & S & S & S & S \\
\noalign{\smallskip} 
soft 9--11 & H & H & H & DH & H & H \\ 
\noalign{\smallskip}
\hline
\noalign{\smallskip}   
$(6,6)$ &  S & S & S & S & S & S \\
\noalign{\smallskip} 
$(5,5)$ &  H & Spl & Spl & Spl & Spl & Spl \\
\noalign{\smallskip} 
$(4,4)$ &  H & H & H & Spl & H & H \\
\noalign{\smallskip}      
$(3,3)$ &  H & H & D & D & H & H \\
\noalign{\smallskip} 
$(2,2)$ &  H & Spl & Spl & Spl & Spl & Spl \\
\noalign{\smallskip} 
$(1,1)$ &  H & Spl & Spl & Spl & Spl & Spl \\
\noalign{\smallskip}
\hline
\end{tabular}
\end{center}
\end{table}
\afterpage{\clearpage}

\section*{Appendix B}

Here we give a more complete account of basic strategy at snackjack.  See Tables~\ref{SJ-BS-1,2,3-deck} for three or fewer decks and Table~\ref{SJ-BS-d-deck} for four or more decks.

\begin{table}[htb]
\caption{\label{SJ-BS-1,2,3-deck}Composition-dependent basic strategy for $d$-deck snackjack, $d=1,2,3$.  Entries containing an en-dash correspond to events that cannot occur, given that the dealer does not have a natural.}
\catcode`@=\active \def@{\hphantom{0}}
\catcode`#=\active \def#{\hphantom{$\;^1$}}
\tabcolsep=.08cm
\begin{center}
\begin{tabular}{ccccccccccccc}
 &  & \multicolumn{3}{c}{single deck} && \multicolumn{3}{c}{double deck} && \multicolumn{3}{c}{triple deck}\\
\noalign{\smallskip}
\hline
\noalign{\smallskip}
player & numbers of & \multicolumn{3}{c}{dealer upcard} && \multicolumn{3}{c}{dealer upcard} && \multicolumn{3}{c}{dealer upcard}\\
\noalign{\smallskip}
total & 1s, 2s, 3s & 1 & 2 & 3 &@& 1 & 2 & 3 &@& 1 & 2 & 3 \\
\noalign{\smallskip}
\hline
\noalign{\smallskip}
@hard 7@& $(0,2,1)$ & @S@& -- & S  && S & S & S && S & S & S \\
        & $(1,0,2)$ & #S#  & #S#  & #S#  && #S# & #S# & #S# && #S# & #S# & #S# \\
        & $(1,3,0)$ & -- & -- & -- && S & S & S && S & S & S \\
        & $(2,1,1)$ & -- & S  & S  && S & S & S && S & S & S \\
        & $(3,2,0)$ & -- & -- & -- && S & S & S && S & S & S \\
        & $(4,0,1)$ & -- & -- & -- && -- & S & S && S & S & S \\
        & $(5,1,0)$ & -- & -- & -- && -- & -- & -- && S & S & S \\
        & $(7,0,0)$ & -- & -- & -- && -- & -- & -- && -- & -- & -- \\
\noalign{\smallskip}
\hline
\noalign{\smallskip}
@hard 6@& $(0,3,0)$ & -- & -- & -- && S & S & S && S & S & S \\
        & $(1,1,1)$ & S  & S  & S  && S & S & S && S & S & S \\
        & $(2,2,0)$ & -- & -- & S  && S & S & S && S & S & S \\
        & $(3,0,1)$ & -- & -- & -- && S & S & S && S & S & S \\
        & $(4,1,0)$ & -- & -- & -- && -- & S & S && S & S & S \\
        & $(6,0,0)$ & -- & -- & -- && -- & -- & -- && -- & S & S \\
\noalign{\smallskip}
\hline
\noalign{\smallskip}   
@hard 5@& $(0,1,1)$ & S  & #H$\,^1$ & H  && H & H & H && H & H & H \\
        & $(1,2,0)$ & -- & -- & H  && S & S & H && H & S & H \\
        & $(2,0,1)$ & -- & S  & H  && H & S & H && H & S & H \\
        & $(3,1,0)$ & -- & -- & -- && H & S & H && H & S & H \\
        & $(5,0,0)$ & -- & -- & -- && -- & -- & -- && H & S & H \\
\noalign{\smallskip}
\hline
\noalign{\smallskip}       
@soft 7@& $(1,0,1)$ & #S$\,^2$  & S  & S  && S & S & S && S & S & S \\
        & $(2,1,0)$ & -- & #S$\,^3$  & S  && S & S & S && S & S & S \\
        & $(4,0,0)$ & -- & -- & -- && -- & S & S && S & S & S \\
\noalign{\smallskip}
\hline
\noalign{\smallskip} 
@soft 6@& $(1,1,0)$ & #H$\,^4$  & D  & D  && #H$\,^1$ & D & D && H & D & D \\
        & $(3,0,0)$ & -- & -- & -- && S & S & S && S & H & S \\    
\noalign{\smallskip}
\hline
\noalign{\smallskip}         
@pair@ & $(0,0,2)$ & S  & Spl & Spl && S & Spl & Spl && S & Spl & Spl \\
        & $(0,2,0)$ & #D$\,^5$  & --  & H   && D & D & D && D & D & D \\
        & $(2,0,0)$ & -- & Spl & Spl && Spl & Spl & Spl && Spl & Spl & Spl \\
\noalign{\smallskip}
\hline
\noalign{\smallskip}
\multicolumn{11}{l}{$^1$Or S.  $^2$Or D.  $^3$Or H. $^4$Or S or D.  $^5$Or Spl.}\\
\multicolumn{5}{l}{}\\
\end{tabular}
\end{center}
\end{table}
\afterpage{\clearpage}

\begin{table}[htb]
\caption{\label{SJ-BS-d-deck}Composition-dependent basic strategy for $d$-deck snackjack, for all integers $d\ge4$.}
\catcode`@=\active \def@{\hphantom{0}}
\catcode`#=\active \def#{\hphantom{$\;^1$}}
\tabcolsep=.08cm
\begin{center}
\begin{tabular}{ccccc}
\hline
\noalign{\smallskip}
player & numbers of & \multicolumn{3}{c}{dealer upcard}\\
\noalign{\smallskip}
total & 1s, 2s, 3s & $1$ & $2$ & $3$ \\
\noalign{\smallskip}
\hline
\noalign{\smallskip}
@hard 7@& $(0,2,1)$ & S & S & S \\
        & $(1,0,2)$ & S & S & S \\
        & $(1,3,0)$ & S & S & S \\
        & $(2,1,1)$ & S & S & S \\
        & $(3,2,0)$ & S & S & S \\
        & $(4,0,1)$ & S & S & S \\
        & $(5,1,0)$ & S & S & S \\
        & $(7,0,0)$ & S & S & S \\
\noalign{\smallskip}
\hline
\noalign{\smallskip}
@hard 6@& $(0,3,0)$ & S & S & S \\
        & $(1,1,1)$ & S & S & S \\
        & $(2,2,0)$ & S & S & S \\
        & $(3,0,1)$ & S & S & S \\
        & $(4,1,0)$ & S & S & S \\
        & $(6,0,0)$ & S & S & S \\
\noalign{\smallskip}
\hline
\noalign{\smallskip}   
@hard 5@& $(0,1,1)$ & H & H & H \\
        & $(1,2,0)$ & H & #H$\,^1$ & H \\
        & $(2,0,1)$ & H & H & H \\
        & $(3,1,0)$ & H & #H$\,^2$ & H \\
        & $(5,0,0)$ & H & #H$\,^3$ & H \\
\noalign{\smallskip}
\hline
\noalign{\smallskip}       
@soft 7@& $(1,0,1)$ & S & S & S \\
        & $(2,1,0)$ & S & S & S \\
        & $(4,0,0)$ & S & S & S \\
\noalign{\smallskip}
\hline
\noalign{\smallskip} 
@soft 6@& $(1,1,0)$ & H & D & D \\
        & $(3,0,0)$ & H & H & H \\    
\noalign{\smallskip}
\hline
\noalign{\smallskip}
@pairs@ & $(0,0,2)$ & #Spl$\,^4$ & Spl & #Spl# \\
        & $(0,2,0)$ & D & D & D \\
        & $(2,0,0)$ & Spl & Spl & Spl \\
\noalign{\smallskip}
\hline
\noalign{\smallskip}
\multicolumn{5}{l}{$^1$S if $4\le d\le6$. $^2$S if $4\le d\le7$.}\\
\multicolumn{5}{l}{$^3$S if $4\le d\le9$. $^4$S if $4\le d\le8$.}\\
\end{tabular}
\end{center}
\end{table}
\afterpage{\clearpage}

\section*{Appendix C}

The polynomial $P$ in \eqref{E(n_1,n_2,n_3)} of degree 8 with 147 terms is given by\bigskip

\noindent $P(n_1,n_2,n_3):=10080 n_1 - 26136 n_1^2 + 26264 n_1^3 - 13538 n_1^4 + 3920 n_1^5 - 
 644 n_1^6 + 56 n_1^7 - 2 n_1^8 - 19248 n_1 n_2 + 26260 n_1^2 n_2 - 
 18802 n_1^3 n_2 + 7494 n_1^4 n_2 - 1638 n_1^5 n_2 + 182 n_1^6 n_2 - 
 8 n_1^7 n_2 + 28360 n_1 n_2^2 - 20516 n_1^2 n_2^2 + 7904 n_1^3 n_2^2 - 
 1874 n_1^4 n_2^2 + 252 n_1^5 n_2^2 - 14 n_1^6 n_2^2 - 21972 n_1 n_2^3 + 
 11254 n_1^2 n_2^3 - 2544 n_1^3 n_2^3 + 284 n_1^4 n_2^3 - 14 n_1^5 n_2^3 + 
 8744 n_1 n_2^4 - 3324 n_1^2 n_2^4 + 484 n_1^3 n_2^4 - 24 n_1^4 n_2^4 - 
 1828 n_1 n_2^5 + 462 n_1^2 n_2^5 - 34 n_1^3 n_2^5 + 192 n_1 n_2^6 - 
 24 n_1^2 n_2^6 - 8 n_1 n_2^7 + 5112 n_1 n_3 - 2980 n_1^2 n_3 - 730 n_1^3 n_3 + 
 1077 n_1^4 n_3 - 371 n_1^5 n_3 + 55 n_1^6 n_3 - 3 n_1^7 n_3 - 3528 n_2 n_3 - 
 1380 n_1 n_2 n_3 + 5484 n_1^2 n_2 n_3 - 2041 n_1^3 n_2 n_3 + 14 n_1^4 n_2 n_3 + 
 79 n_1^5 n_2 n_3 - 8 n_1^6 n_2 n_3 + 6804 n_2^2 n_3 - 9710 n_1 n_2^2 n_3 + 
 807 n_1^2 n_2^2 n_3 + 916 n_1^3 n_2^2 n_3 - 174 n_1^4 n_2^2 n_3 + 
 5 n_1^5 n_2^2 n_3 - 4792 n_2^3 n_3 + 8961 n_1 n_2^3 n_3 - 
 2544 n_1^2 n_2^3 n_3 + 148 n_1^3 n_2^3 n_3 + 8 n_1^4 n_2^3 n_3 + 
 1886 n_2^4 n_3 - 3053 n_1 n_2^4 n_3 + 717 n_1^2 n_2^4 n_3 - 
 43 n_1^3 n_2^4 n_3 - 414 n_2^5 n_3 + 455 n_1 n_2^5 n_3 - 56 n_1^2 n_2^5 n_3 + 
 46 n_2^6 n_3 - 25 n_1 n_2^6 n_3 - 2 n_2^7 n_3 - 14600 n_1 n_3^2 + 
 13906 n_1^2 n_3^2 - 5823 n_1^3 n_3^2 + 1388 n_1^4 n_3^2 - 171 n_1^5 n_3^2 + 
 8 n_1^6 n_3^2 + 3864 n_2 n_3^2 + 12446 n_1 n_2 n_3^2 - 13956 n_1^2 n_2 n_3^2 +\linebreak
 4934 n_1^3 n_2 n_3^2 - 712 n_1^4 n_2 n_3^2 + 36 n_1^5 n_2 n_3^2 - 
 5566 n_2^2 n_3^2 - 1395 n_1 n_2^2 n_3^2 + 3820 n_1^2 n_2^2 n_3^2 - 
 1212 n_1^3 n_2^2 n_3^2 + 100 n_1^4 n_2^2 n_3^2 + 2106 n_2^3 n_3^2 - 
 928 n_1 n_2^3 n_3^2 - 144 n_1^2 n_2^3 n_3^2 + 62 n_1^3 n_2^3 n_3^2 - 
 466 n_2^4 n_3^2 + 289 n_1 n_2^4 n_3^2 - 28 n_1^2 n_2^4 n_3^2 + 
 66 n_2^5 n_3^2 - 24 n_1 n_2^5 n_3^2 - 4 n_2^6 n_3^2 + 10548 n_1 n_3^3 - 
 7393 n_1^2 n_3^3 + 1978 n_1^3 n_3^3 - 278 n_1^4 n_3^3 + 17 n_1^5 n_3^3 - 
 2454 n_2 n_3^3 - 7095 n_1 n_2 n_3^3 +\linebreak
 5424 n_1^2 n_2 n_3^3 - 1198 n_1^3 n_2 n_3^3 + 84 n_1^4 n_2 n_3^3 + 3106 n_2^2 n_3^3 + 
 880 n_1 n_2^2 n_3^3 - 996 n_1^2 n_2^2 n_3^3 + 142 n_1^3 n_2^2 n_3^3 - 
 706 n_2^3 n_3^3 + 24 n_1 n_2^3 n_3^3 + 42 n_1^2 n_2^3 n_3^3 + 56 n_2^4 n_3^3 - 
 7 n_1 n_2^4 n_3^3 - 2 n_2^5 n_3^3 - 3923 n_1 n_3^4 + 2018 n_1^2 n_3^4 - 
 318 n_1^3 n_3^4 + 18 n_1^4 n_3^4 + 786 n_2 n_3^4 + 1860 n_1 n_2 n_3^4 - 
 890 n_1^2 n_2 n_3^4 + 92 n_1^3 n_2 n_3^4 - 936 n_2^2 n_3^4 - 
 112 n_1 n_2^2 n_3^4 + 74 n_1^2 n_2^2 n_3^4 + 156 n_2^3 n_3^4 - 
 4 n_1 n_2^3 n_3^4 - 6 n_2^4 n_3^4 + 777 n_1 n_3^5 - 273 n_1^2 n_3^5 + 
 21 n_1^3 n_3^5 - 114 n_2 n_3^5 - 239 n_1 n_2 n_3^5 + 54 n_1^2 n_2 n_3^5 + 
 126 n_2^2 n_3^5 + 5 n_1 n_2^2 n_3^5 - 12 n_2^3 n_3^5 - 77 n_1 n_3^6 + 
 14 n_1^2 n_3^6 + 6 n_2 n_3^6 + 12 n_1 n_2 n_3^6 - 6 n_2^2 n_3^6 + 3 n_1 n_3^7$.

\section*{Appendix D}
Here we evaluate the snackjack expectation differences when departing from basic strategy with $3,3$ and $2,3$.  We obtain\par\medskip

\noindent$E_{\std,(n_1,n_2,n_3)}((0,0,2),1)-E_{\spl,(n_1,n_2,n_3)}((0,0,2),1)=(12 n_1 - 22 n_1^2 + 12 n_1^3 - 2 n_1^4 - 12 n_2 + 9 n_1^2 n_2 - 3 n_1^3 n_2 + 
 22 n_2^2 - 9 n_1 n_2^2 - 12 n_2^3 + 3 n_1 n_2^3 + 2 n_2^4 - 11 n_1 n_3 + 
 13 n_1^2 n_3 - 4 n_1^3 n_3 + 10 n_2 n_3 - n_1 n_2 n_3 - n_1^2 n_2 n_3 - 
 14 n_2^2 n_3 + 3 n_1 n_2^2 n_3 + 4 n_2^3 n_3 + 4 n_1 n_3^2 - n_1^2 n_3^2 - 
 2 n_2 n_3^2 - n_1 n_2 n_3^2 + 2 n_2^2 n_3^2 - n_1 n_3^3)/[(n_1 + n_2) (n_1 + n_2 + n_3 - 1)_3]$,
\par\medskip
\noindent$E_{\std,(n_1,n_2,n_3)}((0,0,2),2)-E_{\spl,(n_1,n_2,n_3)}((0,0,2),2)=(12 n_1 - 22 n_1^2 + 12 n_1^3 - 2 n_1^4 - 12 n_2 - n_1 n_2 + 11 n_1^2 n_2 - 
 4 n_1^3 n_2 + 22 n_2^2 - 10 n_1 n_2^2 + n_1^2 n_2^2 - 12 n_2^3 + 3 n_1 n_2^3 + 
 2 n_2^4 + 6 n_3 - 21 n_1 n_3 + 23 n_1^2 n_3 - 6 n_1^3 n_3 + 11 n_2 n_3 - 
 2 n_1 n_2 n_3 - 3 n_1^2 n_2 n_3 - 15 n_2^2 n_3 + 4 n_1 n_2^2 n_3 + 4 n_2^3 n_3 - 
 11 n_3^2 + 11 n_1 n_3^2 - 5 n_1^2 n_3^2 - 3 n_2 n_3^2 + n_1 n_2 n_3^2 + 
 3 n_2^2 n_3^2 + 6 n_3^3 - 2 n_1 n_3^3 - n_3^4)/(n_1 + n_2 + n_3)_4$, 
\par\medskip
\noindent$E_{\std,(n_1,n_2,n_3)}((0,0,2),3)-E_{\spl,(n_1,n_2,n_3)}((0,0,2),3)=(2 n_2 + n_1 n_2 - 2 n_1^2 n_2 - 3 n_2^2 + n_1 n_2^2 + n_2^3 + 4 n_1 n_3 - 
 2 n_1^2 n_3 - 3 n_2 n_3 - n_1 n_2 n_3 + 2 n_2^2 n_3 - 2 n_1 n_3^2 + n_2 n_3^2)/[(n_2 + n_3) (n_1 + n_2 + n_3 - 1)_2]$,
\par\medskip
\noindent$E_{\std,(n_1,n_2,n_3)}((0,1,1),1)-E_{\hit,(n_1,n_2,n_3)}((0,1,1),1)=(6 n_1 - 11 n_1^2 + 6 n_1^3 - n_1^4 + 12 n_2 - 27 n_1 n_2 + 19 n_1^2 n_2 - 
 4 n_1^3 n_2 - 22 n_2^2 + 25 n_1 n_2^2 - 7 n_1^2 n_2^2 + 12 n_2^3 - 
 6 n_1 n_2^3 - 2 n_2^4 + 6 n_1^2 n_3 - 2 n_1^3 n_3 - 10 n_2 n_3 + 
 15 n_1 n_2 n_3 - 6 n_1^2 n_2 n_3 + 14 n_2^2 n_3 - 8 n_1 n_2^2 n_3 - 4 n_2^3 n_3 - 
 6 n_1 n_3^2 + 2 n_2 n_3^2 - n_1 n_2 n_3^2 - 2 n_2^2 n_3^2 + 2 n_1 n_3^3)/[(n_1 + n_2) (n_1 + n_2 + n_3 - 1)_3]$,\par\medskip
\noindent$E_{\std,(n_1,n_2,n_3)}((0,1,1),2)-E_{\hit,(n_1,n_2,n_3)}((0,1,1),2)=(6 n_1 - 11 n_1^2 + 6 n_1^3 - n_1^4 + 12 n_2 - 30 n_1 n_2 + 23 n_1^2 n_2 - 
 5 n_1^3 n_2 - 22 n_2^2 + 26 n_1 n_2^2 - 8 n_1^2 n_2^2 + 12 n_2^3 - 
 6 n_1 n_2^3 - 2 n_2^4 - 12 n_3 - 4 n_1 n_3 + 11 n_1^2 n_3 - 3 n_1^3 n_3 - 
 6 n_2 n_3 + 20 n_1 n_2 n_3 - 9 n_1^2 n_2 n_3 + 14 n_2^2 n_3 - 9 n_1 n_2^2 n_3 - 
 4 n_2^3 n_3 + 22 n_3^2 - 4 n_1 n_3^2 - 2 n_1^2 n_3^2 - 4 n_2 n_3^2 - 
 2 n_1 n_2 n_3^2 - 2 n_2^2 n_3^2 - 12 n_3^3 + 2 n_1 n_3^3 + 2 n_2 n_3^3 + 2 n_3^4)/(n_1+n_2+n_3)_4$, and
\par\medskip
\noindent$E_{\std,(n_1,n_2,n_3)}((0,1,1),3)-E_{\hit,(n_1,n_2,n_3)}((0,1,1),3)=(-2 n_2 + 3 n_1 n_2 - n_1^2 n_2 + 3 n_2^2 - 2 n_1 n_2^2 - n_2^3 + 2 n_1 n_3 - 
 n_1^2 n_3 + 2 n_2 n_3 - 3 n_1 n_2 n_3 - 2 n_2^2 n_3 - n_1 n_3^2)/[(n_2 + n_3) (n_1 + n_2 + n_3 - 1)_2]$.
\par\bigskip

Finally, the snackjack expectations with $1,2$ are\medskip

\noindent$E_{\std,(n_1,n_2,n_3)}((1,1,0),1) = -n_1 (-2 n_2 + n_1 n_2 + n_2^2 + n_3 + 2 n_2 n_3 - 
n_3^2)/[(n_1 + n_2) (n_1 + n_2 + n_3 - 1)_2]$,
\par\medskip
\noindent$E_{\hit,(n_1,n_2,n_3)}((1,1,0),1) = (24 n_1 - 50 n_1^2 + 35 n_1^3 - 10 n_1^4 + n_1^5 + 24 n_2 - 
     52 n_1 n_2 + 53 n_1^2 n_2 - 22 n_1^3 n_2 + 3 n_1^4 n_2 - 50 n_2^2 + 
     53 n_1 n_2^2 - 24 n_1^2 n_2^2 + 4 n_1^3 n_2^2 + 35 n_2^3 - 22 n_1 n_2^3 + 
     4 n_1^2 n_2^3 - 10 n_2^4 + 3 n_1 n_2^4 + n_2^5 - 34 n_1 n_3 + 
     46 n_1^2 n_3 - 21 n_1^3 n_3 + 3 n_1^4 n_3 - 2 n_2 n_3 + 6 n_1 n_2 n_3 - 
     13 n_1^2 n_2 n_3 + 4 n_1^3 n_2 n_3 + 6 n_2^2 n_3 - n_1 n_2^2 n_3 + 
     n_1^2 n_2^2 n_3 - 5 n_2^3 n_3 + n_1 n_2^3 n_3 + n_2^4 n_3 + 23 n_1 n_3^2 - 
     15 n_1^2 n_3^2 + 3 n_1^3 n_3^2 - 5 n_2 n_3^2 + 6 n_1 n_2 n_3^2 + 
     6 n_2^2 n_3^2 - 3 n_1 n_2^2 n_3^2 - n_2^3 n_3^2 - 8 n_1 n_3^3 + 
     2 n_1^2 n_3^3 + n_2 n_3^3 - n_1 n_2 n_3^3 - n_2^2 n_3^3 + 
     n_1 n_3^4)/[(n_1 + n_2) (n_1 + n_2 + n_3 - 1)_4]$,
\par\medskip
\noindent$E_{\dbl,(n_1, n_2, n_3)}((1,1,0),1) = 
  2 (-6 n_1 + 11 n_1^2 - 6 n_1^3 + n_1^4 + 6 n_2 - 3 n_1 n_2 - 2 n_1^2 n_2 + 
      n_1^3 n_2 - 11 n_2^2 + 7 n_1 n_2^2 - n_1^2 n_2^2 + 6 n_2^3 - 2 n_1 n_2^3 -
       n_2^4 + 7 n_1 n_3 - 7 n_1^2 n_3 + 2 n_1^3 n_3 - 5 n_2 n_3 + 4 n_1 n_2 n_3 -
       n_1^2 n_2 n_3 + 7 n_2^2 n_3 - 3 n_1 n_2^2 n_3 - 2 n_2^3 n_3 - 4 n_1 n_3^2 +
       n_1^2 n_3^2 + n_2 n_3^2 - n_2^2 n_3^2 + 
      n_1 n_3^3)/[(n_1 + n_2) (n_1 + n_2 + n_3 - 1)_3]$,
\par\medskip
\noindent$E_{\std,(n_1, n_2, n_3)}((1,1,0),2) = (n_1 n_2 - n_1 n_2^2 + 2 n_3 - n_1 n_3 + 3 n_2 n_3 - n_1 n_2 n_3 - 2 n_2^2 n_3 - 3 n_3^2 + n_1 n_3^2 - n_2 n_3^2 + n_3^3)/[(n_1 + n_2 + n_3)_3]$,
\par\medskip
\noindent$E_{\hit,(n_1, n_2, n_3)}((1,1,0),2) = (24 n_1 - 50 n_1^2 + 35 n_1^3 - 10 n_1^4 + n_1^5 + 24 n_2 - 64 n_1 n_2 + 
    72 n_1^2 n_2 - 30 n_1^3 n_2 + 4 n_1^4 n_2 - 50 n_2^2 + 60 n_1 n_2^2 - 
    33 n_1^2 n_2^2 + 6 n_1^3 n_2^2 + 35 n_2^3 - 23 n_1 n_2^3 + 5 n_1^2 n_2^3 - 
    10 n_2^4 + 3 n_1 n_2^4 + n_2^5 + 24 n_3 - 64 n_1 n_3 + 72 n_1^2 n_3 - 
    30 n_1^3 n_3 + 4 n_1^4 n_3 + 24 n_2 n_3 + 20 n_1 n_2 n_3 - 37 n_1^2 n_2 n_3 + 
    9 n_1^3 n_2 n_3 - 22 n_2^2 n_3 + n_1 n_2^2 n_3 + 5 n_1^2 n_2^2 n_3 + 
    4 n_2^3 n_3 - 50 n_3^2 + 60 n_1 n_3^2 - 33 n_1^2 n_3^2 + 6 n_1^3 n_3^2 - 
    34 n_2 n_3^2 + 4 n_1 n_2 n_3^2 + 5 n_1^2 n_2 n_3^2 + 27 n_2^2 n_3^2 - 
    6 n_1 n_2^2 n_3^2 - 4 n_2^3 n_3^2 + 35 n_3^3 - 23 n_1 n_3^3 + 
    5 n_1^2 n_3^3 + 11 n_2 n_3^3 - n_1 n_2 n_3^3 - 5 n_2^2 n_3^3 - 10 n_3^4 + 
    3 n_1 n_3^4 - n_2 n_3^4 + n_3^5)/(n_1 + n_2 + n_3)_5$,
\par\medskip
\noindent$E_{\dbl,(n_1, n_2, n_3)}((1,1,0),2) = 
  2 (-6 n_1 + 11 n_1^2 - 6 n_1^3 + n_1^4 + 6 n_2 - n_1 n_2 - 5 n_1^2 n_2 + 
      2 n_1^3 n_2 - 11 n_2^2 + 7 n_1 n_2^2 - n_1^2 n_2^2 + 6 n_2^3 - 
      2 n_1 n_2^3 - n_2^4 - 6 n_3 + 13 n_1 n_3 - 12 n_1^2 n_3 + 3 n_1^3 n_3 - 
      9 n_2 n_3 + 4 n_1 n_2 n_3 + n_1^2 n_2 n_3 + 12 n_2^2 n_3 - 4 n_1 n_2^2 n_3 - 
      3 n_2^3 n_3 + 11 n_3^2 - 9 n_1 n_3^2 + 3 n_1^2 n_3^2 + 3 n_2 n_3^2 - 
      n_1 n_2 n_3^2 - 3 n_2^2 n_3^2 - 6 n_3^3 + 2 n_1 n_3^3 + 
      n_3^4)/(n_1 + n_2 + n_3)_4$,
\par\medskip
\noindent$E_{\std,(n_1, n_2, n_3)}((1,1,0),3) = n_2 (1 - n_2 + n_3)/[(n_2 + n_3) (n_1 + n_2 + n_3 - 1)]$,
\par\medskip
\noindent$E_{\hit,(n_1, n_2, n_3)}((1,1,0),3) = (3 n_1 n_2 - 4 n_1^2 n_2 + n_1^3 n_2 - n_1 n_2^2 + n_1^2 n_2^2 + 6 n_1 n_3 - 5 n_1^2 n_3 + n_1^3 n_3 - n_2 n_3 - 9 n_1 n_2 n_3 + 4 n_1^2 n_2 n_3 + 
n_2^2 n_3 + 2 n_1 n_2^2 n_3 - 5 n_1 n_3^2 + 2 n_1^2 n_3^2 + n_2 n_3^2 + 
3 n_1 n_2 n_3^2 - n_2^2 n_3^2 + n_1 n_3^3)/[(n_2 + n_3) (n_1 + n_2 + n_3 - 1)_3]$, and
\par\medskip
\noindent$E_{\dbl,(n_1, n_2, n_3)}((1,1,0),3) = -2 (2 n_2 - n_1^2 n_2 - 3 n_2^2 + n_1 n_2^2 + n_2^3 + 2 n_1 n_3 - n_1^2 n_3 - n_2 n_3 - n_1 n_2 n_3 + n_2^2 n_3 - n_1 n_3^2)/[(n_2 + n_3) (n_1 + n_2 + n_3 - 1)_2]$.


\begin{thebibliography}{00}
\bibitem{Ep67}Epstein, R. A.  \textit{The Theory of Gambling and Statistical Logic}.  Academic Press, New York, 1967.

\bibitem{Ep13}Epstein, R. A.  \textit{The Theory of Gambling and Statistical Logic}, Second Edition.  Academic Press, Waltham, MA, 2013.

\bibitem{Et10}Ethier, S. N.  \textit{The Doctrine of Chances: Probabilistic Aspects of Gambling}.  Springer-Verlag, Berlin--Heidelberg, 2010.

\bibitem{EL18}Ethier, S. N. and Lee, J.  The flashing Brownian ratchet and Parrondo's paradox.  \textit{R. Soc.\ Open Sci.} \textbf{5} (2018) (171685) 1--13.  \url{http://rsos.royalsocietypublishing.org/content/5/1/171685}.

\bibitem{EL05}Ethier, S. N. and Levin, D. A.  On the fundamental theorem of card counting with application to the game of trente et quarante.  \textit{Adv. Appl. Probab.} \textbf{37} (2005) 90--107. 

\bibitem{G76}Griffin, P.  The rate of gain in player expectation for card games characterized by sampling without replacement and an evaluation of card counting systems. In Eadington, W. R. (Ed.), \textit{Gambling and Society: Interdisciplinary Studies on the Subject of Gambling}, pp. 429--442. Charles C. Thomas, Springfield, IL, 1976.

\bibitem{G99}Griffin, P. A.  \textit{The Theory of Blackjack: The Compleat Card Counter's Guide to the Casino Game of 21}, Sixth Edition.  Huntington Press, Las Vegas, 1999.

\bibitem{MBG75}Manson, A. R., Barr, A. J., and Goodnight, J. H.  Optimum zero-memory strategy and exact probabilities for 4-deck blackjack. \textit{Amer.\ Statist.} \textbf{29} (1975) 84--88.  Correction \textbf{29} 175.

\bibitem{M09}Marzuoli, A. Toy models in physics and the reasonable effectiveness of mathematics. \textit{Scientifica Acta} \textbf{3} (1) (2009) 13--24.

\bibitem{S18}Schlesinger, D.  \textit{Blackjack Attack:  Playing the Pros' Way}, Third Edition, revised.  Huntington Press, Las Vegas, 2018.

\bibitem{T00}Thorp, E. O. Does basic strategy have the same expectation for each round?  In Vancura, O., Cornelius, J. A., and Eadington, W. R. (eds.) \textit{Finding the Edge: Mathematical Analysis of Casino Games}, pp.~115--132.  Institute for the Study of Gambling and Commercial Gaming, University of Nevada, Reno, 2000.

\bibitem{TW73}Thorp, E. O. and Walden, W. E. The fundamental theorem of card counting with applications to trente-et-quarante and baccarat.  \textit{Int.\ J. Game Theory} \textbf{2} (1973) 109--119.

\bibitem{We18}Werthamer, N. R. \textit{Risk and Reward: The Science of Casino Blackjack}, Second Edition.  Springer International Publishing AG, Cham, Switzerland, 2018.


\end{thebibliography}
\end{document}